\documentclass[11pt, a4paper]{article}

\skip\footins=8mm plus 1mm

\usepackage{latexsym,amsmath,amsfonts,amssymb,amsthm}

\usepackage{longtable, booktabs}
\usepackage{supertabular}
\usepackage{rotating}
\usepackage{multicol}
\usepackage{multirow}

\usepackage{hyperref}
\usepackage[labelsep=period]{caption}

\newtheorem{theorem}{Theorem}[section]

\newtheorem{lemma}[theorem]{Lemma}

\newtheorem{corollary}[theorem]{Corollary}
\newtheorem{remark}[theorem]{Remark}
\newtheorem{definition}[theorem]{Definition}

\def\proof{\par\noindent{\textbf{Proof.}~}}

\long\def\delete#1{}

\newcommand{\beq}{\begin{equation}}
\newcommand{\eeq}{\end{equation}}

\newcommand{\be}{\begin{equation}}
\newcommand{\ee}{\end{equation}}

\def\la{\langle} \def\ra{\rangle}     \def\id{{\rm id}}
 \def\ord{{\rm ord}} \def\G1{G^\mathcal{C}} \def\magma{{\sc Magma} } \def\row{{\rm row}}
\def\diag{{\rm diag}}      

           \def\PG{{\rm PG}}

\def\BB{\mathcal{B}} \def\CC{\mathcal{C}} \def\DD{\mathcal{D}}   
\def\HH{\mathcal{H}} \def\LL{\mathcal{L}}   \def\PP{\mathcal{P}}

\def\De{\Delta} \def\Ga{\Gamma}  \def\Si{\Sigma} \def\Om{\Omega}  

\def\a{\alpha} \def\b{\beta} \def\d{\delta} \def\g{\gamma}  \def\l{\lambda} \def\s{\sigma} \def\t{\tau}
\def\ve{\varepsilon} \def\vp{\varphi} 

  \def\Sym{{\rm Sym}}        \def\Aut{{\rm Aut}}
  \def\soc{{\rm soc}}        \def\PSL{{\rm PSL}}
  \def\PSU{{\rm PSU}}  \def\PGU{{\rm PGU}}  \def\SL{{\rm SL}}   \def\GL{{\rm GL}}  \def\GU{{\rm GU}}
\def\SU{{\rm SU}}  \def\AGL{{\rm AGL}}        \def\AG{{\rm AG}}   \def\PG{{\rm PG}}
    \def\GammaL{{\rm\Gamma L}}   \def\AGammaL{{\rm A\Gamma L}}   
\def\PGammaU{{\rm P\Gamma U}}        \def\Sp{{\rm Sp}}  
            \def\G{{\rm G}}   \def\R{{\rm R}}
\def\Sz{{\rm Sz}}               \def\Co{{\rm Co}}
\def\Cos{{\rm Cos}}              \def\HS{{\rm HS}}

\def\bzero{\mathbf{0}}  \def\ba{\mathbf{a}}   \def\bb{\mathbf{b}}   \def\bc{\mathbf{c}}     \def\bfe{\mathbf{e}}            
\def\bx{\mathbf{x}}   \def\by{\mathbf{y}}   \def\bz{\mathbf{z}}      
      \def\bw{\mathbf{w}}

\begin{document}

\title{Vertex-imprimitive symmetric graphs with exactly one edge between any two distinct blocks}

\author{\renewcommand{\thefootnote}{\arabic{footnote}}Teng Fang\footnotemark[1] , Xin Gui Fang\footnotemark[2] , Binzhou Xia\footnotemark[3] , Sanming Zhou\footnotemark[4]}

\footnotetext[1]{Center for Applied Mathematics at Tianjin University, Tianjin 300072, P. R. China}

\footnotetext[2]{School of Mathematical Sciences, Peking University, Beijing 100871, P. R. China}

\footnotetext[3]{School of Mathematics and Statistics, University of Western Australia, Crawley 6009, WA, Australia}

\footnotetext[4]{School of Mathematics and Statistics, The University of Melbourne, Parkville, VIC 3010, Australia}

\renewcommand{\thefootnote}{}
\footnotetext{{\em E--mail addresses}: \texttt{tfang@tju.edu.cn} (T. Fang), \texttt{xgfang@math.pku.edu.cn} (X. G. Fang), \texttt{binzhou.xia@uwa.edu.au} (B. Xia), \texttt{sanming@unimelb.edu.au} (S. Zhou).}

\date{}

\openup 0.2\jot
\maketitle

\begin{abstract}
A graph $\Ga$ is called $G$-symmetric if it admits $G$ as a group of automorphisms acting transitively on the set of ordered pairs of adjacent vertices. We give a classification of $G$-symmetric  graphs $\Ga$ with $V(\Ga)$ admitting a nontrivial $G$-invariant partition $\BB$ such that there is exactly one edge of $\Ga$ between any two distinct blocks of $\BB$. This is achieved by giving a classification of $(G, 2)$-point-transitive and $G$-block-transitive designs $\DD$ together with $G$-orbits $\Om$ on the flag set of $\DD$ such that $G_{\s, L}$ is transitive on $L \setminus \{\s\}$ and $L \cap N = \{\s\}$ for distinct $(\s, L), (\s, N) \in \Om$, where $G_{\s, L}$ is the setwise stabilizer of $L$ in the stabilizer $G_{\s}$ of $\s$ in $G$. Along the way we determine all imprimitive blocks of $G_{\s}$ on $V \setminus \{\s\}$ for every $2$-transitive group $G$ on a set $V$, where $\s \in V$.
\smallskip


\end{abstract}

\section{Introduction}
\label{sec:introduction}

Intuitively, a graph is symmetric if all its arcs have the same status in the graph, where an \emph{arc} is an ordered pair of adjacent vertices. Since Tutte's seminal work \cite{Tutte}, symmetric graphs have long been important objects of study in graph theory due to their intrinsic beauty and wide applications (see \cite{Praeger00} for an excellent overview of the area). In this paper we give a classification of those symmetric graphs with an automorphism group acting transitively on the arc set and imprimitively on the vertex set such that there is exactly one edge between any two blocks of the underlying invariant partition.

A finite graph $\Ga$ with vertex set $V(\Ga)$ is called {\em $G$-symmetric} (or {\em $G$-arc-transitive}) if it admits $G$ as a group of automorphisms (that is, $G$ acts on $V(\Ga)$ and preserves the adjacency relation of $\Ga$) such that $G$ is transitive on $V(\Ga)$ and transitive on the set of arcs of $\Ga$. (A graph is \emph{symmetric} if it is $\Aut(\Ga)$-symmetric, where $\Aut(\Ga)$ is the full automorphism group of $\Ga$.) The group $G$ is said to be imprimitive on $V(\Ga)$ if $V(\Ga)$ admits a nontrivial {\em $G$-invariant partition} $\BB = \{B, C, \ldots\}$, that is, $1 < |B| < |V(\Ga)|$ and $B^g := \{\a^g\;|\; \a \in B\} \in \BB$ for any $g \in G$ and $B \in \BB$. In this case $(\Ga, G, \BB)$ is said to be a {\em symmetric triple}. The {\em quotient graph} of $\Ga$ relative to $\BB$, denoted by $\Ga_{\BB}$, is defined to be the graph with vertex set $\BB$ such that $B, C \in \BB$ are adjacent if and only if there exists at least one edge of $\Ga$ with one end-vertex in $B$ and the other in $C$. (As usual we assume that $\Ga_{\BB}$ has at least one edge so that each block of $\BB$ is an independent set of $\Ga$.) For adjacent $B, C \in \BB$, define $\Ga[B, C]$ to be the bipartite subgraph of $\Ga$ with bipartition $\{\Ga(C) \cap B, \Ga(B) \cap C\}$, where $\Ga(B)$ is the set of vertices of $\Ga$ with at least one neighbour in $B$. Since $\Ga_{\BB}$ can be easily seen to be $G$-symmetric, this bipartite graph is independent of the choice of adjacent $B, C$ up to isomorphism. Denote by $\Ga_{\BB}(B)$ the neighbourhood of $B$ in $\Ga_{\BB}$, and by $\Ga_{\BB}(\a)$ the set of blocks of $\BB$ containing at least one neighbour of $\a \in V(\Ga)$ in $\Ga$. Denote
\beq
\label{equ:br}
v := |B|,\;\; r := |\Ga_{\BB}(\a)|,\;\; b := |\Ga_{\BB}(B)|,\;\; k := |\Ga(C) \cap B|.
\eeq
Since $\Ga$ is $G$-symmetric and $\BB$ is $G$-invariant, these parameters are independent of the choice of $\a \in V(\Ga)$ and adjacent $B, C \in \BB$.

Various possibilities for $\Ga[B, C]$ can happen. In the ``densest" case where $\Ga[B, C] \cong K_{v,v}$ is a complete bipartite graph, $\Ga$ is uniquely determined by $\Ga_{\BB}$, namely, $\Ga \cong \Ga_{\BB}[K_v]$ is the lexicographic product of $\Ga_{\BB}$ by the complete graph $K_v$. The ``sparsest" case where $\Ga[B, C] \cong K_{2}$ (that is, $k=1$) can also happen; in this case $\Ga$ is called a \emph{spread} of $\Ga_{\BB}$ in \cite{GL}, where it was shown that spreads play a significant role in the study of edge-primitive graphs. Spreads have also arisen from some other classes of symmetric graphs (see \cite{Li-Praeger-Zhou98, Zhou98, Zhou-EJC}), and a study of them was undertaken in \cite[Section 4]{Zhou-flag}. Spreads of cycles and complete graphs with $r=1$ were briefly discussed in \cite[Section 4]{Gardiner-Praeger95}, where Gardiner and Praeger remarked that when $k=1$ and $\Ga_{\BB}$ is a complete graph ``it is not at all clear what one can say about $\Ga$ in general".

In response to the remark above, in this paper we give a classification of all spreads of complete graphs.

\begin{theorem}
\label{THM:MAIN}
All symmetric triples $(\Ga, G, \BB)$ with $G\leq \Aut(\Ga)$ such that there is exactly one edge of $\Ga$ between any two distinct blocks of $\BB$ are classified in this paper and will be described in Sections \ref{sec:almost simple case} and \ref{sec:affine case}.
\end{theorem}

Several interesting families of symmetric triples $(\Ga, G, \BB)$ arise from this classification. In particular, we obtain four infinite families of connected symmetric spreads of complete graphs (see Lemmas \ref{lem:Sz(q)}, \ref{lem:R(q)}, \ref{lem:R(q)-1} and \ref{lem:con-G-flag-G0}). Such graphs are mutually non-isomorphic as they have different orders or valencies.

As shown in \cite[Theorem 4.2]{Gardiner-Praeger95}, in the degenerate case where $r=1$, we have $\Ga \cong {v+1 \choose 2} \cdot K_2$ (the graph of ${v+1 \choose 2}$ independent edges) and $G$ can be any 2-transitive group of degree $v+1$. So we will only consider the general case where $r > 1$.

A major tool in the proof of Theorem \ref{THM:MAIN} is the ``flag graph construction" \cite{Zhou-flag}, which implies that the classification of symmetric spreads of complete graphs is equivalent to that of $(G,2)$-point-transitive and $G$-block-transitive designs together with certain $G$-orbits on their flag sets. We give the latter classification in the following theorem, but postpone related definitions and results on flag graphs to Section \ref{subsec:flag}. (We use $\soc(G)$ to denote the socle of a group $G$, that is, the product of its minimal normal subgroups. We use $G_\bzero$ to denote the stabilizer of the zero vector $\bzero$ when $G$ acts on a vector space $V$, and $\bfe_1 = (1, 0, \ldots, 0)$ the vector of $V$ with $0$ at every coordinate except the first one.)

\begin{theorem}
\label{THM:MAIN RESULT}
Let $\DD$ be a $(G, 2)$-point-transitive and $G$-block-transitive $2$-$(|V|,r+1,\l)\;(r>1)$ design with point set $V$, where $G\leq \Sym(V)$. Then there exists at most one 1-feasible $G$-orbit on the set of flags of $\DD$. Moreover, all possibilities for $(\DD, G)$ such that such a 1-feasible $G$-orbit $\Om$ exists, and the unique $G$-flag graph $\Ga(\DD,\Om,\Psi)$ associated with each $(\DD, G)$ together with its connectedness, are given in Tables \ref{tab0}-\ref{tab0a}, where $\Psi$ is the set of ordered pairs of flags $((\s, L), (\t, N)) \in \Om \times \Om$ such that $\s \ne \t$ and $\s, \t \in L \cap N$. Furthermore, $\Om$ is explicitly given for each $(\DD, G)$ in all cases except (h)-(j) in Table \ref{tab0a}.
\end{theorem}

\begin{table}
\begin{center}
  \begin{tabular}{l|l|l|l|l}
\hline
     & $G$  & $\DD$ & $\Ga(\DD,\Om,\Psi)$  & Details \\ \hline
(a) & $A_7$ & $\PG(3,2)$ & $35 \cdot K_3$ & \textsection\ref{sec:almost simple case} \\ \hline
(b) &  $\soc(G)=\PSL(d,q)$ & $\PG(d-1, q)$ & $\frac{(q^d-1)(q^{d-1}-1)}{(q^2-1)(q-1)} \cdot K_{q+1}$ & \textsection\ref{sec:almost simple case} \\
& $d \ge 3$ & & &\\ \hline
(c) & $G \le \PGammaU(3,q)$ & $2$-$(q^3+1,q+1,1)$ &  $(q^4-q^3+q^2) \cdot K_{q+1}$ & L\ref{lem:PSU(3,q), nontrivial imprimitive blocks} \\
     & $\soc(G)=\PSU(3,q)$ & &  & \\
     & $q \geq 3$ & &  & \\ \hline
(d) & $\soc(G)=\Sz(q)$ & $2$-$(q^2+1, q+1, q+1)$ & C, ord $=q(q^2+1)$ & L\ref{lem:Sz(q)} \\
     & $q=2^{2e+1}>2$ &  & and val $=q$  & \\ \hline
(e) & $\soc(G)=\R(q)$ & $2$-$(q^3+1, q+1, 1)$ & $(q^4-q^3+q^2) \cdot K_{q+1}$ & \textsection\ref{subsec:Ree} \\
     & $q=3^{2e+1}\geq3$ &  &  & \\ \hline
     & $\soc(G)=\R(q)$ & $2$-$(q^3+1, q+1, q+1)$ & C, ord $=q^2(q^3 + 1)$ and & L\ref{lem:R(q)} \\
     & $q=3^{2e+1} \ge 3$ &  & val $=q$ if $q > 3$; three & \\
     & & & components if $q=3$ & \\ \hline
     & $\soc(G)=\R(q)$ & $2$-$(q^3+1, q^2+1, q^2+1)$ & C, ord $=q(q^3 + 1)$ and & L\ref{lem:R(q)-1} \\
     & $q=3^{2e+1} \geq 3$ &  & val $=q^2$ & \\ \hline
     & $\R(3)$ & $2$-$(28, 10, 10)$ & Three components & L\ref{lem:R(q)-2} \\ \hline
\end{tabular}
  \caption{Theorem \ref{THM:MAIN RESULT}: almost simple case. Acronym: L = Lemma, T = Table, C = Connected, D = Disconnected, ord = Order, val = Valency}
\label{tab0}
  \end{center}
\end{table}

\begin{table}
\begin{center}
  \begin{tabular}{l|l|c|l|l}
\hline
     & $G$  & $\DD$ & $\Ga(\DD,\Om,\Psi)$  & Details \\ \hline
(f) & $G \le \AGammaL(1,q)$ & $2$-$(q, p^t, 1)$; $\DD$ has a block $L$ & $\frac{q(q-1)}{p^t (p^t-1)}\cdot K_{p^t}$ & \textsection\ref{subsec:AGammaL(1, q)} \\
     & $q=p^{d}$ & such that $L$ is a subfield of $\mathbb{F}_{q}$  & & \\\hline
     & $G \le \AGammaL(1,q)$ & $2$-$(q, |L|, |L|)$; $\DD$ has a block & C iff $p\equiv-1\;(\bmod\;4)$, & L\ref{lem:con-G-flag-G0} \\
     & $q=p^{d}$ & $L$ with $0, 1 \in L$ and $L \setminus \{0\}$ & $d$ is odd, and $P \le \mathbb{F}_q^\times$ & \\
     & & the union of some cosets of & with index $2$; in this case & \\
     & & a subgroup of $\mathbb{F}_{q}^\times$ & ord $=2q$ and val $=\frac{q-1}{2}$ & \\ \hline
(g) & $G\leq\AGammaL(n,q)$ & $2$-$(q^n, |L|, \l)$, $\l=1$ or $|L|$; & D & L\ref{lem:Sp(n,q)} \\
    & $G_\bzero\unrhd\Sp(n,q)$ & $\DD$ has a block $L\subseteq\la\bfe_1\ra$ & &   \\
    & $V=\mathbb{F}_q^n$, $n \ge 4$ even & with $\{a\in\mathbb{F}_q^\times\;|\;a\bfe_1\in L\}$ & & \\
    & & a subgroup of $\mathbb{F}_q^\times$ & & \\ \hline
    & $G\leq\AGammaL(n,q)$ & As above & D & L\ref{lem:SL(2,q)=Sp(2,q)} \\
    & $G_\bzero\unrhd\SL(n,q)$ &  & & \textsection\ref{ssec:SLnq} \\
    & $V=\mathbb{F}_q^n$, $n \ge 2$ & & & \\ \hline
    & $G\leq\AGammaL(6,q)$ & As above ($n=6$) & D & L\ref{lem:G2q} \\
    & $G_\bzero\unrhd G_{2}(q)$ &  & &  \\
    & $V=\mathbb{F}_q^6$, $q > 2$ even & & & \\ \hline
(h) & $G\leq\AGL(6,3)$ & $2$-$(3^6, r+1, \l)$, $\l = 1$ or $r+1$; & D & T\ref{tab1} \\
    & $G_\bzero\cong\SL(2,13)$ & $\DD$ has a block $L$ such that &  &   \\
    & $V=\mathbb{F}_3^6$ & $L\setminus\{\bzero\}$ is an orbit of some & & \\
    &  &  $H \leq G_\bzero$ on $V$ with $G_{\bzero,\bx}\leq H$ & & \\
    &  & for some $\bx \in L\setminus\{\bzero\}$ & & \\ \hline
(i) & $G\leq \AGL(2,p)$ & $2$-$(p^2, r+1, \l)$, $\l = 1$ or $r+1$; & D & T\ref{tab:p=5}-\ref{tab:p=59} \\
    & $G_{\bf 0}\unrhd\SL(2,3)$ or & $\DD$ has a block $L$ with the &  &   \\
    & $G_\bzero\unrhd\SL(2,5)$ & same property as in (h)  & & \\
    & $p = 5, 7, 11, 19, 23$ &  & & \\
    & $29, 59$, $V = \mathbb{F}_p^2$ &  & & \\ \hline
(j) & $G\leq \AGL(4,3)$ & $2$-$(3^4, r+1, \l)$, $\l = 1$ or $r+1$; & D & T\ref{tab2}-\ref{tab3} \\
    & $G_{\bf 0}\unrhd\SL(2,5)$ or & $\DD$ has a block $L$ with the &  &  \\
    & $G_\bzero\unrhd E$, $V = \mathbb{F}_3^4$ & same property as in (h) & & \\ \hline
\end{tabular}
  \caption{Theorem \ref{THM:MAIN RESULT}: affine case. Acronym: L = Lemma, T = Table, C = Connected, D = Disconnected, ord = Order, val = Valency}
\label{tab0a}
  \end{center}
\end{table}

Theorem \ref{THM:MAIN} follows from Theorem \ref{THM:MAIN RESULT} and Corollary \ref{corol:k=1&complete quotient}, but details of the corresponding flag graphs $\Ga(\DD,\Om,\Psi)$ will be given during the proof of Theorem \ref{THM:MAIN RESULT}. The number of pairs $(\DD, G)$ in each of (h)-(j) above will be computed by \magma but their structures will not be given due to space limit. As a byproduct of the proof of Theorem \ref{THM:MAIN RESULT}, we give (or enumerate in cases (h)-(j)) all imprimitive blocks of $G_{\s}$ on $V \setminus \{\s\}$ for every 2-transitive group $G$ on $V$.

The reader is referred to \cite{Li-Praeger-Zhou98, Zhou98, Zhou-EJC, Zhou-flag} for several studies on vertex-imprimitive symmetric graphs using the geometric approach developed in \cite{Gardiner-Praeger95}. Together with \cite{CZ14, SZ13, Zhou-EJC}, in \cite{FXZ} the authors of the present paper completed the classification of symmetric triples $(\Ga, G, \BB)$ such that $k=v-1 \ge 2$ and $\Ga_{\BB}$ is a complete graph.

\section{Preliminaries}

\subsection{Notation and terminology}
\label{subsec:def}

The reader is referred to \cite{Dixon-Mortimer} and \cite{Beth-Jung-Lenz} for notation and terminology on permutation groups and block designs, respectively. Unless stated otherwise, all designs are assumed to have no repeated blocks, and each block of a design is identified with the set of points incident with it.

Let $G$ be a group acting on a set $\Om$. That is, for any $\a \in \Om$ and $g \in G$ there corresponds a point in $\Om$ denoted by $\a^g$, such that $\a^{1_G} = \a$ and $(\a^g)^h = \a^{gh}$ for any $\a \in \Om$ and $g, h \in G$, where $1_G$ is the identity element of $G$. Let $P_i$ be a point or subset of $\Om$ for $i=1,2,\ldots,n$. Define $(P_1,P_2,\ldots,P_n)^g:=(P_1^g,P_2^g,\ldots,P_n^g)$ for $g\in G$, where $P_i^g:=\{\a^g\;|\; \a\in P_i\}$ if $P_i$ is a subset of $\Om$. Denote $P_i^G:=\{P_i^g\;|\; g\in G\}$. In particular, $\a^G$ is the $G$-orbit on $\Om$ containing $\a$. Define $G_{P_1,P_2,\ldots,P_n}:=\{g\in G\;|\; P_i^g=P_i,\;i=1,\ldots,n\}\leq G$. In particular, if $\a$ is a point and $P$ a subset of $\Om$, then $G_{\a}$ is the stabilizer of $\a$ in $G$, $G_P$ is the setwise stabilizer of $P$ in $G$, and $G_{\a, P}$ is the setwise stabilizer of $P$ in $G_{\a}$.

Let $G$ and $H$ be groups acting on $\Om$ and $\De$, respectively. These two actions are said to be {\em permutation isomorphic} if there exist a bijection $\rho: \Om \rightarrow \De$
and an isomorphism $\eta: G \rightarrow H$ such that $\rho(\a^g) = (\rho(\a))^{\eta(g)}$ for $\a \in \Om$ and $g \in G$. If in addition $G=H$ and $\eta$ is the identity automorphism of $G$, then the two actions are said to be {\em permutation equivalent}. It is immediate from the definition that if $\vp: G\rightarrow \Sym(\Om)$ and $\psi: H\rightarrow \Sym(\Om)$ are monomorphisms, then the corresponding actions of $G$ and $H$ on $\Om$ are permutation isomorphic if and only if $\vp(G)$ and $\psi(H)$ are conjugate in $\Sym(\Om)$. Let $\Ga$ and $\Si$ be $G$-symmetric graphs. If there exists a graph isomorphism $\rho: V(\Ga) \rightarrow V(\Si)$ such that the actions of $G$ on $V(\Ga)$ and $V(\Si)$ are permutation equivalent with respect to $\rho$, then $\Ga$ and $\Si$ are said to be {\em $G$-isomorphic} with respect to the {\em $G$-isomorphism} $\rho$, and we denote this fact by $\Ga \cong_G \Si$.

\subsection{Flag graphs}
\label{subsec:flag}

Let $(\Ga, G, \BB)$ be a symmetric triple. As in \cite{Gardiner-Praeger95}, define $\DD(B) := (B, \Ga_{\BB}(B))$ to be the incidence structure with ``point set'' $B$ and ``block set'' $\Ga_{\BB}(B)$ such that $\a \in B$ and $C \in \Ga_{\BB}(B)$ are incident if and only if $C \in \Ga_{\BB}(\a)$. It can be verified \cite{Gardiner-Praeger95} that $\DD(B)$ is a $1$-$(v, k, r)$ design with $b$ blocks, and is independent of $B$ up to isomorphism. Denote by $\DD^{*}(B) := (\Ga_{\BB}(B), B)$ the dual $1$-design of $\DD(B)$. We may identify the ``blocks'' $\a \in B$ of $\DD^{*}(B)$ with the subsets $\Ga_{\BB}(\a)$ of the ``point set'' $\Ga_{\BB}(B)$ of $\DD^{*}(B)$, and we call two such ``blocks'' $\Ga_{\BB}(\b)$, $\Ga_{\BB}(\g)$ {\em repeated} if $\b, \g \in B$ are distinct but $\Ga_{\BB}(\b) = \Ga_{\BB}(\g)$. For $\a \in V(\Ga)$, let $B(\a)$ denote the unique block of $\BB$ containing $\a$ and set $L(\a) := \{B(\a)\} \cup \Ga_{\BB}(\a)$. Let $\LL$ be the set of all $L(\a)$, $\a \in V(\Ga)$, with repeated ones identified. One can see that the action of $G$ on $\BB$ induces an action of $G$ on $\LL$ defined by $L(\a)^g := L(\a^g)$, for $\a \in V(\Ga)$ and $g \in G$. Define
$$
\DD(\Ga, \BB) := (\BB, \LL)
$$
to be the incidence structure with incidence relation the set-theoretic inclusion. Then $\DD(\Ga, \BB)$ is a 1-design with block size $r+1$ which admits $G$ as a point- and block-transitive group of automorphisms (\cite[Lemma 3.1]{Zhou-flag}). In the case when $\DD^{*}(B)$ has no repeated blocks (which occurs particularly when $k=1$ by \cite[Lemma 4.1(a)]{Zhou-flag}), $\Ga$ can be reconstructed from $\DD(\Ga, \BB)$ by the following construction (see \cite{Zhou-EJC, Zhou-flag} and \cite[Construction 5.2]{IPZ}).

\begin{definition}
\label{def:strongly feasible}
{\em (\cite[Definitions 2.1 and 4.2]{Zhou-flag}) Let $\DD$ be a $G$-point-transitive and $G$-block-transitive $1$-design with block size at least $2$. Let $\s$ be a point of $\DD$, $\Om$ a $G$-orbit on the set of flags of $\DD$ and $\Om(\s)$ the set of flags of $\Om$ with point entry $\s$. The set $\Om$ is said to be {\em feasible} if
\begin{itemize}
\item[(a)] $|\Om(\s)|\geq2$; and
\item[(b)] for some (and hence all) flag $(\s, L) \in \Om$, $G_{\s, L}$ is transitive on $L \setminus \{\s\}$.
\end{itemize}
If $\Om$ also satisfies
\begin{itemize}
\item[(c)] $L \cap N = \{\s\}$, for distinct $(\s, L), (\s, N) \in \Om(\s)$,
\end{itemize}
then $\Om$ is said to be {\em 1-feasible}.

Given a feasible $G$-orbit $\Om$ on the set of flags of $\DD$, denote
\begin{equation}
\label{equ:Xi}
\Xi(\DD, \Om):=\{((\s,L), (\t,N)) \in \Om\times\Om\;|\;\s \neq \t\;\text{and\;} \s, \t \in L \cap N\}.
\end{equation}
If $\Psi$ is a {\em self-paired} $G$-orbit (that is, $((\s, L), (\t, N)) \in \Psi$ implies $((\t, N)$, $(\s, L)) \in \Psi$) on $\Xi(\DD, \Om)$, then the {\em $G$-flag graph} of $\DD$ with respect to $(\Om, \Psi)$, denoted by $\Ga(\DD, \Om, \Psi)$, is defined to be the graph with vertex set $\Om$ and arc set $\Psi$.
}
\end{definition}

The following was proved in \cite[Theorem 1.1]{Zhou-flag}: If $\DD^{*}(B)$ contains no repeated blocks, then $\Ga$ is $G$-isomorphic to a $G$-flag graph of $\DD(\Ga, \BB)$ with respect to some $(\Om, \Psi)$. Conversely, any $G$-flag graph $\Ga(\DD, \Om, \Psi)$ is a $G$-symmetric graph admitting
\begin{equation*}
\BB(\Om) := \{\Om(\s)\;|\; \s\; \text{is a point of} \;\DD\}
\end{equation*}
as a $G$-invariant partition such that the corresponding $\DD^{*}(\Om(\s))$ contains no repeated blocks. The case where $k =1$ occurs if and only if $\Om$ is 1-feasible, and in this case $G$ is faithful on $V(\Ga)$ if and only if it is faithful on the point set of $\DD$ (see \cite[Theorem 4.3]{Zhou-flag} and the remark below it). In the case where $k =1$ and $\Ga_{\BB}$ is a complete graph, we have $\Ga_{\BB} \cong K_{vr+1}$ as $\Ga_{\BB}$ has valency $vr$. Since $\Ga_{\BB}$ is $G$-symmetric, this occurs precisely when $G$ is $2$-transitive on $\BB$. Hence in this case $\DD(\Ga, \BB)$ is a $(G, 2)$-point-transitive and $G$-block-transitive $2$-$(vr+1, r+1, \l)$ design for some integer $\l \geq 1$. Conversely, if $\DD$ is a $(G, 2)$-point-transitive and $G$-block-transitive $2$-$(vr+1, r+1, \l)$ design, then for any $G$-flag graph $\Ga = \Ga(\DD, \Om, \Psi)$ of $\DD$, we have $\Ga_{\BB(\Om)} \cong K_{vr+1}$. Thus \cite[Theorem 4.3]{Zhou-flag} has the following consequence.

\begin{corollary}
\label{corol:k=1&complete quotient}
(\cite[Corollary 4.4]{Zhou-flag})
Let $v \geq 2$ and $r \geq 1$ be integers, and let $G$ be a group. Then the following statements are equivalent:
\begin{itemize}
\item[\rm (a)] $\Ga$ is a $G$-symmetric graph of valency $r$ admitting a nontrivial $G$-invariant
partition $\BB$ of block size $v$ such that $\Ga[B,C]\cong K_2$ for any two distinct blocks $B,C$ of $\BB\;($so $\Ga_{\BB} \cong K_{vr+1})$;

\item[\rm (b)] $\Ga \cong_G \Ga(\DD,\Om,\Psi)$, for a $(G, 2)$-point-transitive and $G$-block-transitive $2$-$(vr+1, r+1, \l)$ design $\DD$, a 1-feasible $G$-orbit $\Om$ on the set of flags of $\DD$, and a self-paired $G$-orbit $\Psi$ on $\Xi(\DD, \Om)$.
\end{itemize}
\noindent In addition, $G$ is faithful on $V(\Ga)$ if and only if it is faithful on the point set of $\DD$.

Moreover, for any point $\s$ of $\DD$, the set of points of $\DD$
other than $\s$ admits a $G_{\s}$-invariant partition of block
size $r$, namely, $\{L \setminus \{\s\}\;|\; (\s, L) \in \Om\}$. Hence
$\DD$ is not $(G, 3)$-point-transitive when $r \ge 2$. Furthermore,
either $\l = r+1$, or $\l = 1$ and $\Ga \cong (v(vr+1)/(r+1))\cdot K_{r+1}$.
\end{corollary}

In the case where $r=1$, by the $G$-flag graph construction, $V(\Ga)$ can be identified with $V^{(2)}:=\{(i,j)\;|\;0\leq i,j\leq v,i\neq j\}$. Thus $\BB=\{B_0,B_1,B_2,\ldots,B_v\}$ where $B_i=\{(i,j)\;|\;0\leq j\leq v,j\neq i\}$, $\Ga={v+1 \choose 2}\cdot K_2$ with edge set $\{\{(i,j), (j,i)\}\;|\;0\leq i,j\leq v,i\neq j\}$, and $G$ is any $2$-transitive group on $\{0,1,2,\ldots,v\}$ acting on $V^{(2)}$ coordinate-wise (see also \cite[Theorem 4.2]{Gardiner-Praeger95}).

In what follows we assume $r > 1$. Since the graphs in Theorem \ref{THM:MAIN} are precisely those in Corollary \ref{corol:k=1&complete quotient} (a), to prove Theorem \ref{THM:MAIN} it suffices to classify $(\DD,\Om,\Psi)$ in part (b) of this corollary with $G$ faithful on the point set of $\DD$. In the remainder of the paper we give this classification and thus prove Theorem \ref{THM:MAIN RESULT}.

\subsection{Preliminary results}
\label{sec:facts and necessary conditions}

\begin{lemma}
\label{lem:2 lambda}
Let $\DD$ be a $(G, 2)$-point-transitive and $G$-block-transitive $2$-$(|V|,r+1,\l)$ design with point set $V$.
Then there is at most one 1-feasible $G$-orbit on the flag set of $\DD$. Moreover, if such an orbit exists, say $\Om = (\d, L)^G$, then either (a) $G_L$ is transitive on $L$ $($or equivalently $G_L \not \le G_\d)$, $\lambda = 1$, and $\Om$ is the set of all flags of $\DD$; or (b) $G_L$ is not transitive on $L$ $($or equivalently $G_L \le G_\d)$ and $\lambda = r+1$.
\end{lemma}
\proof
Let $\s,\t \in V$ be distinct points. Denote by $L_1, \ldots, L_{\lambda}$ the $\lambda$ blocks of $\DD$ containing both $\s$ and $\t$. Denote
$$
\Om_i := (\s, L_i)^G,\quad i=1, \ldots, \lambda.
$$
Since $G$ is $2$-transitive on $V$, the sets $\Om_1, \ldots, \Om_{\l}$ are all possible $G$-orbits on the flag set of $\DD$, possibly with $\Om_i = \Om_j$ for distinct $i$ and $j$. Suppose that $\Om=(\d, L)^G$ is 1-feasible. Then by (c) in Definition \ref{def:strongly feasible} we have $|V|=vr+1$ for some integer $v$.

First assume that $G_L$ is transitive on $L$. Then $(\eta,L)\in\Om$ for every $\eta\in L$. Hence, by the transitivity of $G$ on the block set of $\DD$, $\Om= \Om_1 = \cdots = \Om_\lambda$ is the flag set of $\DD$ and thus $\l=1$ by (c) in Definition \ref{def:strongly feasible}.

Next assume that $G_L$ is not transitive on $L$. Then $G_L \le G_\d$ by (b) in Definition \ref{def:strongly feasible}. Let $\eta$ be a fixed point of $V$. For each $\pi \in V \setminus \{ \eta \}$, by (c) in Definition \ref{def:strongly feasible} there is only one flag in $\Om(\pi)$ whose block entry contains $\eta$. On the other hand, if $(\t_1, M)$ and $(\t_2, M)$ are flags in $\Om$ for distinct $\t_1$, $\t_2\in V$ and some $M \in L^G$, then $\t_1, \t_2\in M$ and there exists $g \in G$ such that $(\t_1, M) = (\t_2, M)^g$. Thus $g \in G_M$ and $\t_1 = \t_2^g$. Since $\Om$ satisfies (b) in Definition \ref{def:strongly feasible}, $G_M$ is transitive on $M$, but this contradicts our assumption. Hence there are $|\Om(\eta)| + (|V|-1)=v + rv = (r+1)v$ blocks of $\DD$ containing $\eta$. By the relations between parameters of the $2$-design $\DD$, we get $\lambda = r+1$.

Suppose that $\Om_i \ne \Om_j$ and both of them are 1-feasible. Since $\DD$ is $G$-block-transitive, there exists a point $\xi$ of $\DD$ such that $(\xi, L_j) \in \Om_i$. The assumption $\Om_i\neq\Om_j$ implies $\s \neq \xi$ and $G_{L_j} = G_{\xi, L_j} \le G_{\xi}$ (for otherwise $G_{L_j}$ is transitive on $L_j$ and $\Om_i=\Om_j$ is the flag set of $\DD$ by the proof above). Since $\xi \in L_j \setminus \{ \s \}$, $G_{\s, L_j} \le G_{L_j} \le G_{\xi}$ and $|L_j| = r+1 \geq 3$, $G_{\s, L_j}$ cannot be transitive on $L_j \setminus \{ \s \}$, which contradicts the assumption that $\Om_j$ is 1-feasible. Hence there is at most one 1-feasible $G$-orbit on the flag set of $\DD$.
\qed

\begin{lemma}
\label{lem:1 basic}
Let $\DD$ be a $(G, 2)$-point-transitive and $G$-block-transitive $2$-$(|V|,r+1,\l)$ design with point set $V$.
Suppose that there is a 1-feasible $G$-orbit $\Om=(\s, L)^G$ on the flag set of $\DD$. Let $P := L \setminus \{\s\}$ and let $H$ be a transitive subgroup of $G_{\s}$ on $V \setminus \{\s\}$. Then $\Xi(\DD, \Om)$ (see \eqref{equ:Xi}) is a self-paired $G$-orbit on the set of ordered pairs of distinct flags in $\Om$, $P$ is an imprimitive block of $H$ on $V \setminus \{\s\}$, and $P$ is the union of some $H_\t$-orbits (including the $H_\t$-orbit $\{\t\}$ of length 1), where $\t\in P$ is a fixed point.
\end{lemma}

\proof
By (c) in Definition \ref{def:strongly feasible} we know that $\PP = \{M \setminus\{\d\}\;|\; (\d, M) \in \Om\}$ is a $G_{\s}$-invariant partition of $V \setminus \{\s\}$. Since $H \le G_{\s}$ is transitive on $V \setminus \{\s\}$, it follows that $P$ is an imprimitive block of $H$ on $V \setminus \{\s\}$. Since $\t\in P$, $P$ is $H_\t$-invariant and so is the union of some $H_\t$-orbits. Let $((\d, M), (\pi, N))$, $((\xi,\overline{M}), (\eta,\overline{N}))\in \Xi(\DD, \Om)$. Then there exists $g\in G$ such that $\d^g=\xi$ and $\pi^g=\eta$. Hence $\xi,\eta\in M^g\cap N^g$. By (c) in Definition \ref{def:strongly feasible}, we have $\overline{M}=M^g$ and $\overline{N}=N^g$, and thus $\Xi(\DD, \Om)$ is a self-paired $G$-orbit on the set of ordered pairs of distinct flags in $\Om$.
\qed

\smallskip
From now on we use the following abbreviations when $\DD$ and $\Om$ are clear from the context:
$$
\Xi := \Xi(\DD, \Om),\; \Ga(\DD,\Om,\Xi) := \Ga(\DD,\Om,\Xi(\DD, \Om)).
$$

Given a group $G$, a subgroup $T$ of $G$, and an element $g\in G$ with $g \not \in N_{G}(T)$ and $g^2\in T \cap T^g$, define the coset graph $\Cos(G,T,TgT)$ to be the graph with vertex set $[G:T] := \{Tx\;|\; x \in G\}$ and edge set $\{\{Tx,Ty\}\;|\;xy^{-1}\in TgT\}$. It is well known (see e.g. \cite{Praeger97}) that $\Cos(G,T,TgT)$ is a $G$-symmetric graph with $G$ acting on $[G:T]$ by right multiplication, and $\Cos(G,T,TgT)$ is connected if and only if $\la T,g\ra=G$. Conversely, any $G$-symmetric graph $\Ga$ is $G$-isomorphic to $\Cos(G,T,TgT)$ (see e.g. \cite{Praeger97}), where $g$ is an element of $G$ interchanging two adjacent vertices $\a$ and $\b$ of $\Ga$ and $T:=G_\a$, and the required $G$-isomorphism is given by $V(\Ga) \rightarrow [G:T], \g \mapsto Tx$, with $x \in G$ satisfying $\a^x = \g$. Based on this one can prove the following result.

\begin{lemma}
\label{lem:connectedness of G-flag graph}
Let $((\s,L),(\t,N))\in\Xi$ and $T:=G_{\s, L}$. Let $g \in G$ satisfy $(\s, \t)^g = (\t, \s)$ and set $H:=\la T,g\ra$. Then $\rho: \Om \rightarrow [G:T], \g \mapsto Tx$, with $x \in G$ satisfying $(\s, L)^x = \g$, defines a $G$-isomorphism from $\Ga(\DD,\Om, \Xi)$ to  $\Cos(G,T,TgT)$, under which the preimage of the subgraph $\Cos(H,T,TgT)$ of $\Cos(G,T,TgT)$ is the connected component of $\Ga(\DD,\Om, \Xi)$ containing the vertex $(\s,L)$.
\end{lemma}

In the proof of Theorem \ref{THM:MAIN RESULT}, we will use the following results whose proofs are straightforward and hence omitted.

\begin{lemma}
\label{lem:6 num1}
Let $\ell$ be a nonnegative integer and $q>2$ a prime power such that $3 \mid (q+1)$.
\begin{itemize}
\item[{\em (a)}] If $(\ell (q^2-1)/{3} + q) \mid q^3$, then $\ell = 0$ or $3q$;

\item[{\em (b)}] if $(\ell (q^2-1)/3 + 1) \mid q^3$, then $\ell = 0$ or $3$.
\end{itemize}
\end{lemma}

\begin{lemma}
\label{lem:8 num3}
Let $\ell$ and $n$ be positive integers and $q>1$ a prime power. If $(\ell(q-1)+1) \mid q^n$, then $\ell = ({q^{i}-1})/({q-1})$ for some $i = 1, 2, \ldots, n$.
\end{lemma}

\begin{lemma}
\label{lem:9 num4}
Let $\ell\geq0$ be an integer, $n$ a positive integer, and $q$ an odd power of $3$. Then $(\ell(q-1)+({q-1})/{2}+1)$ $\nmid q^n$.
\end{lemma}

\subsection{Methodology}
\label{sec:our method}

Let $\DD$ be a $(G, 2)$-point-transitive and $G$-block-transitive $2$-$(|V|,r+1,\l)\;(r>1)$ design with point set $V$, where $G\leq \Sym(V)$. By Lemma \ref{lem:1 basic}, for a 1-feasible $G$-orbit $\Om = (\s, L)^G$ on the set of flags of $\DD$, $L \setminus \{\s\}$ must be an imprimitive block of $G_\s$ on $V\setminus\{\s\}$. The major task in the proof of Theorem \ref{THM:MAIN RESULT} is to determine such blocks by examining subgroups $H$ of $G_\s$, with the help of the following methods.

\smallskip
(i) Suppose that $H \le G_\s$ is transitive on $V \setminus \{\s\}$. For each nontrivial imprimitive block $P$ of $H$ on $V \setminus \{\s\}$, we will check whether $P$ is also an imprimitive block of $G_\s$ on $V \setminus \{\s\}$. By Lemma \ref{lem:1 basic}, $P$ is the union of some $H_\t$-orbits on $V \setminus \{\s\}$, where $\t \in P$.

\smallskip
(ii) If there is a point $\t \in V \setminus \{\s\}$ such that $H_\t=G_{\s,\t}$, then  by \cite[Theorem 1.5A]{Dixon-Mortimer}, $P:=\t^H$ is an imprimitive block of $G_\s$ on $V \setminus \{\s\}$.

\smallskip
If $P$ is an imprimitive block of $G_\s$ on $V \setminus \{\s\}$ obtained from (i) or (ii), then define
\begin{align*}
\DD:=(V, L^G),\; \text{where \;} L:= P \cup \{\s\},
\end{align*}
to be the incidence structure with point set $V$ and block set $L^G$, with incidence relation the set-theoretic inclusion. By \cite[Proposition 4.6]{Beth-Jung-Lenz}, $\DD$ is a $2$-$(|V|, |P|+1, \l)$ design admitting $G$ as a 2-point-transitive and block-transitive group of automorphisms. Moreover, $\Om:= (\s, L)^G$ is the unique 1-feasible $G$-orbit on the flag set of $\DD$, and $(\s,L)$ is adjacent to $r := |P|$ vertices $(\eta_1, M_1)$, $\ldots$, $(\eta_r, M_r)$ in the $G$-flag graph $\Ga(\DD,\Om,\Xi)$, where $\{\eta_1,\ldots, \eta_{r}\}=P$ and $\s\in M_i\setminus \{\eta_i\}$, $i=1,2,\ldots, r$. The value of $\l$ and the connectedness of $\Ga(\DD,\Om,\Xi)$ will be determined for each $P$.

Since $G$ is $2$-transitive on $V$, it is either almost simple (with socle a nonabelian simple group) or affine (with socle an abelian group). We will deal with these two cases in the next two sections, with results summarized in Tables \ref{tab0} and \ref{tab0a}, respectively.

\section{Almost simple case}
\label{sec:almost simple case}

In this section we assume that $G\leq\Sym(V)$ is $2$-transitive on $V$ of degree $u:=|V|$ with $\soc(G)$ a nonabelian simple group. It is well known (\cite{Kantor85}, \cite[p.196]{Cameron99}, \cite{Cameron81}) that $\soc(G)$ and $u$ are as follows:
\begin{description}
\item[\rm (i)] $\soc(G) = A_u$, $u \geq 5$;

\item[\rm (ii)] $\soc(G) = \PSL(d, q)$, $d \geq 2$, $q$ is a prime power and $u = (q^d-1)/(q-1)$, where $(d, q)\neq$ $(2,2), (2,3)$;

\item[\rm (iii)] $\soc(G) = \PSU(3, q)$, $q \geq 3$ is a prime power and $u = q^3 + 1$;

\item[\rm (iv)]  $\soc(G) = \Sz(q)$, $q = 2^{2e+1} > 2$ and $u = q^2+1$;

\item[\rm (v)] $\soc(G) = \R(q)'$, $q = 3^{2e+1} $ and $u = q^3+1$;

\item[\rm (vi)] $G = \Sp_{2d}(2)$, $d \geq 3$ and $u = 2^{2d-1} \pm 2^{d-1}$;

\item[\rm (vii)] $G=\PSL(2,11)$, $u = 11$;

\item[\rm (viii)] $\soc(G) = M_u$, $u = 11,12,22,23,24$;

\item[\rm (ix)] $G = M_{11}$, $u = 12$;

\item[\rm (x)] $G=A_7$, $u=15$;

\item[\rm (xi)] $G = \HS$, $u = 176$;

\item[\rm (xii)] $G = \Co_3$, $u = 276$.
\end{description}

Let $\s, \t$ be distinct points in $V$. In cases (i), (viii) and (ix), since $\soc(G)$ is $3$-transitive, a $2$-design as in Lemma \ref{lem:1 basic} admitting $G$ as a group of automorphisms does not exist. In case (vii), since $G_{\s,\t}$ has orbit-lengths $3$ and $6$ on $V \setminus \{\s,\t \}$ (see \cite{Kantor85}), again a $2$-design as in Lemma \ref{lem:1 basic} does not exist. In case (xi), since $G_{\s,\t}$ has orbit-lengths $12$, $72$ and $90$ on $V \setminus \{\s,\t \}$ (see \cite{Kantor85}), there is no $2$-design as in Lemma \ref{lem:1 basic}. Similarly, in case (xii) such a $2$-design does not exist as $G_{\s,\t}$ has orbit-lengths $112$ and $162$ on $V \setminus \{\s,\t \}$ (see \cite{Kantor85}).

In case (ii) with $d=2$, let $\s=\la \bfe_1\ra$ and $\t=\la \bfe_2\ra$, where $\bfe_1=(1,0)$ and $\bfe_2=(0,1)$. If $\vp\in \PSL(2,q)_{\s,\t}$, then $\bfe_1^\vp=(\mu,0)$ and $\bfe_2^\vp=(0,1/\mu)$ for some $\mu\in \mathbb{F}_q^\times$. Let $\eta=\la (x,1)\ra\in V\setminus\{\s,\t\}$, where $x\in\mathbb{F}_q^\times$. Then $\eta^\vp=\la(x\mu^2, 1)\ra$ and the $\PSL(2,q)_{\s,\t}$-orbit on $V$ containing $\eta$ has length $q-1$ when $q$ is even and $(q-1)/2$ when $q$ is odd. Thus, if $P$ is a nontrivial imprimitive block of $\PSL(2,q)_\s$ on $V\setminus\{\s\}$ containing $\t$, then $q$ is odd and $|P|=1+(q-1)/2=(q+1)/2$. Since $|P|$ divides $u-1=q$, we must have $q=1$, a contradiction. Hence this case does not produce any $G$-flag graph.

In case (ii) with $d \ge 3$, $G_{\s,\t}$ has orbit-lengths $q-1$ and $u-(q+1)$ on $V \setminus \{\s,\t\}$. Thus $\DD=\PG(d-1,q)$ and $\l=1$ by \cite{Kantor85}. The imprimitive block of $G_\s$ on $V \setminus \{\s\}$ containing $\t$ has size $q$ and so $\Ga(\DD,\Om,\Xi)\cong ((q^d-1)(q^{d-1}-1)/(q^2-1)(q-1))\cdot K_{q+1}$, contributing to (b) in Table \ref{tab0}.

In case (vi), $G_\s$ acts on $V \setminus \{\s\}$ as $O^\pm(2d,2)$ does on its singular vectors (see \cite{Kantor85}), and $G_{\s,\t}$ has orbit-lengths $2(2^{d-1} \mp 1)(2^{d-2} \pm 1)$ and $2^{2d-2}$ on $V \setminus \{\s,\t \}$. (We may assume $d \ge 3$ since the case $d < 3$ is covered by other cases in the classification of finite $2$-transitive groups.) Since the size of a $G_{\s,\t}$-orbit on $V \setminus \{\s,\t \}$ plus $1$ cannot divide $u-1$, there is no $2$-design as in Lemma \ref{lem:1 basic} admitting $G$ as a group of automorphisms.

In case (x), since $G_{\s,\t}$ has orbit-lengths $1$ and $12$ on $V \setminus \{\s,\t \}$, we have $\DD=\PG(3,2)$ and $\l=1$ (see \cite{Kantor85}). Since the imprimitive block of $G_\s$ on $V \setminus \{\s\}$ containing $\t$ has size $2$, we have $\Ga(\DD,\Om,\Xi)\cong 35\cdot K_3$, contributing to (a) in Table \ref{tab0}.

\subsection{$\soc(G) = \PSU(3, q)$, $q \geq 3$ a prime power and $u = q^3 + 1$}

This case yields (c) in Table \ref{tab0}. We use the following permutation representation of $\PSU(3, q)$ (see \cite[pp.248--249]{Dixon-Mortimer}). Let $W$ be a $3$-dimensional vector space over $\mathbb{F}_{q^2}$. Using $\xi \mapsto \overline{\xi} = \xi^q$ to denote the automorphism of $\mathbb{F}_{q^2}$ of order $2$, we define a hermitian form $\varphi : W \times W \rightarrow \mathbb{F}_{q^2}$ by $\varphi(w, z) = \xi_1\overline{\eta_3} + \xi_2\overline{\eta_2} + \xi_3\overline{\eta_1}$, where $w = (\xi_1, \xi_2, \xi_3)$ and $z = (\eta_1, \eta_2, \eta_3)$ are vectors in $W$. One can see that for this form the set of $1$-dimensional isotropic subspaces is given by
$$
V=\{\langle (1,0,0) \rangle\} \cup \{\langle (\a, \b, 1) \rangle \;|\; \a + \overline{\a} + \b \overline{\b} =0, \a, \b \in \mathbb{F}_{q^2} \}.
$$
(A vector $w \in W$ is \emph{isotropic} if $\varphi(w, w) = 0$.) We have $|V| = q^3 +1$ and $\PSU(3, q)$ is $2$-transitive on $V$. Let $\pi: {\GU(3,q)} \rightarrow {\PGU(3,q)}= {\GU(3,q)}/Z({\GU(3,q)})$ be the natural homomorphism, where $Z({\GU(3,q)})$ is the center of ${\GU(3,q)}$. Denote
\begin{equation*}
\begin{aligned}
t_{\a,\b}:=
\begin{bmatrix}
1 & -\overline{\b} & \a \\
0 & 1 & \b \\
0 & 0 & 1
\end{bmatrix}
\end{aligned}
,\;\,
\begin{gathered}
h_{\gamma,\delta}:=
\begin{bmatrix}
\gamma & 0 & 0 \\
0 & \delta & 0 \\
0 & 0 & \overline{\gamma}^{-1}
\end{bmatrix},
\;\, \mbox{where $\a,\b, \g,\d \in \mathbb{F}_{q^2}$}.
\end{gathered}
\end{equation*}
\noindent If $\d\overline{\d}=1$, $\g \neq 0$ and $\a + \overline{\a} + \b \overline{\b} =0$, then $t_{\a,\b}$ and $h_{\gamma,\delta}$ define elements of $\GU(3, q)$. There are $q^3$ matrices of type $t_{\a,\b}$ and $(q^2-1)(q+1)$ of type $h_{\g,\d}$. Let $\bfe_1 = (1,0,0)$ and $\bfe_3 = (0,0,1)$. Then ${\PGU(3,q)}_{\la \bfe_1 \ra} = \{\pi(h_{\g,\d}) \pi(t_{\a,\b})\;|\; \a,\b, \g,\d \in \mathbb{F}_{q^2}, \d\overline{\d}=1, \g \neq 0, \a + \overline{\a} + \b \overline{\b} =0\}$ and ${\GU(3,q)}_{\la \bfe_1 \ra, \la \bfe_3 \ra} = \{ h_{\g,\d}\;|\; \g,\d \in \mathbb{F}_{q^2}, \d\overline{\d}=1, \g \neq 0\}$. Obviously $t_{\a,\b} \in \SU(3,q)$, and $h_{\g,\d} \in \SU(3,q)$ if and only if $\d = \g^{q-1}$. Moreover, $h_{\g,\d} \in \SU(3,q)$ is a scalar matrix if and only if $\g^{q-2} =1$.

In the rest of this section we denote $\PSU(3,q)$ by $J$.

\begin{lemma}
\label{lem:PSU}
(\cite[Lemma 3.2]{FXZ}) Let $\la (\eta_1,\eta_2,1) \ra$ $\in V \setminus \{\la \bfe_1\ra, \la \bfe_3 \ra\}$. Denote by $Q$ the $J_{\la \bfe_1 \ra, \la \bfe_3 \ra}$-orbit containing $\la (\eta_1,\eta_2,1) \ra$. If $\eta_2 = 0$, then $|Q|=q-1$; if $\eta_2 \neq 0$, then
$$
|Q|=|J_{\la \bfe_1 \ra, \la \bfe_3 \ra}|=
\begin{cases}
q^2-1 \text{,\quad if\; $3 \nmid (q+1),$} \\
(q^2-1)/3 \text{,\quad if \;$3 \mid (q+1)$.}
\end{cases}
$$
\end{lemma}

Suppose that $P$ is a nontrivial imprimitive block of $J_{\la \bfe_1 \ra}$ on $V \setminus \{\la \bfe_1 \ra\}$ containing $\la \bfe_3 \ra$. Then $P \setminus \{\la \bfe_3 \ra\}$ is the union of some $J_{\la \bfe_1 \ra, \la \bfe_3 \ra}$-orbtis on $V \setminus \{\la \bfe_1 \ra, \la \bfe_3 \ra\}$. By Lemmas \ref{lem:6 num1} and \ref{lem:PSU}, we have $|P| = q$ or $q^2$.

\begin{lemma}
\label{lem:PSU(3,q), nontrivial imprimitive blocks}
Let $P$ be a nontrivial imprimitive block of $J_{\la \bfe_1 \ra}$ on $V \setminus \{\la \bfe_1 \ra\}$ containing $\la \bfe_3 \ra$. Then $P=\{\la (\a,0,1)\ra\;|\;\a+\overline{\a}=0\}$. Moreover, let $\DD:=(V,L^J)$ and $\Om:=(\la \bfe_1\ra, L)^J$, where $L:=P\cup\{\la \bfe_1\ra\}$, and let $G\leq \PGammaU(3, q)$ with $\soc(G)=J$. Then $\DD$ is a $2$-$(q^3+1, q+1, 1)$ design admitting $G$ as a group of automorphisms, $\Om$ is a 1-feasible $G$-orbit on the flag set of $\DD$, and $\Ga(\DD,\Om,\Xi)\cong (q^4-q^3+q^2)\cdot K_{q+1}$.
\end{lemma}

\proof
We first prove that $|P| \ne q^2$. Suppose otherwise. Denote the $q$ solutions in $\mathbb{F}_{q^2}$ of the equation $x + \overline{x} = 0$ by $\ve_0 = 0$, $\ve_1$, $\ldots$, $\ve_{q-1}$. Then $\la (\ve_1, 0,1) \ra$, $\ldots$, $\la (\ve_{q-1}, 0,1) \ra$ is a $J_{\la \bfe_1 \ra, \la \bfe_3 \ra}$-orbit on $V \setminus \{\la \bfe_1 \ra, \la \bfe_3 \ra\}$. By Lemma \ref{lem:6 num1}, $\la (\ve_i, 0,1)\ra$ is not contained in $P$ for $i > 0$. Now $\Sigma:=\{P^g\;|\;g \in J_{\la \bfe_1 \ra}\}$ is a system of blocks of $J_{\la \bfe_1 \ra}$ on $V \setminus \{\la \bfe_1 \ra\}$ with $|\Sigma|=q$, and $T := \la \pi(t_{\a,\b})\;|\; \a + \overline{\a} + \b \overline{\b} =0\ra$ is transitive on $\Sigma$. Actually, $T$ is a normal subgroup of $J_{\la \bfe_1 \ra}$ acting regularly on $V \setminus \{\la \bfe_1 \ra\}$ (see \cite[p.249]{Dixon-Mortimer}). Thus, the stabilizer of $P$ in $T$ has order $q^2$, that is, $|T_P|=q^2$. For $\pi(t_{\a_1, \b})$, $\pi(t_{\a_2, \b}) \in T_P$, we have ${\la (0, 0,1) \ra}^{\pi(t_{\a_1,\b}) \pi(t_{\a_2, \b}^{-1})}= {\la(\a_1, \b,1)\ra}^{\pi(t_{-\a_2 - \b\overline{\b}, -\b})} =\la(\a_1-\a_2, 0,1)\ra \in P$. Since $\la(\ve_i, 0,1)\ra$ is not contained in $P$ for $i > 0$, we have $\a_1 = \a_2$ and $\pi(t_{\a_1,\b}) = \pi(t_{\a_2, \b})$. Therefore,
\begin{equation}
\label{equ:PSU}
\{\b \;| \;\la (\a,\b,1)\ra \in P \} = \{\b \;|\; \pi(t_{\a,\b}) \in T_P\} = \mathbb{F}_{q^2}.
\end{equation}

For any $\la(\eta_1, \eta_2,1)\ra$, $\la(\xi_1, \xi_2,1)\ra \in P$, $\eta_2$, $\xi_2 \neq 0$, since $P$ is an imprimitive block of $J_{\la \bfe_1 \ra}$ on $V \setminus \{\la \bfe_1 \ra\}$, both $\pi(t_{\eta_1,\eta_2})$ and $\pi(t_{\xi_1, \xi_2})$ fix $P$ setwise. Thus
$$
{\la(0,0,1)\ra}^{\pi(t_{\eta_1,\eta_2})\pi(t_{\xi_1,\xi_2})}={\la (\eta_1, \eta_2,1)\ra}^{\pi(t_{\xi_1, \xi_2})}=\la(\eta_1+\xi_1-\overline{\xi_2}\eta_2,\eta_2+\xi_2,1)\ra \in P,
$$
$$
{\la(0,0,1)\ra}^{\pi(t_{\xi_1, \xi_2})\pi(t_{\eta_1,\eta_2})}={\la(\xi_1, \xi_2,1)\ra}^{\pi(t_{\eta_1, \eta_2})}=\la(\xi_1+\eta_1-\overline{\eta_2}\xi_2,\xi_2+\eta_2,1)\ra \in P.
$$
Hence $\eta_1+\xi_1-\overline{\xi_2}\eta_2= \xi_1+\eta_1-\overline{\eta_2}\xi_2$ by \eqref{equ:PSU}, that is,  ${(\xi_2/{\eta_2})}^{q-1} = 1$, which implies $(\xi_2/{\eta_2}) \in {\rm Fix}_f(\mathbb{F}_{q^2})$, where $f$ denotes the automorphism of $\mathbb{F}_{q^2}$ of order $2$. Choosing $\eta_2 =1$, we have $\xi_2 \in {\rm Fix}_f(\mathbb{F}_{q^2})$ and so $\mathbb{F}_{q^2} \subseteq {\rm Fix}_f(\mathbb{F}_{q^2})$, a contradiction. Thus $|P|\neq q^2$.

So we have $|P|=q$. By Lemma \ref{lem:PSU}, $P=\{\la (\a,0,1)\ra\;|\;\a+\overline{\a}=0\}$. Since $G\leq \PGammaU(3, q)$, $P$ is also an imprimitive block of $G_{\la \bfe_1\ra}$ on $V \setminus \{\la \bfe_1 \ra\}$, and moreover $\Om=(\la \bfe_1\ra, L)^G$ and $L^G=L^J$. Hence $\DD$ is a $2$-$(q^3+1, q+1, \l)$ design admitting $G$ as a group of automorphisms with $\l=1$ (see \cite[p.249]{Dixon-Mortimer} and \cite{Kantor85}), $\Om$ is a 1-feasible $G$-orbit on the flag set of $\DD$, and $\Ga(\DD,\Om,\Xi)\cong (q^4-q^3+q^2)\cdot K_{q+1}$ by Corollary \ref{corol:k=1&complete quotient}.
\qed

\subsection{$\soc(G)=\Sz(q)$, $q = 2^{2e+1} > 2$ and $u = q^2+1$}

This case yields (d) in Table \ref{tab0}. We use the following permutation representation of $\Sz(q)$ (see \cite[p.250]{Dixon-Mortimer}). The mapping $\s : \xi \mapsto \xi^{2^{e+1}}$ is an automorphism of $\mathbb{F}_q$ and $\s^2$ is the Frobenius automorphism $\xi \mapsto \xi^{2}$. Define
\begin{equation}
\label{equ:point set of Sz(q)}
V:= \{(\eta_1, \eta_2, \eta_3) \in \mathbb{F}_q^3 \;|\; \eta_3 = \eta_1\eta_2+\eta_1^{\s+2}+\eta_2^{\s} \} \cup \{ \infty \}.
\end{equation}
Then $|V| = q^2+1$. For $\a, \b, \kappa \in \mathbb{F}_q$ with $\kappa \neq 0$, define the following permutations of $V$ fixing $\infty$:
\begin{align*}
t_{\a,\b}&:(\eta_1, \eta_2, \eta_3) \mapsto \left(\eta_1 + \a, \eta_2 + \b + \a^{\s}\eta_1, \mu\right), \\
n_\kappa &:(\eta_1, \eta_2, \eta_3) \mapsto (\kappa\eta_1, \kappa^{\s+1}\eta_2, \kappa^{\s+2}\eta_3),
\end{align*}
where $\mu=\eta_3+\a\b+\a^{\s+2}+\b^\s$ $+\a\eta_2+\a^{\s+1}\eta_1+\b\eta_1$. Define the involution $w$ fixing $V$ by
$$
w : (\eta_1, \eta_2, \eta_3) \leftrightarrow \left(\frac{\eta_2}{\eta_3}, \frac{\eta_1}{\eta_3},\frac{1}{\eta_3}\right) \text{ for $\eta_3 \neq 0$}, \;\;
\infty \leftrightarrow (0,0,0)=: \bzero.
$$
Then $\Sz(q)$ is generated by $w$ and all $t_{\a,\b}$ and $n_\kappa$. We have ${\Sz(q)}_{\infty}=\la t_{\a,\b}, n_\kappa\;|$ $\a, \b, \kappa \in \mathbb{F}_q$, $\kappa \neq 0 \ra$ and ${\Sz(q)}_{\infty, \bzero}=\la n_\kappa \;|$ $\kappa \in \mathbb{F}_q, \kappa \neq 0 \ra$, the latter being a cyclic group.

\begin{lemma}
\label{lem:orbit-length, Sz(q)}
Every ${\Sz(q)}_{\infty,\bzero}$-orbit on $V\setminus \{\infty,\bzero\}$ has length $q-1$.
\end{lemma}

\proof
Since $2^{e+1} + 1$ and $2^{2e+1} - 1$ are coprime and $\mathbb{F}_q^{\times}$ is a cyclic group of order $2^{2e+1} - 1$, the mapping $\mathbb{F}_q^{\times} \rightarrow \mathbb{F}_q^{\times}$, $z \mapsto z^{\s +1}$ is a group automorphism. Thus, if $\eta_1\neq0$ or $\eta_2\neq0$, then $(a\eta_1, a^{\s+1}\eta_2, a^{\s+2}\eta_3)= (b\eta_1, b^{\s+1}\eta_2, b^{\s+2}\eta_3)$ if and only if $a=b$, and so the result follows.
\qed

\begin{lemma}
\label{lem:Sz(q)}
Let $G$ be any subgroup of $\Sym(V)$ containing $\Sz(q)$ as a normal subgroup. Suppose that $P$ is a nontrivial imprimitive block of $G_{\infty}$ on $V\setminus \{\infty\}$ containing $\bzero$. Then the following hold:
\begin{itemize}
\item[\rm (a)] $P =\{(0,\eta,\eta^\s)\in V\,|\,\eta\in \mathbb{F}_q\}$ and ${\Sz(q)}_{\infty,P}=\la t_{0,\xi}, n_\kappa\, |\,\kappa\in \mathbb{F}_q^{\times}, \xi \in \mathbb{F}_q \ra$ (\cite[Lemma 3.6]{FXZ});
\item[\rm (b)] setting $L:=P \cup \{\infty\}$, $\DD := (V, L^{\Sz(q)}) = (V, L^{G})$ is a $2$-$(q^2+1, q+1, q+1)$ design admitting $G$ as a $2$-point-transitive and block-transitive group of automorphisms, and $\Om := {(\infty, L)}^{\Sz(q)} = {(\infty, L)}^{G}$ is a 1-feasible $G$-orbit on the flag set of $\DD$; moreover, the $G$-flag graph of $\DD$ with respect to $(\Om, \Xi)$ is the same as the $\Sz(q)$-flag graph $\Ga(\DD,\Om,\Xi)$;
\item[\rm (c)] $\Ga(\DD,\Om,\Xi)$ is connected with order $|\Om|=q(q^2+1)$ and valency $q$.
\end{itemize}
\end{lemma}

\proof
We prove (b) and (c) only, as (a) was proved in \cite[Lemma 3.6]{FXZ}. Since $\Sz(q)$ is $2$-transitive on $V$, $\DD$ is a $2$-$(q^2+1, q+1, \l)$ design admitting $\Sz(q)$ as a $2$-point-transitive and block-transitive group of automorphisms. Since $w$ does not stabilize $L$, $\l\neq1$ and thus $\l = q+1$ by Lemma \ref{lem:2 lambda}.

Since $\Sz(q)$ is a normal subgroup of $G$ and $\Sz(q)$ has index $2e+1$ in its normalizer $Q$ in $\Sym(V)$ (\cite[p.197, Table 7.4]{Cameron99}), $Q/\Sz(q)$ is a cyclic group of order $2e+1$ and $G = \la \Sz(q), \zeta \ra$, where $\zeta$ is an automorphism of $\mathbb{F}_q$ inducing a permutation of $V$ that fixes $\infty$ and acts on the elements of $V\setminus\{\infty\}$ componentwise. Hence $P = \{(0, \eta, \eta^\s)\in V\;|\;\eta \in \mathbb{F}_q\}$ is a nontrivial imprimitive block of $G_{\infty}$ on $V \setminus \{\infty\}$. Moreover, ${(\infty, L)}^G ={(\infty, L)}^{\Sz(q)}$ and $L^G = L^{\Sz(q)}$. Hence $\DD$ admits $G$ as an automorphism group and $\Om$ is a 1-feasible $G$-orbit on the flag set of $\DD$.

Denote $H:=\la t_{0,\xi}, w\; |\;\xi \in \mathbb{F}_q\ra$. For any $(\eta_1, \eta_2, \eta_3) \in V\setminus\{\infty, \bzero\}$, if $\eta_1=0$ then $\bzero^{t_{0, \eta_2}}=(\eta_1, \eta_2, \eta_3)$, and if $\eta_1\neq0$ then $\bzero^{t_{0, \theta} w t_{0, \eta_2}}=(\eta_1, \eta_2, \eta_3)$, where $\theta/\theta^\s=\eta_1$. Hence $H$ is transitive on $V$, and thus $(q^2+1)q$ divides $|H|$. So $|H|$ does not divide $q^2(q-1)$, $2(q-1)$, $4(q+\sqrt{2q}+1)$ or $4(q-\sqrt{2q}+1)$. Thus, by \cite[p.137, Theorem 9]{Suzuki62}, $|H|=(s^2+1)s^2(s-1)$, where $s^m=q$ for some positive integer $m$. It follows that $m=1$, $|H|=(q^2+1)q^2(q-1)$, and thus $\Sz(q)=H$. Therefore, $\Sz(q)=\la{\Sz(q)}_{\infty,P}, w\ra$ and so $\Ga(\DD,\Om,\Xi)$ is connected by Lemma \ref{lem:connectedness of G-flag graph}.
\qed

\subsection{$\soc(G)=\R(q)$, $q = 3^{2e+1} >3$, $u = q^3+1$; or $G =\R(3)$, $\R(3)' \cong \PSL(2,8)$, $u = 28$}
\label{subsec:Ree}

This case yields (e) in Table \ref{tab0}. We use the following permutation representation of $\R(q)$ (see \cite[p.251]{Dixon-Mortimer}). The mapping $\s : \xi \mapsto \xi^{3^{e+1}}$ is an automorphism of $\mathbb{F}_q$ and $\s^2$ is the Frobenius automorphism $\xi \mapsto \xi^{3}$. The set $V$ of points on which $\R(q)$ acts consists of $\infty$ and the set of 6-tuples $(\eta_1, \eta_2, \eta_3, \l_1, \l_2, \l_3)$ with $\eta_1, \eta_2, \eta_3 \in \mathbb{F}_q$ and
\begin{align}
\label{equ:point set of R(q)}
\begin{cases}
\l_1= \eta_1^2{\eta_2} - {\eta_1}\eta_3 + \eta_2^\s - \eta_1^{\s +3}, \\
\l_2= \eta_1^{\s}\eta_2^{\s} - \eta_3^{\s} + \eta_1 \eta_2^2 + {\eta_2}\eta_3- \eta_1^{2\s +3}, \\
\l_3= {\eta_1}\eta_3^{\s} - \eta_1^{\s+1} \eta_2^{\s} + \eta_1^{\s+3} {\eta_2} + \eta_1^2\eta_2^2 -\eta_2^{\s+1}-\eta_3^2+ \eta_1^{2\s +4}.
\end{cases}
\end{align}
Thus $|V| = q^3+1$. For $\a, \b, \g, \kappa \in \mathbb{F}_q$ with $\kappa \neq 0$, define the following permutations of $V$ fixing $\infty$:
\begin{align*}
t_{\a,\b,\g}&:(\eta_1, \eta_2, \eta_3, \l_1, \l_2, \l_3) \mapsto\\
&(\eta_1 + \a, \eta_2 + \b + \a^{\s}\eta_1, \eta_3+\g-\a\eta_2+\b\eta_1-{\a}^{\s+1}\eta_1, \mu_1, \mu_2, \mu_3), \\
n_\kappa &:(\eta_1, \eta_2, \eta_3, \l_1, \l_2, \l_3) \mapsto
(\kappa\eta_1, \kappa^{\s+1}\eta_2, \kappa^{\s+2}\eta_3, \kappa^{\s+3}\l_1, \kappa^{2\s+3}\l_2, \kappa^{2\s+4}\l_3),
\end{align*}
\noindent where $\mu_1, \mu_2, \mu_3$ can be calculated from the formulas in \eqref{equ:point set of R(q)}. Define the involution $w$ fixing $V$ by
\begin{align*}
w : (\eta_1, \eta_2, \eta_3, \l_1, \l_2, \l_3) & \leftrightarrow \left(\frac{\l_2}{\l_3}, \frac{\l_1}{\l_3}, \frac{\eta_3}{\l_3}, \frac{\eta_2}{\l_3}, \frac{\eta_1}{\l_3},\frac{1}{\l_3}\right) \text{ for $\l_3 \neq 0$,} \\
\infty & \leftrightarrow (0,0,0,0,0,0)=: \bzero.
\end{align*}
\noindent (We use the corrected definition \cite{Dixon-Mortimer-errata} of the action of $w$ on $V$ \cite[p.251]{Dixon-Mortimer}.) The Ree group $\R(q)$ is the group generated by $w$ and all $t_{\a,\b,\g}$ and $n_\kappa$. We have ${\R(q)}_{\infty} = \la t_{\a,\b,\g}, n_\kappa \;|\; \a, \b, \g, \kappa \in \mathbb{F}_q, \kappa \neq 0 \ra$ and ${\R(q)}_{\infty, \bzero}$ is the cyclic group $\la n_\kappa \;| \;\kappa \in \mathbb{F}_q, \kappa \neq 0 \ra$. Since the first three coordinates in each element of $V$ determine the other three, we present an element of $V$ by $(\eta_1, \eta_2, \eta_3, \ldots)$. The following lemma was used in \cite{FXZ}, and its proof is given below for later reference.

\begin{lemma}
\label{lem:orbit-length, R(q)}
(\cite[Lemma 3.10]{FXZ}) Let $(\eta_1, \eta_2, \eta_3, \ldots)\in$ $V\setminus \{\infty,\bzero\}$. Then
\begin{equation*}
|{(\eta_1, \eta_2, \eta_3, \ldots)}^{{\R(q)}_{\infty,\bzero}}|=
\begin{cases}
q-1 \text{,\quad if \;$\eta_1 \neq 0$ or $\eta_3 \neq 0$,} \\
(q-1)/{2} \text{,\quad if \;$\eta_1 = \eta_3 =0$.}
\end{cases}
\end{equation*}
\end{lemma}

\proof
Since $\id: \mathbb{F}_q^{\times} \rightarrow \mathbb{F}_q^{\times}$, $\xi \mapsto \xi$ and $\varphi: \mathbb{F}_q^{\times} \rightarrow \mathbb{F}_q^{\times}$, $\xi \mapsto \xi^{\s+2}$ are both group automorphisms, if $\eta_1 \neq 0$ or $\eta_3 \neq 0$, then $(a\eta_1, a^{\s+1}\eta_2, a^{\s+2}\eta_3,\ldots)= (b\eta_1, b^{\s+1}\eta_2, b^{\s+2}\eta_3,\ldots)$ if and only if $a=b$.

Let $\d$ be a generator of the cyclic group $\mathbb{F}_q^{\times}$. Since $\d^{\s+1} = \d^{3^{e+1}+1}$ and $\gcd(3^{e+1}+1, q-1) = 2$, we have $|\d^{\s+1}|= ({q-1})/{2}$, and thus
\begin{align}
\label{equ:L1 and L2}
L_1 := {(0,1,0,\ldots)}^{{\R(q)}_{\infty,\bzero}}\; \text{ and }\; L_2 := {(0,\d,0,\ldots)}^{{\R(q)}_{\infty,\bzero}}
\end{align}
\noindent are $\R(q)_{\infty,\bzero}$-orbits on $V\setminus \{\infty,\bzero\}$ of length $({q-1})/{2}$.
\qed

\smallskip
By Lemma \ref{lem:orbit-length, R(q)}, ${\R(q)}_{\infty,\bzero}$ has two orbits of length $({q-1})/{2}$ and $q(q+1)$ orbits of length $q-1$ on $V\setminus \{\infty,\bzero\}$.

In the remainder of this section, let $G\leq\Sym(V)$ contain $\R(q)$ as a normal subgroup and $P$ be a nontrivial imprimitive block of ${G}_{\infty}$ on $V\setminus \{\infty\}$ containing $\bzero$. Since $\R(q)$ has index $2e+1$ in its normalizer $Q$ in $\Sym(V)$ (\cite[p.197, Table 7.4]{Cameron99}), $Q/\R(q)$ is a cyclic group of order $2e+1$ and
\be
\label{rqzeta}
G = \la \R(q), \zeta \ra,
\ee
where $\zeta$ is an automorphism of $\mathbb{F}_q$ inducing a permutation of $V$ that fixes $\infty$ and acts on the elements of $V\setminus\{\infty\}$ componentwise. We have
\be
\label{eq:rq}
G_\infty=\la \R(q)_\infty, \zeta\ra,\; G_{\infty, \bzero}=\la \R(q)_{\infty, \bzero}, \zeta\ra.
\ee
By Lemma \ref{lem:1 basic}, $P \setminus  \{\bzero\}$ is the union of some ${\R(q)}_{\infty,\bzero}$-orbits on $V\setminus \{\infty,\bzero\}$. By Lemmas \ref{lem:8 num3} and \ref{lem:9 num4}, we have $|P|= q$ or $|P| = q^2$. If $|P|=q^2$, then by Lemma \ref{lem:9 num4}, either $L_1 \cup L_2\subseteq P$ or $(L_1 \cup L_2)\cap P=\emptyset$, where $L_1$ and $L_2$ are as defined in (\ref{equ:L1 and L2}).

\medskip
\textsf{Case (i):} $|P|=q$, $P\setminus\{\bzero\}=L_1\cup L_2$

\medskip
Now $P=\{(0,\eta,0,\ldots)\;|\;\eta\in \mathbb{F}_q\}$. Since $G_{\infty} =\la t_{\a,\b,\g}, n_\kappa, \zeta \;| \;\a,\b,\g\in \mathbb{F}_q,\kappa\in \mathbb{F}_q^{\times}\ra$ by (\ref{eq:rq}) and ${(0,\eta,0,\ldots)}^{t_{\a,\b,\g}}= (\a, \eta +\b, \g- \a \eta,\ldots)$, $P$ is an imprimitive block of $G_{\infty}$ on $V\setminus\{\infty\}$. Denote $L := P \cup \{\infty\}$, and define $\DD :=(V, L^{G})=(V, L^{\R(q)})$ and $\Om := {(\infty, L)}^{G}={(\infty, L)}^{\R(q)}$. Then $\DD$ is a $2$-$(q^3+1, q+1, \l)$ design admitting $G$ as a group of automorphisms, and $\Om$ is a 1-feasible $G$-orbit on the flag set of $\DD$. It can be verified that $w$ stabilizes $L$ and thus $\l=1$. It follows that $\Ga(\DD,\Om,\Xi)\cong (q^4-q^3+q^2)\cdot K_{q+1}$, yielding the first possibility in (e) in Table \ref{tab0}.

\medskip
\textsf{Case (ii):} $|P|=q$, $P\setminus\{\bzero\}$ is an ${\R(q)}_{\infty,\bzero}$-orbit on $V\setminus \{\infty,\bzero\}$ of length $q-1$

\medskip
Now each element of $P$ is of the form $(\kappa\eta_1, \kappa^{\s+1}\eta_2, \kappa^{\s+2}\eta_3, \ldots)$ for some $\kappa \in \mathbb{F}_q$, where $(\eta_1, \eta_2, \eta_3, \ldots)$ is a fixed point in $P\setminus\{\bzero\}$. Suppose that $P\cap P^{t_{\a,\b,\g}}\neq \emptyset$, that is, for some $\kappa_0,\kappa_1 \in \mathbb{F}_q$,
\begin{align}
\label{equ:R(q) nonempty intersection}
{(\kappa_1\eta_1, \kappa_1^{\s+1}\eta_2, \kappa_1^{\s+2}\eta_3, \ldots)}^{t_{\a,\b,\g}}=(\kappa_0\eta_1, \kappa_0^{\s+1}\eta_2, \kappa_0^{\s+2}\eta_3, \ldots).
\end{align}
Then $\kappa_1\eta_1+\a =\kappa_0\eta_1$, $\kappa_1^{\s+1}\eta_2 + \b + \a^{\s}\kappa_1\eta_1 = \kappa_0^{\s+1}\eta_2$ and $\kappa_1^{\s+2}\eta_3+\g-\a\kappa_1^{\s+1}\eta_2+\b\kappa_1\eta_1-{\a}^{\s+1}\kappa_1\eta_1=\kappa_0^{\s+2}\eta_3$, or equivalently,
\begin{align}
\label{equ:value of alpha,beta,gamma}
\begin{cases}
\a =(\kappa_0-\kappa_1)\eta_1,\\
\b = (\kappa_0^{\s+1}-\kappa_1^{\s+1})\eta_2 - (\kappa_0^{\s} - \kappa_1^{\s})\kappa_1 \eta_1^{\s+1},\\
\g =(\kappa_0^{\s+2}-\kappa_1^{\s+2})\eta_3+(\kappa_0\kappa_1^{\s+1}-\kappa_1\kappa_0^{\s+1})\eta_1\eta_2+ (\kappa_0^{\s}-\kappa_1^{\s})\kappa_0\kappa_1\eta_1^{\s+2}.
\end{cases}
\end{align}
Hence, if $\a,\b,\g$ are given by \eqref{equ:value of alpha,beta,gamma} in terms of $\eta_1,\eta_2,\eta_3$, then \eqref{equ:R(q) nonempty intersection} holds. Since $P$ is an imprimitive block of ${\R(q)}_{\infty}$ on $V\setminus \{\infty\}$, we need to verify that $P^{t_{\a,\b,\g}} = P$, that is, for every $\ell \in \mathbb{F}_q$, the following equation system has a solution $x \in \mathbb{F}_q$:
\begin{align}
\label{equ:R(q) any l}
\begin{cases}
(x-\ell)\eta_1 = \a,\\
(x^{\s+1}-\ell^{\s+1})\eta_2 - (x^{\s} - \ell^{\s})l \eta_1^{\s+1} = \b,\\
(x^{\s+2}-\ell^{\s+2})\eta_3+(x\ell^{\s+1}-\ell x^{\s+1})\eta_1\eta_2+ (x^{\s}-\ell^{\s})x\ell\eta_1^{\s+2} = \g.
\end{cases}
\end{align}

\begin{lemma}
\label{lem:R(q), eta1=eta2=0}
If \eqref{equ:R(q) any l} has a solution for every $\ell \in \mathbb{F}_q$, then $\eta_1 = \eta_2 = 0$.
\end{lemma}
\proof
If $P^{t_{\eta,\theta,\xi}} \cap P = \emptyset$ for any $t_{\eta,\theta,\xi} \neq \id$, then different $t_{\eta,\theta,\xi}$ must map $P$ to different elements in $P^{{\R(q)}_\infty}$, and thus $q^3 = |\la t_{\eta,\theta,\xi} \;|\; \eta,\theta,\xi \in \mathbb{F}_q \ra| \leq |P^{{\R(q)}_{\infty}}| = q^2$, a contradiction. Hence we can assume that at most two of $\a,\b$ and $\g$ are $0$ in \eqref{equ:R(q) nonempty intersection}. We claim that $\eta_1 = 0$. Suppose otherwise. Then $x= {\a}/{\eta_1}+\ell$, $\a\neq0$ (for otherwise $x=\ell$ and $\b=\g=0$ by \eqref{equ:R(q) any l}), and the second equation of \eqref{equ:R(q) any l} becomes ${\frac{\a\eta_2}{\eta_1}}\ell^{\s}+\left(\frac{\a^{\s}\eta_2}{\eta_1^{\s}} - \a^{\s} \eta_1\right)\ell + \frac{\a^{\s+1}\eta_2}{\eta_1^{\s+1}} - \b = 0$, which holds for every $\ell \in \mathbb{F}_q$. On the other hand, the polynomial of $\ell$ on the left-hand side should have at most $3^{e+1}$ roots if it is nonzero. Thus, if $q>3$ (so that $q > 3^{e+1}$), then this polynomial must be the zero polynomial and hence ${\a\eta_2}/{\eta_1}={\a^{\s}\eta_2}/{\eta_1^{\s}} - \a^{\s} \eta_1=0$, which implies $\a=0$, a contradiction. Assume $q=3$. Then $\s = \id$, and ${\frac{\a\eta_2}{\eta_1}}\ell + \left(\frac{\a\eta_2}{\eta_1} - \a \eta_1\right)\ell + \frac{\a^{2}\eta_2}{\eta_1^{2}} - \b = 0$ holds for every $\ell \in \mathbb{F}_3$, which implies ${\a\eta_2}/{\eta_1}+ ({\a\eta_2}/{\eta_1}- \a\eta_1)=0$. Hence $\eta_2=-1$ and the third equation of \eqref{equ:R(q) any l} becomes $\a \ell^2+\a^2\ell/\eta_1+\g-\a\eta_3/\eta_1=0$, which cannot hold for every $\ell\in \mathbb{F}_3$, a contradiction. Therefore, $\eta_1=0$. Consequently, we have $\a=0$ by \eqref{equ:value of alpha,beta,gamma}.

If $\eta_3=0$, then $\g=0$ by the third equation of \eqref{equ:value of alpha,beta,gamma}, and the second equation of \eqref{equ:R(q) any l} becomes $(x^{\s+1}- \ell^{\s+1}) \eta_2 =\b$, which has a solution $x\in \mathbb{F}_q$ for every $\ell \in \mathbb{F}_q$. By our assumption, $\b \neq 0$. For $\ell=0$ there exists $t \in \mathbb{F}_q^{\times}$ such that $t^{\s+1} =\b/{\eta_2}$. For $\ell=t$ there exists $z \in \mathbb{F}_q$ such that $z^{\s+1} = t^{\s+1}+\b/{\eta_2}= t^{\s+1}+t^{\s+1}= -t^{\s+1}$. Thus ${({z}/{t})}^{\s+1} =-1$. Let $\d$ be a generator of the cyclic group $\mathbb{F}_q^{\times}$ and set ${z}/{t}= {\d}^n$ ($n>0$). Then $-1={(\d^n)}^{\s+1} = {(\d^{\s+1})}^n$ has order $2$ in $\mathbb{F}_q^{\times}$. Since $|\d^{\s+1}|= {(q-1)}/{2}$, we have $2= |{(\d^n)}^{\s+1}|= {(q-1)}/{(2\cdot \gcd({(q-1)}/{2}, n))}$ and $\gcd({(q-1)}/{2}, n)= {(q-1)}/{4}$. But $4 \nmid (q-1)$ as $q$ is an odd power of $3$, a contradiction. Therefore $\eta_3\neq0$, and the third equation of \eqref{equ:R(q) any l} becomes $x^{\s+2} ={\g}/{\eta_3} + \ell^{\s+2}$. Since $x{(x^{\s+2})}^{\s-2}= {(x^{\s^2-3})}=1$ for $x \in \mathbb{F}_q^{\times}$, we have $x={{(x^{\s+2})}^2}/{{(x^{\s+2})}^{\s}}= {{({\g}/{\eta_3}+\ell^{\s+2})}^2}/{{({\g}/{\eta_3} + \ell^{\s+2})}^{\s}}$ when ${\g}/{\eta_3} + \ell^{\s+2}\neq0$, and $x^{\s+1}({\g}/{\eta_3} + \ell^{\s+2})= {({\g}/{\eta_3} + \ell^{\s+2})}^{\s}$. The latter also holds when ${\g}/{\eta_3} + \ell^{\s+2}=0$. Multiplying $({\g}/{\eta_3} + \ell^{\s+2})$ on both sides of the second equation of  \eqref{equ:R(q) any l}, we obtain
\begin{equation}
\label{equ:2nd equ of}
\b\ell^{\s+2} + \frac{\g \eta_2}{\eta_3} \ell^{\s+1} + \frac{\b\g}{\eta_3} - \frac{\g^{\s} \eta_2}{\eta_3^{\s}}=0,
\end{equation}
\noindent which holds for every $\ell \in \mathbb{F}_q$. Since the polynomial of $\ell$ on the left-hand side of this equation has at most $3^{e+1} + 2$ roots if it is nonzero, we must have $\b=0$ and $\g\eta_2/{\eta_3} =0$ if $q > 3$. If $q=3$, then $\s=\id$ and \eqref{equ:2nd equ of} becomes $\b\ell + \frac{\g \eta_2}{\eta_3} \ell^2 + \frac{\b\g}{\eta_3} - \frac{\g \eta_2}{\eta_3}=0$, and so we still have $\b=0$ and $\g\eta_2/{\eta_3} =0$. By our assumption, $\g \neq 0$ and thus $\eta_2 = 0$.
\qed

\begin{lemma}
\label{lem:R(q)}
Let $G = \la \R(q), \zeta\ra$ (as given in \eqref{rqzeta}) be a subgroup of $\Sym(V)$ containing $\R(q)$ as a normal subgroup.
Suppose that $P$ is an imprimitive block of $G_\infty$ on $V\setminus\{\infty\}$ containing $\bzero$ such that $P\setminus\{\bzero\}$ is an ${\R(q)}_{\infty,\bzero}$-orbit on $V\setminus \{\infty,\bzero\}$ of length $q-1$. Then the foollowing hold:
\begin{itemize}
\item[\rm (a)] $P =\{(0,0,\eta,\ldots)\in V\;|\;\eta\in \mathbb{F}_q\}$ and $G_{\infty,P}= \la n_\kappa, t_{0,0,\xi}, \zeta \;|\; \xi \in \mathbb{F}_q, \kappa \in \mathbb{F}_q^{\times}\ra$;
\item[\rm (b)] setting $L := P \cup \{\infty\}$, $\DD:=(V, {L}^{\R(q)}) = (V, {L}^{G})$ is a $2$-$(q^3+1, q+1, q+1)$ design admitting $G$ as a $2$-point-transitive and block-transitive group of automorphisms, and $\Om:={(\infty, L)}^{\R(q)}={(\infty, L)}^{G}$ is a 1-feasible $G$-orbit on the flag set of $\DD$; moreover, the $G$-flag graph of $\DD$ with respect to $(\Om, \Xi)$ is the same as the $\R(q)$-flag graph $\Ga(\DD,\Om,\Xi)$ with respect to $(\Om, \Xi)$;
\item[\rm (c)] $\Ga(\DD,\Om,\Xi)$ is connected with order $|\Om|=q^2(q^3+1)$ and valency $q$ when $q>3$ and has three connected components when $q=3$.
\end{itemize}
\end{lemma}

\proof
(a) By Lemma \ref{lem:orbit-length, R(q)} and Lemma \ref{lem:R(q), eta1=eta2=0}, $P =\{(0,0,\eta,\ldots)\in V\;|\;\eta\in \mathbb{F}_q\}$. By \eqref{eq:rq} one can check that $P$ is indeed an imprimitive block of $G_\infty$ on $V\setminus\{\infty\}$, and $G_{\infty,P}= \la n_\kappa, t_{0,0,\xi}, \zeta \;|\; \xi \in \mathbb{F}_q, \kappa \in \mathbb{F}_q^{\times}\ra$.

(b) Since $w$ does not stabilize $L$, we have $\l = q+1$ for $\DD$.

(c) We first recall the following known result (see \cite[p.60, Theorem C]{Kleidman88} or \cite[p.3758, Lemma 2.2]{XGfang-Praeger}): For any subgroup $H$ of $\R(q)$, either $|H|=(s^3+1)s^3(s-1)$, where $s^m=q$ for some positive integer $m$, or $|H|$ divides $q^3(q-1)$, $12(q+1)$, $q^3-q$, $6(q+\sqrt{3q}+1)$, $6(q-\sqrt{3q}+1)$, $504$ or $168$.

By Lemma \ref{lem:connectedness of G-flag graph}, it suffices to prove $\R(q)=\la t_{0,0,\xi}, n_\kappa, w\;|\;\xi\in \mathbb{F}_q, \kappa\in \mathbb{F}_q^\times\ra$ when $q>3$ and $\R(3)'=\la t_{0,0,\xi}, n_{-1}, w\;|\;\xi\in \mathbb{F}_3\ra$. Denote $H:=\la t_{0,0,\xi},n_\kappa, w\; |\;\xi \in \mathbb{F}_q, \kappa\in \mathbb{F}_q^\times \ra$. Since for $\xi,\eta\in \mathbb{F}_q^\times$ and $\theta\in \mathbb{F}_q$,
\begin{align}
\label{equ:R(q) is generated by}
&\bzero\stackrel{t_{0,0,\xi}}{\mapsto}(0,0,\xi,0,-\xi^\s,-\xi^2)\stackrel{w}{\mapsto}\left(\xi^{\s-2},0,-\frac{1}{\xi},0,0,-\frac{1}{\xi^2}\right)\stackrel{t_{0,0,\theta}}{\mapsto}\left(\xi^{\s-2},0,\theta-\frac{1}{\xi},\ldots\right),\notag\\
&(\eta,0,0,-\eta^{\s+3},-\eta^{2\s+3},-\eta^{2\s+4})\stackrel{w}{\mapsto}\left(-\frac{1}{\eta},-\frac{1}{\eta^{\s+1}},0,\ldots\right)
\stackrel{t_{0,0,\theta}}{\mapsto}\left(-\frac{1}{\eta},-\frac{1}{\eta^{\s+1}},\theta,\ldots\right),
\end{align}
we see that $\{(\zeta,0,\eta,\ldots)\in V\;|\;\zeta,\eta\in \mathbb{F}_q\}\cup\{(\zeta,-\zeta^{\s+1},\theta,\ldots)\in V\;|\;\zeta\in \mathbb{F}_q^\times,\theta\in \mathbb{F}_q\}$ is included in the $H$-orbit containing $\infty$ and thus $|H|=|\infty^H||H_\infty|\geq (2q^2-q+1)q(q-1)$.

When $q\geq 27$, we have $|H|=(s^3+1)s^3(s-1)$ by the above-mentioned result, where $s^m=q$ for some positive odd integer $m$. It follows that $m=1$ and $H=\R(q)$. When $q=3$, we use the permutation representation of $\R(3)$ as a primitive group of degree $28$ in the database of primitive groups in \magma \cite{magma}. Now $\R(3)$ acts on $\Delta:=\{1,2,\ldots,28\}$, and the two actions of $\R(3)$ on $V$ and $\Delta$ are permutation isomorphic. Let $J$ be the normal subgroup of $\R(3)_1$ (the stabilizer of $1\in \Delta$ in $\R(3)$) which is regular on $\Delta\setminus\{1\}$, and let $Z$ be the centre of $J$. Then $H$ is (permutation) isomorphic to $\la Z,\R(3)_{1,2}, \t\ra$ for some involution $\t \in \R(3)$. Computation in \magma shows that $\la Z,\R(3)_{1,2}, \t\ra$ has order $6,18$ or $504$ for any involution $\t$ of $\R(3)$. Hence $|H|=504$ as $|H|\geq 16\cdot 3\cdot 2=96$. Since $\R(3)'$ is the only subgroup of $\R(3)$ of order $504$, we have $H=\R(3)'$.
\qed

\medskip
\textsf{Case (iii):} $|P|=q^2$, $L_1\cup L_2 \subseteq P$
\medskip

Now $\{(0,\eta,0,\ldots) \;| \;\eta \in \mathbb{F}_q\} \subseteq P$. Since ${(0,\eta,0,\ldots)}^{t_{\a,\b,\g}}= (\a, \eta +\b, \g- \a \eta,\ldots)$, we have $\la t_{0,\b,0} \;|$ $\b \in \mathbb{F}_q\ra \leq {G}_{\infty, P}$ and $H:= \la t_{0,\b,0}, n_\kappa \;|\;\b \in \mathbb{F}_q, \kappa \in \mathbb{F}_q^{\times}\ra \leq {G}_{\infty, P}$.

\begin{lemma}
\label{lem:R(q)-1}
Let $G = \la \R(q), \zeta\ra$ (as given in \eqref{rqzeta}) be a subgroup of $\Sym(V)$ containing $\R(q)$ as a normal subgroup. Suppose that $P$ is an imprimitive block of $G_\infty$ on $V\setminus\{\infty\}$ containing $\bzero$ such that $|P|=q^2$ and $L_1\cup L_2\subseteq P$. Then the following hold:
\begin{itemize}
\item[\rm (a)] $P= \{(0,\eta_2,\eta_3,\ldots)\; |\; \eta_2, \eta_3 \in \mathbb{F}_q\}$ and $G_{\infty,P}= \la n_\kappa, t_{0,\xi,\eta}, \zeta \;|\; \xi, \eta \in \mathbb{F}_q, \kappa \in \mathbb{F}_q^{\times}\ra$;
\item[\rm (b)] setting $L := P \cup \{\infty\}$, $\DD:=(V, {L}^{\R(q)})=(V, {L}^{G})$ is a $2$-$(q^3+1, q^2+1, q^2+1)$ design admitting $G$ as a $2$-point-transitive and block-transitive group of automorphisms, and $\Om:={(\infty, L)}^{\R(q)}={(\infty, L)}^{G}$ is a 1-feasible $G$-orbit on the flag set of $\DD$; moreover, the $G$-flag graph of $\DD$ with respect to $(\Om, \Xi)$ is the same as the $\R(q)$-flag graph $\Ga(\DD,\Om,\Xi)$ with respect to $(\Om, \Xi)$;
\item[\rm (c)] $\Ga(\DD,\Om,\Xi)$ is connected with order $|\Om|=q(q^3+1)$ and valency $q^2$.
\end{itemize}
\end{lemma}

\proof
(a) Since $|P|=q^2$ by our assumption, it suffices to prove $\eta_1 = 0$ for every $(\eta_1,\eta_2,\eta_3,\ldots) \in P$. Suppose otherwise. By the action of $H$ we may assume $(1,0,\ve_0,\ldots) =: \rho \in P$ for some $\ve_0 \in \mathbb{F}_q$. Since $|{\rho}^{H}|= {|H|}/{|H_{\rho}|} =|H|= q(q-1)$ and ${\rho}^{H} \cap \{(0,\eta,0,\ldots)\;|\;\eta \in \mathbb{F}_q\} =\emptyset$, we have $P = {\rho}^{H} \cup \{(0,\eta,0,\ldots) \;|\; \eta \in \mathbb{F}_q\} = \{(\eta_1,\eta_2,\eta_3,\ldots) \in V\;| \;\eta_3= \eta_1^{\s+2} \ve_0 + \eta_1\eta_2\}$. So $(0,1,0,\ldots)$ and $(1,0,\ve_0,\ldots)$ are points in $P$. If ${(0,1,0,\ldots)}^{t_{x,y,z}}= (1,0,\ve_0,\ldots)$, then $x= 1,y=-1,z=1+\ve_0$. On the other hand, $(0,-1,0,\ldots) \in P$, and ${(0,-1,0,\ldots)}^{t_{1,-1,1+\ve_0}}= (1,-2,2+\ve_0,\ldots) = (1,1,2+\ve_0,\ldots) \notin P$. Therefore, $P$ is not an imprimitive block of ${\R(q)}_{\infty}$ on $V\setminus \{\infty\}$, a contradiction.

(b) It is clear that $\DD$ is a $2$-$(q^3+1, q^2+1, \l)$ design admitting $G$ as a $2$-point-transitive and block-transitive group of automorphisms, and $\Om$ is a 1-feasible $G$-orbit on the flag set of $\DD$. If $\l =1$, then ${\R(q)}_L$ is $2$-transitive on $L$ and thus $|{\R(q)}_L|= |L|\cdot|{\R(q)}_{\infty, L}|= (q^2+1)q^2(q-1)$. But $|{\R(q)}_L|$ should divide $|\R(q)|= (q^3+1)q^3(q-1)$, a contradiction. Therefore, by Lemma \ref{lem:2 lambda}, $\l = q^2+1$.

(c) By Lemma \ref{lem:connectedness of G-flag graph}, it suffices to prove $\R(q) = H:= \la t_{0,\xi,\eta}, w\;|\;\xi,\eta\in \mathbb{F}_q\ra$. In fact, similar to \eqref{equ:R(q) is generated by}, we can see that $H$ is transitive on $V$, and thus $|H|=|V||H_\infty|$ is divisible by $(q^3+1)q^2$. Similar to the proof of part (c) of Lemma \ref{lem:R(q)}, when $q\geq 27$, we have $|H|=(s^3+1)s^3(s-1)$, where $s^m=q$ for some odd positive integer $m$. It follows that $m=1$ and $H=\R(q)$. When $q=3$, the action of $\R(3)$ on $V$ is permutation isomorphic to the primitive action of $\R(3)$ on $\Delta:=\{1,2,\ldots,28\}$ (see \magma \cite{magma}). Let $J$ be the normal subgroup of $\R(3)_1$ which is regular on $\Delta\setminus\{1\}$. $J$ has two subgroups of order $9$ which are normal in $\R(3)_1$. One of them, say $X$, is elementary abelian, while the other is cyclic. So $H$ is (permutation) isomorphic to $\widetilde{H}:=\la X, \t\ra$ for some involution $\t \in \R(3)$. Computation in \magma shows that $|\widetilde{H}|=18$ or $1512$  for any involution $\t$ in $\R(3)$. Since $|H|\geq 28\cdot 9$, it follows that $H=\R(3)$.
\qed

\medskip
\textsf{Case (iv):} $|P|=q^2$, $(L_1 \cup L_2)\cap P=\emptyset$
\medskip

\begin{lemma}
\label{lem:R(q)-2}
Let $\soc(G)=\R(q)$ or $G=\R(3)$. Suppose that $P$ is an imprimitive block of $G_\infty$ on $V\setminus\{\infty\}$ containing $\bzero$ such that $|P|=q^2$ and $(L_1 \cup L_2)\cap P=\emptyset$. Then the following hold:
\begin{itemize}
\item[\rm (a)] $G=\R(3)$, $P=\{(a,-a^2,c,\ldots)\in V\;|\;a,c \in \mathbb{F}_3\}$, and $G_{\infty,P}= \la n_{-1}, t_{x,-x^2,z}\;|\; x,z \in \mathbb{F}_3\ra$;
\item[\rm (b)] setting $L:=P \cup \{\infty\}$, $\DD :=(V, L^{\R(3)})$ is a $2$-$(28, 10, 10)$ design admitting $\R(3)$ as a $2$-point-transitive and block-transitive group of automorphisms, and $\Om :={(\infty, L)}^{\R(3)}$ is a 1-feasible $\R(3)$-orbit on the flag set of $\DD$;
\item[\rm (c)] the $\R(3)$-flag graph $\Ga(\DD,\Om,\Xi)$ has three connected components.
\end{itemize}
\end{lemma}

\proof
(a) Since ${(0,0,0,\ldots)}^{t_{\a,\b,\g}}= (\a,\b,\g,\ldots)$, if $(\a,\b,\g,\ldots)\in P$, then $t_{\a,\b,\g}$ stabilizes $P$. Let $(\eta_1,\eta_2,\eta_3,\ldots)\in P$ be a fixed element with $\eta_1\neq0$. Since $N:=\la n_\kappa\;|\;\kappa\in \mathbb{F}_q^\times\ra$ stabilizes $P$, we may assume that $\eta_1=1$. Then $(0,0,\g,\ldots)\in P$ for any $\g\in \mathbb{F}_q^\times$ as $t_{1,\eta_2,\eta_3}^3=t_{0,0,-1}$ and $N$ stabilizes $P$. Thus, if $\b\neq0$, then $(0,\b,\g,\ldots)\notin P$, for otherwise $(0,\b,0,\ldots)=(0,\b,\g,\ldots)^{t_{0,0,-\g}}\in P$, which contradicts the assumption $(L_1 \cup L_2)\cap P=\emptyset$.

Let $(1,\xi_2,\xi_3,\ldots)$ and $(1,\theta_2,\theta_3,\ldots)$ be two points in $P$. Then $(1,\xi_2,\xi_3,\ldots)^N=(1,\theta_2,\theta_3,\ldots)^N$ holds if and only if $\xi_2=\theta_2$ and $\xi_3=\theta_3$. Moreover, $(-1,\xi_2,-\xi_3,\ldots)\in P$ and $t_{-1,\xi_2,-\xi_3}$ stabilizes $P$. Hence $(1,\eta_2,\eta_3,\ldots)^{t_{-1,\xi_2,-\xi_3}}=(0,\eta_2+\xi_2-1,\eta_3-\xi_3+\eta_2+\xi_2-1,\ldots)\in P$ and so $\eta_2+\xi_2-1 = 0$. Setting $\xi_2=\eta_2$, we obtain $\eta_2=-1$. Therefore, $P=\{(0,0,\g,\ldots)\;|\;\g\in \mathbb{F}_q\}\cup \left(\cup_{c\in \mathbb{F}_q} (1,-1,c,\ldots)^N\right) = \{(a,-a^{\s+1},c,\ldots)\in V\;|\;a,c \in \mathbb{F}_q\}$.
On the other hand, since $(a,-a^{\s+1},c,\ldots)^{t_{x,y,z}}=(a+x,-a^{\s+1}+y+x^\s a,c+z+xa^{\s+1}+ya-x^{\s+1}a,\ldots)$, if $(a_0,-a_0^{\s+1},c,\ldots)^{t_{x,y,z}}\in P$ for some $a_0\in \mathbb{F}_q$, then $y=x^\s a_0-a_0^\s x-x^{\s+1}$. Thus, $t_{x,y,z}$ stabilizes $P$ if and only if $x^\s a_0-a_0^\s x-x^{\s+1}=x^\s a-a^\s x-x^{\s+1}$ for any $a\in \mathbb{F}_q$, that is, $(a-a_0)^\s x=(a-a_0)x^\s$ for any $a\in \mathbb{F}_q$. Hence $P=\{(a,-a^{\s+1},c,\ldots)\in V\;|\;a,c \in \mathbb{F}_q \}$ is an imprimitive block of ${\R(q)}_{\infty}$ on $V\setminus \{\infty\}$ if and only if $q=3$.

(b) $\DD$ is a $2$-$(28, 10, \l)$ design for some $\l \ge 1$. If $\l =1$, then ${\R(3)}_L$ is $2$-transitive on $L$ and thus $|{\R(3)}_L|= |L|\cdot|{\R(3)}_{\infty, L}|$. But $|{\R(3)}_L|$ should divide $|\R(3)|= 28\cdot 54$, a contradiction. Therefore, by Lemma \ref{lem:2 lambda}, $\l = 10$.

(c) Set $H:=\la n_{-1}, t_{x,-x^2,z}, w\;|\; x,z \in \mathbb{F}_3\ra$. Then $|H|=|\infty^H||H_\infty|\geq 16\cdot 3\cdot 2=96$ by \eqref{equ:R(q) is generated by}. Similar to the proof of part (c) of Lemma \ref{lem:R(q)-1}, in view of the action of $\R(3)$ on $\Delta := \{1,2,\ldots,28\}$, let $J$ be the normal subgroup of $\R(3)_1$ which is regular on $\Delta\setminus\{1\}$. Then $H$ is (permutation) isomorphic to $\widetilde{H} := \la Y, \R(3)_{1,2}, \t\ra$ for some involution $\t \in \R(3)$, where $Y$ is the cyclic subgroup of $J$ of order $9$ which is normal in $\R(3)_1$. Using \magma, we find that $|\widetilde{H}|=18$ or $504$ for any involution $\t$ of $\R(3)$. Therefore, $H=\R(3)'$ and $\Ga(\DD,\Om,\Xi)$ has three connected components by Lemma \ref{lem:connectedness of G-flag graph}.
\qed

\section{Affine case}
\label{sec:affine case}

In this section we assume that $G$ is a $2$-transitive permutation group acting on a set $V$ which we always assume to be a suitable vector space over some finite field, and $\soc(G)$ is abelian. Let $u:=|V|= p^d$ be the degree of this permutation representation, where $p$ is a prime and $d \geq 1$. It is known (see \cite{Kantor85}, \cite[p.194]{Cameron99}, \cite{Cameron81}, \cite[p.386]{Huppert82}) that $u$ and the stabilizer $G_{\bf 0}$ of the zero vector ${\bf 0}$ in $G$ are as follows.

\begin{description}
\item[\rm (i)] $G_{\bf 0} \leq \GammaL(1, q)$, $q=p^d$;

\item[\rm (ii)] $G_{\bf 0} \unrhd \SL(n, q)$, $n \geq 2$, $q^n = p^d$;

\item[\rm (iii)] $G_{\bf 0} \unrhd \Sp(n, q)$, $n \geq 4$, $n$ is even, $q^n = p^d$;

\item[\rm (iv)]  $G_{\bf 0} \unrhd G_2(q)$, $q^6 = p^d$, $q >2$, $q$ is even;

\item[\rm (v)] $G_{\bf 0} = G_2(2)' \cong \PSU(3, 3)$, $u = 2^6$;

\item[\rm (vi)] $G_{\bf 0} \cong A_6$ or $A_7$, $u = 2^4$;

\item[\rm (vii)] $G_{\bf 0} \cong \SL(2,13)$, $u = 3^6$;

\item[\rm (viii)] $G_{\bf 0} \unrhd \SL(2, 5)$ or $G_{\bf 0} \unrhd \SL(2, 3)$, $d=2$, $p = 5, 7, 11, 19, 23, 29$, or $59$;

\item[\rm (ix)] $d=4$, $p=3$. $G_{\bf 0} \unrhd \SL(2, 5)$ or $G_{\bf 0} \unrhd E$, where $E$ is an extraspecial group of order $32$.
\end{description}

Since $G$ is $2$-transitive on $V$, in each case $G_{\bf 0}$ is transitive on $V \setminus \{{\bf 0}\}$. Denote by $T$ the subgroup of $\Sym(V)$ consisting of all translations of $V$, so that $G = T \rtimes G_{\bf 0}$.

\subsection{$G_{\bf 0} \leq \GammaL(1, q)$, $q=p^d$}
\label{subsec:AGammaL(1, q)}

This case gives (f) in Table \ref{tab0a}. Now $G$ acts on $V=\mathbb{F}_q$, and a typical element of $G$ is of the form
$$
\t({a, c, \vp}): \mathbb{F}_q \rightarrow \mathbb{F}_q, z \mapsto az^{\vp} + c,
$$
where $a \in \mathbb{F}_q^{\times}$, $c \in \mathbb{F}_q$ and $\vp \in \Aut(\mathbb{F}_q) = \la \zeta\ra$, with
$$
\zeta: \mathbb{F}_q \rightarrow \mathbb{F}_q,\; z \mapsto z^p
$$
the Frobenius map. Denote
$$
t{(a, j)} := \t(a, 0, \zeta^j),
$$
where $j$ is an integer. Given $\d \in \Aut(\mathbb{F}_q)$ and integer $i \ge 0$, we use $[\d,i]$ to denote $({p^{ni}-1})/({p^n-1})$ and $\d-1$ to denote $p^n-1$, where $n$ is the smallest positive integer such that $\d=\zeta^n$. Thus, for $i>0$ and $x\in \mathbb{F}_q^\times$, $x^{[\d,i]}$ is the product of $x^{\d^{i-1}}$, $x^{\d^{i-2}}$, $\ldots$, $x^\d$, $x$ in $\mathbb{F}_q^{\times}$. The following two lemmas can be easily proved.

\begin{lemma}
\label{lem:12 sequence}
Let $H$ be a subgroup of $\mathbb{F}_q^{\times}$. Let $x \in \mathbb{F}_q^{\times} \setminus H$ and $\d \in \Aut(\mathbb{F}_q)$. In the sequence: $H$, $Hx^{[\d,1]}$, $Hx^{[\d,2]}$, $\ldots$, $Hx^{[\d,n]}$, $\ldots$, if $j$ is the smallest positive integer such that $Hx^{[\d,j]}$ equals a previous item, then $Hx^{[\d,j]} = H$.
\end{lemma}

\begin{lemma}
\label{lem:13 num6}
Let $H$ be a subgroup of $\mathbb{F}_q^{\times}$ and $x \in \mathbb{F}_q^{\times}$. Then $x \in H$ if and only if $x^{|H|} = 1$.
\end{lemma}

Let $a$ and $c$ be coprime positive integers. Denote by $\ord_p(a)$ the exponent of a prime $p$ in $a$, namely the largest nonnegative integer $i$ such that $p^i\mid a$. If $j$ is the smallest positive integer such that $c \mid (a^j-1)$, then we say that $a$ has order $j$ ($\bmod \;c$). Denote by $S(n)$ the set of prime divisors of a positive integer $n$.

\begin{lemma}
\label{lem:14 num7}
Let $m>1$ and $a>1$ be coprime integers.
\begin{itemize}
\item[\rm (a)] If $m$ is odd, then $a$ has order $m\;(\bmod \;(a-1)m)$ if and only if $S(m)\subseteq S(a-1)$;
\item[\rm (b)] if $m$ is even, then $a$ has order $m\;(\bmod \;(a-1)m)$ if and only if $S(m)\subseteq S(a-1)$, and either $a\equiv 1\;(\bmod \;4)$ or $4\nmid m$.
\end{itemize}
\end{lemma}

\proof
Let $p$ be a prime not dividing $a$ and suppose $a$ has order $f$ ($\bmod \;p$). Suppose that $n$ is a positive integer and $f\mid n$. Then by \cite[pp.355--356, 1), 2), 3), 4)]{Artin}, we have $\ord_p(a^{n}-1) =\ord_p(a^f-1) +\ord_p(n)$ when $p$ or $n$ is odd, and $\ord_2(a^{n}-1) =\ord_2(a-1) +\ord_2(a+1) +\ord_2(n)-1$ when $p=2$ and $n$ is even.

First assume that $a$ has order $m$ ($\bmod \;(a-1)m$). Suppose $S(m)\nsubseteq S(a-1)$ and let $p$ be the largest one in $S(m)\setminus S(a-1)$. Then $p$ is odd and $p\nmid a$ as $\gcd(a, m)=1$ by our assumption. Suppose that $a$ has order $f$ ($\bmod \;p$). Then $f \mid m$ as $p$ divides $a^m-1$. Set $m_1=m/p$. Let $\xi\in S((a-1)m)=S(a-1)\cup S(m)$ and suppose that $a$ has order $f_\xi$ ($\bmod \;\xi$). If $\xi\in S(m)\setminus S(a-1)$, then $\xi$ is odd, $f_\xi\mid m_1$ since $f_\xi<\xi\leq p$, and $\ord_{\xi}(a^{m_1}-1) =\ord_{\xi}(a^{f_\xi}-1) +\ord_{\xi}(m_1)\geq \ord_{\xi}(m)=\ord_{\xi}((a-1)m)$. If $\xi\in S(a-1)\setminus S(m)$, then $\ord_\xi(a^{m_1}-1)\geq\ord_\xi(a-1)=\ord_\xi((a-1)m)$. If $\xi\in S(m)\cap S(a-1)$, then $\ord_\xi(a^{m_1}-1)\geq\ord_\xi(a-1)+\ord_\xi(m_1)=\ord_\xi(a-1)+\ord_\xi(m)=\ord_\xi((a-1)m)$. It follows that $(a-1)m\mid (a^{m_1}-1)$, a contradiction. Therefore $S(m)\subseteq S(a-1)$. If in addition $m=2\ell$, we need to prove either $a\equiv 1\;(\bmod\;4)$ or $4\nmid m$. Suppose $a\equiv -1\;(\bmod\;4)$ and $4\mid m$, then $\ord_\xi(a^{\ell}-1)=\ord_\xi(a-1)+\ord_\xi(\ell)=\ord_\xi(a-1)+\ord_\xi(m)$ when $\xi\in S(m)$ is odd, and
$$
\ord_2(a^{\ell}-1)=\ord_2(a-1)+\ord_2(a+1)+\ord_2(\ell)-1\geq\ord_2(a-1)+\ord_2(m),
$$
which implies $(a-1)m\mid (a^{\ell}-1)$, a contradiction.

Next assume that $S(m)\subseteq S(a-1)$. If $m$ is even, we also assume $4\nmid m$ or $a\equiv 1\;(\bmod \;4)$. Apparently $(a-1)m\mid(a^m-1)$. Let $t$ be a positive integer such that $(a-1)m\mid(a^t-1)$. Then for any odd $p\in S(m)$, we have $\ord_p(t)\geq\ord_p(m)$. Hence if $m$ is odd, then $m\leq t$. Suppose that $m$ is even. Then $t$ has to be even (otherwise $\ord_2(a^t-1)=\ord_2(a-1)$ and $a^t-1$ can not be divided by $(a-1)m$). If $4\nmid m$, then obviously $m\leq t$. If $a\equiv 1\;(\bmod \;4)$, then $\ord_2(a^t-1)=\ord_2(a-1)+\ord_2(t)\geq \ord_2((a-1)m)$ and thus $\ord_2(t)\geq \ord_2(m)$, which implies $m\leq t$. Therefore, $a$ has order $m$ ($\bmod \;(a-1)m$).
\qed

\smallskip
Next we prove some results on the structure of $G_0$ and determine all possible imprimitive blocks of $G_0$ on $\mathbb{F}_q^\times$ containing $1$.

Denote by $s$ the smallest positive integer such that $t(a, s) \in G_0$ for some $a \in \mathbb{F}_q^{\times}$. Then $s$ must be a divisor of $d$ and $\{\ell>0\;|\; t({a,\ell}) \in G_0\; \text{for some}\;a\in\mathbb{F}_q^{\times}\} = \{js\;|\; j= 1, 2, \ldots\}$, since $t({a_1,{j_1}})t({a_2,{j_2}})= t({a_2{a_1^{\zeta^{j_2}}},{j_1+ j_2}}) \in G_0$ and $t({a_i,{j_i}})^{-1} = t((1/a_i)^{\zeta^{-j_i}},{-j_i}) \in G_0$ for $t({a_i,{j_i}})\in G_0$, $i=1,2$. If $s=d$, then $G_0 \leq \GL(1, q)$ and so $G_0= \GL(1, q)$ as $G_0$ is transitive on $\mathbb{F}_q^\times$. For any integer $i$, set
\begin{align*}
A_i := \{t({a,is})\;|\;t({a, is}) \in G_0\},\; H_i := \{a\;|\;t({a,is}) \in G_0\}.
\end{align*}
Then $A_0$ is a normal cyclic subgroup of $G_0$, $H:=H_0$ is a cyclic subgroup of $\mathbb{F}_q^{\times}$, and $A_i=A_j$ if and only if $d\mid (i-j)s$. Since $A_i t{(x,s)}^j\subseteq A_{i+j}$ for any $i, j$ and $t(x,s) \in A_1$, $|A_i|$ is a constant and thus $A_i t{(x,s)}^j=A_{i+j}$. Setting
$$
\vp:=\zeta^s,
$$
we then have
$$
A_i = A_0 t{(x,s)}^i,\;H_i = Hx^{[\vp,i]},\; i=1,2,\ldots,d/s-1,
$$
and $A_{{d}/{s}} =A_0 t({x,s})^{{d}/{s}} = A_0$ and $H_{{d}/{s}} = Hx^{[\vp,d/s]} = H$. Since
$$
G_0=A_0\cup A_1\cup\cdots\cup A_{{d}/{s}-1}\; \mbox{(disjoint union)}
$$
is transitive on $\mathbb{F}_q^{\times}$, we have
$$
\mathbb{F}_q^{\times}=H\cup H_1\cup H_2\cup\cdots\cup H_{{d}/{s}-1}.
$$

Denote by $m$ the smallest positive integer such that $t({1, ms}) \in G_{0,1}$. Then $m\leq d/s$, $G_{0,1} =\la t({1,ms})\ra$ and $|G_{0,1}| ={d}/({ms})$. Let $x\in H_1$. Since $1\in H_m=Hx^{[\vp,m]}$, we have $Hx^{[\vp,m]} = H$. In the case when $m>1$, if $Hx^{[\vp,j]} = H$ for some $j < m$, then $H_j = H$ and $t({1,js}) \in A_j \subseteq G_0$, which contradicts the definition of $m$. Hence by Lemma \ref{lem:12 sequence}, in the sequence: $H$, $Hx^{[\vp,1]}$, $Hx^{[\vp,2]}$, $\ldots$, $Hx^{[\vp,m-1]}$, $Hx^{[\vp,m]}$, $\ldots$, the first $m$ items are pairwise distinct, and the subsequent items repeat the previous ones. Since $G_0$ is transitive on $\mathbb{F}_q^{\times}$, we have
\be
\label{equ:5 r H}
\mathbb{F}_q^{\times} = H \cup Hx^{[\vp,1]} \cup \cdots \cup Hx^{[\vp,m-1]},\;
|\mathbb{F}_q^{\times}:H| = m,\; m \mid [\vp,m].
\ee

\begin{lemma}
\label{lem:structure of G0}
Let $G \leq \AGammaL(1, q)$ act $2$-transitively on $\mathbb{F}_q$, where $q = p^d$ with $p$ a prime, and let $s$ and $m$ be as above.
\begin{itemize}
\item[\rm (a)] If $m=1$, then $G_0$ is the group generated by $\GL(1, q)$ and $\t({1,0,\zeta^s})$. Conversely, if $G_0$ contains $\GL(1, q)$, then $m=1$;
\item[\rm (b)] if $m>1$, then $p^s$ has order $m\;(\bmod \;m(p^s-1))$, $S(m) \subseteq S(p^s-1)$, and $G_0$ has $\phi(m)$ possibilities, where $\phi(m)$ is the number of positive integers less than $m$ and coprime to $m$.
\end{itemize}
\end{lemma}

\proof
Let $A_i, H_i$ and $H$ be defined as above, and $\vp:=\zeta^s$.

(a) If $m=1$, then $H=\mathbb{F}_q^{\times}$ and $G_0$ is generated by $\GL(1, q)$ and $\t({1,0,\zeta^s})$. Conversely, if $G_0$ contains $\GL(1, q)$, then $m=1$ by the definition of $m$.

(b) Assume $m>1$. Since $G_0$ is transitive on $\mathbb{F}_q^\times$, we have $Hx \neq H$, $Hx^{[\vp,2]}\neq H$, $\ldots$, $Hx^{[\vp,m-1]} \neq H$ by \eqref{equ:5 r H} and Lemma \ref{lem:12 sequence}, where $x\in H_1$. By Lemma \ref{lem:13 num6}, this is equivalent to saying $|H|=({q-1})/{m}$ and $x^{|H|} \neq 1$, $x^{[\vp,2]|H|}\neq 1$, $\ldots$, $x^{[\vp,m-1]|H|}\neq 1$. Denote the set of solutions in $\mathbb{F}_q^{\times}$ of each of the equations
\begin{align*}
\a^{|H|} = 1, \a^{[\vp,2]|H|}= 1, \ldots, \a^{[\vp,m-1]|H|} = 1
\end{align*}
by $E_1$, $E_2$, $\ldots$, $E_{m-1}$, respectively. Then $E_i$ ($i=1, 2, \ldots, m-1$) is a cyclic subgroup of $\mathbb{F}_q^{\times}$ with $|E_i| = \gcd(q-1, [\vp,i]|H|)=|H|\cdot \gcd(m,[\vp,i])$, and $E_i/H$ is a subgroup of $\mathbb{F}_q^\times/H$ of order $\gcd(m,[\vp,i])$. Hence the existence of $x$ satisfying \eqref{equ:5 r H} implies that
$m \mid (q-1)$ and $\mathbb{F}_q^{\times}\neq\cup^{m-1}_{i=1}E_i$ (that is, $\rho \notin\cup^{m-1}_{i=1}E_i$, where $\rho$ is a generator of $\mathbb{F}_q^{\times}$), or equivalently $m \nmid ({p^{si}-1})/({p^s-1})$, $i=1, 2, \ldots, m-1$, and $m \mid ({p^{sm}-1})/({p^s-1})$. Thus $p^s$ has order $m\;(\bmod \;m(p^s-1))$, and $S(m) \subseteq S(p^s-1)$ by Lemma \ref{lem:14 num7}.

Moreover, by \eqref{equ:5 r H}, $\{[\vp,1], \ldots, [\vp,m]\}$ is a complete residue system modulo $m$. It follows that $\cup^{m-1}_{i=1}(E_i/H)$ is the set of all non-generators of $\mathbb{F}_q^\times/H$. Let $\xi$ be a fixed generator of $\mathbb{F}_q^\times$. Then $\mathbb{F}_q^{\times} \setminus\cup^{m-1}_{i=1}E_i = \cup_{i=1}^{\phi(m)} H\xi^{\ell_i}$, where $\{\ell_1=1,\ell_2,\ldots,\ell_{\phi(m)}\}$ is a reduced residue system modulo $m$, and $G_0$ is the group generated by $\{t(a,0)\; |\; a \in H\}$ and $t(\xi^{\ell_i},s)$  for some $i\in\{1, 2, \ldots, \phi(m)\}$.
\qed

\begin{lemma}
\label{lem:11 blocks of G0}
Let $G \leq \AGammaL(1, q)$ act $2$-transitively on $\mathbb{F}_q$, where $q = p^d$ with $p$ a prime. Let $s, m, \vp, H, A_i$ and $H_i$ be as above. A subset $P$ of $\mathbb{F}_q^{\times}$ is an imprimitive block of $G_0$ on $\mathbb{F}_q^{\times}$ containing $1$ if and only if it is of the form
\be
\label{eq:GammaL}
P = K \cup Kw^{[\psi,1]}\cup Kw^{[\psi,2]} \cup \cdots\cup Kw^{[\psi,m/e-1]},
\ee
where $e$ is a divisor of $m$, $\psi:=\zeta^{es}$, and $K$ is a subgroup of $H$ such that $w^{[\psi,m/e]}\in K$ for some $t(w,es)\in A_e$.
\end{lemma}

\proof
Suppose that $P$ is an imprimitive block of $G_0$ on $\mathbb{F}_q^\times$ containing $1$. Then $G_{0,1} \leq G_{0, P} \leq G_0 \leq \GammaL(1, q)$. Denote by $e$ the smallest positive integer such that $t({c, es})\in G_{0, P}$ for some $c \in\mathbb{F}_q^{\times}$. Then $G_{0, P}\subseteq A_0\cup A_e\cup A_{2e}\cup A_{3e}\cup \cdots$. For any integer $i$, set
$$
C_i := \{t(a,ies)\;|\;t({a,ies}) \in G_{0,P}\},\;K_i := \{a\;|\;t({a,ies}) \in G_{0,P}\}.
$$
Then $K:=K_0$ is a subgroup of $H$. Set $\psi:=\vp^e = \zeta^{es}$. For any $t({w,es}) \in G_{0, P}$, we have
\begin{equation}
\label{equ:Hb Hw}
A_{je} = A_0 t(w,es)^j,\;
H_{je} =H w^{[\psi,j]},\;
C_{j} = C_0 t(w,es)^j,\;
K_{j} =K w^{[\psi,j]},\;
j=1,2,\ldots,{m}/{e}-1.
\end{equation}
Let $\ell$ be the smallest positive integer such that $Kw^{[\psi,\ell]} = K$. Then $t({1,\ell es})\in G_{0,1}$. Since $G_{0,1} \leq G_{0, P}$, by the definition of $m$, we have $m = e \ell$ and so $e$ divides $m$. By Lemma \ref{lem:12 sequence}, in the sequence: $K$, $Kw^{[\psi,1]}$, $Kw^{[\psi,2]}$, $\ldots$, $Kw^{[\psi,m/e-1]}$, $Kw^{[\psi,m/e]}$, $\ldots$, the first ${m}/{e}$ items are pairwise distinct, and the subsequent items repeat the previous ones. Since $G_{0, P}$ is transitive on $P$, $P$ must be of the form \eqref{eq:GammaL} and moreover $Kw^{[\psi,m/e]} = K$.

We now prove that any subset $P$ of $\mathbb{F}_q^\times$ defined in \eqref{eq:GammaL} is an imprimitive block of $G_0$ on $\mathbb{F}_q^{\times}$. Let $e, \psi$ and $K$ be as stated in the lemma. By Lemma \ref{lem:12 sequence} and the definition of $m$, in the sequence: $K$, $Kw^{[\psi,1]}$, $Kw^{[\psi,2]}$, $\ldots$, $Kw^{[\psi,m/e-1]}$, $Kw^{[\psi,m/e]}$, $\ldots$, the first ${m}/{e}$ terms must be pairwise distinct, and the subsequent terms repeat the previous ones as $w^{[\psi,m/e]}\in K$. For any $t(z,ns) \in G_0$ ($n>0$), if $e\nmid n$, then $K$, $Kw^{[\psi,1]}$, $Kw^{[\psi,2]}$, $\ldots$, $Kw^{[\psi,m/e-1]}$ are mapped into $H_n$, $H_{e+n}$, $H_{2e+n}$, $\ldots$, $H_{m-e+n}$ respectively by $t(z,ns)$, and hence $P\cap P^{t(z,ns)}=\emptyset$. If $e\mid n$, say $n=ej$, then $Kw^{[\psi,i]}$ is mapped into $H_{e(i +j)}$ by $t(z,ns)$. Since $t(z,ns)\in A_n=A_0 t(w,es)^{j}$, we have $z\in H_n=Hw^{[\psi,j]}$ and $z=yw^{[\psi,j]}$ for some $y\in H$, and $(Kw^{[\psi,i]})^{t(z,ns)}=yw^{[\psi,j]}K(w^{[\psi,i]})^{\psi^{j}}=Kyw^{[\psi,i+j]}$. Hence $(Kw^{[\psi,i]})^{t(z,ns)}=Kw^{[\psi,i+j]}$ if $y\in K$, and $(Kw^{[\psi,i]})^{t(z,ns)}\cap Kw^{[\psi,i]}=\emptyset$ if $y\notin K$. Therefore, $y\in K$ implies $P^{t(z,ns)}=P$, and $y\notin K$ implies $P\cap P^{t(z,ns)}=\emptyset$. Thus $P$ is an imprimitive block of $G_0$ on $\mathbb{F}_q^{\times}$.
\qed

\begin{lemma}
\label{prop:P an impri block, L(1,q)}
Let $G, s, m, \vp, A_i, H_i$ and $H$ be as in Lemma \ref{lem:11 blocks of G0}, and by the proof of Lemma \ref{lem:structure of G0} we may assume $H_1=H\xi^\ell$, where $\xi$ is a fixed generator of $\mathbb{F}_q^\times$ and $\ell$ is a positive integer coprime to $m$. Then a subset $P$ of $\mathbb{F}_q^\times$ containing $1$ is an imprimitive block of $G_0$ on $\mathbb{F}_q^\times$ if and only if either $P$ is a subgroup of $\mathbb{F}_q^\times$, or there exist $j=nm+\ell [\vp,e]$ and $K\leq H$ with $|H/K|$ dividing $j[\psi,m/e]/{m}$ such that $P$ is given by \eqref{eq:GammaL}, where $n\geq0$, $e$ is a divisor of $m$, $\psi:=\zeta^{es}$ and $w:=\xi^j$.
\end{lemma}

\proof
By Lemma \ref{lem:11 blocks of G0}, we may assume $P$ is given by \eqref{eq:GammaL}. If $m=1$, then $e=m=1$, $w=w^{[\psi,m/e]}\in K$, and thus $P=K$ is a subgroup of $\mathbb{F}_q^{\times}$.

Assume $m>1$ in the remaining proof. The subgroup $K$ of $H$ in Lemma \ref{lem:11 blocks of G0} must allow the existence of $w\in H_e$ such that $w^{[\psi,m/e]}\in K$. We may assume $w=\xi^{nm+\ell [\vp,e]}$ for some $n\geq0$ as $H_e=H\xi^{\ell [\vp,e]}$ and $H=\la \xi^m\ra$. Set $j=nm+\ell [\vp,e]$. By Lemma \ref{lem:14 num7}, $e\mid [\vp,e]$ and $m/e$ divides $[\psi,m/e]$. Thus $e\mid j$ and $m\mid j[\psi,m/e]$, and
\be
\label{equ:w[psi, ] in K}
w^{[\psi,m/e]}\in K\Leftrightarrow\xi^{j[\psi,m/e]|K|}=1\Leftrightarrow(q-1)\mid j[\psi,{m}/{e}]|K|\Leftrightarrow\frac{|H|}{|K|}\mid \frac{j[\psi,m/e]}{m}.
\ee
Thus $P$ in \eqref{eq:GammaL} is an imprimitive block of $G_0$ on $\mathbb{F}_q^\times$ if and only if $K$ is a subgroup of $H$ such that $|H/K|$ divides $j[\psi,m/e]/{m}$.
\qed

\begin{remark}
\label{example:P not necessarily a subgroup}
{\em We show that $P$ in Lemma \ref{lem:11 blocks of G0} needs not be a subgroup of $\mathbb{F}_q^\times$. Let $\xi, \ell, j$ and $w$ be as in Lemma \ref{prop:P an impri block, L(1,q)}. Let $K\leq H$ be such that $|H/K|$ divides $j[\psi,m/e]/{m}$. Then $P$ given by \eqref{eq:GammaL} is an imprimitive block of $G_0$ on $\mathbb{F}_q^\times$. Moreover, $P$ is a subgroup of $\mathbb{F}_q^\times$ if and only if $(Kw)^{m/e}=K$, that is, if and only if $|H/K|$ divides $j/e$.

Fix $s,m$ and $e$. Let $\mu$ be an odd prime divisor of $[\psi,m/e]/(m/e)$, and let $t:=\ord_\mu(j/e)\geq0$. Suppose that $p$ has order $f$ ($\bmod \;\mu$) so that $\ord_\mu(p^d-1)=\ord_\mu(p^f-1)+\ord_\mu(d)$. Choose $d$ such that $\ord_\mu((p^d-1)/m)\geq t+1$. Let $K$ be the subgroup of $H$ of index $\mu^{t+1}$ in $H$. Now $\mu^{t+1}\mid j[\psi,m/e]/m$, while $j/e$ cannot be divided by $\mu^{t+1}$. So $P$ is an imprimitive block of $G_0$ on $\mathbb{F}_q^\times$, but not a subgroup of $\mathbb{F}_q^\times$.
}
\end{remark}

In the rest of this section, let $G, H, K, P, m, s, \vp, e$ be as in Lemma \ref{lem:11 blocks of G0}. Set $L:=P\cup\{0\}$ and
$$
\DD:=(\mathbb{F}_q, L^G),\; \Om:=(0, L)^G.
$$
Then $\DD$ is a $2$-$(q, |P|+1, \l)$ design admitting $G$ as a $2$-point-transitive and block-transitive  group of automorphisms, and $\Om$ is a 1-feasible $G$-orbit on the flag set of $\DD$. Next we consider the connectedness of the $G$-flag graphs.

If $\l=1$, then by \cite[Proposition 4.1]{Kantor85}, $L$ is a subfield of $\mathbb{F}_q$. Conversely, if $L$ is a subfield of $\mathbb{F}_q$, then the elements in $G$ interchanging $0$ and $1$ must stabilize $L$ and thus $\l=1$. In the case when $L$ is a subfield of $\mathbb{F}_q$ with $|L|=p^t$, we have $\Ga(\DD,\Om,\Xi)\cong (p^d-1)p^d/((p^t-1)p^t)\cdot K_{p^t}$, yielding the first possibility in (f) in Table \ref{tab0a}.

In what follows we assume that $L$ is not a subfield of $\mathbb{F}_q$ so that $\l > 1$. Denote by $\widetilde{T}_1$ and $\widetilde{T}_2$ the subgroups of the addition group $(\mathbb{F}_q, +)$ generated by $\{a-1\;|\;a\in P\}$ and $\{a\;|\;a\in P\}$, respectively. Define
\begin{align*}
T_1:=\{\t(1,c,\id)\;|\;c\in\widetilde{T}_1\},\; T_2:=\{\t(1,c,\id)\;|\;c\in\widetilde{T}_2\}.
\end{align*}
\noindent Since $|H|=(q-1)/{m}=\frac{q-1}{\vp^m-1}\frac{[\vp,m]}{m}(\vp-1)$ is even when $p$ is odd, we have $\t(-1,0,\id)\in G_0$ and $\t(-1,1,\id)\in G$. Set $-P:=\{-a\,|\,a\in P\}$. We observe that for $p>2$, $P = -P$ (or equivalently $\t(-1,0,\id)$ stabilizes $P$) if and only if $|K|$ is even.

\begin{lemma}
\label{lem:generated by G0P and an elt intechanging 0 and 1}
Let $J:=\la G_{0,P}, \t(-1,1,\id)\ra$. Then $J=(T_1\rtimes G_{0,P})\rtimes \la\t(-1,1,\id)\ra$ if $P\neq -P$, and $J=T_2\rtimes G_{0,P}$ if $P=-P$.
\end{lemma}
\proof
Let $\kappa:=\t(-1,1,\id)$. It is not difficult to verify that $G_{0,P}$ normalizes $T_i$, $i=1,2$. If $\s=\t(a,0,\theta)\in G_{0,P}$, then $\kappa\s^{-1}\kappa\s=\t(1,a-1,\id)$. Hence $T_1\leq J$.

Every element $\s_1\kappa\s_2\cdots\s_{n-1}\kappa\s_n$ of $J$ is in $J_1$ or $J_1 \kappa$, where $J_1:=T_1\rtimes G_{0,P}$ and $\s_i\in G_{0,P}$, $i=1,2,\ldots,n$. We can see that $\kappa$ normalizes $J_1$, and when $P\neq -P$ we have $\kappa\notin J_1$ as $\t(-1,0,\id)\notin G_{0,P}$. When $P=-P$, we have $\t(-1,0,\id)\in G_{0,P}$ and thus $T_2\leq J$.
\qed

\begin{lemma}
\label{lem:con-G-flag-G0}
Suppose that $1 < |P| < q-1$ and $|T_i|=p^{c_i}$, $i=1,2$. Then $\Ga(\DD,\Om,\Xi)$ has $|H/K| ep^{d-c_1}/2$ connected components when $P\neq -P$, and $|H/K| ep^{d-c_2}$ connected components when $P=-P$. Therefore, $\Ga(\DD,\Om,\Xi)$ is connected if and only if $p\equiv-1\;(\bmod\;4)$, $d$ is odd, and $P$ is the subgroup of $\mathbb{F}_q^\times$ of index $2$. Moreover, if $\Ga(\DD,\Om,\Xi)$ is connected, then it has order $|\Om|=2q$ and valency $(q-1)/2$.
\end{lemma}

\proof
By Lemma \ref{lem:connectedness of G-flag graph}, the number of connected components of $\Ga(\DD,\Om,\Xi)$ is equal to
$$
|G:J|=\frac{|G|}{2^{2-i}p^{c_i}|G_{0,P}|}=\frac{|G:G_{0,P}|}{2^{2-i}p^{c_i}}=\frac{p^d(p^d-1)}{2^{2-i}p^{c_i}|K|m/e}=\frac{
|H/K|ep^{d-c_i}}{2^{2-i}},
$$
where $i=1$ if $P\neq -P$ and $i=2$ if $P=-P$.

If $P=-P$, then $\Ga(\DD,\Om,\Xi)$ must be disconnected, for otherwise we would have $c_2=d$, $K=H$, $e=1$ and so $P=\mathbb{F}_q^\times$, a contradiction.

Assume that $P\neq -P$ and $\Ga(\DD,\Om,\Xi)$ is connected. Then $p>2$, $|K|$ is odd, $c_1=d$ and $|H/K|\cdot e=2$. If $|H/K|=1$ and $e=2$, then $|K|=|H|=(q-1)/m=\frac{q-1}{\vp^m-1}\frac{[\vp,m]}{m}(\vp-1)$ is even, a contradiction. So we have $|H/K|=2$ and $e=1$. By Lemma \ref{prop:P an impri block, L(1,q)}, $j[\vp,m]/m$ is even, where $j=nm+\ell$ for some $n\geq0$ and $\ell$ coprime to $m$. Since $\ord_2(q-1)=\ord_2(|H|)+\ord_2(m)=1+\ord_2(m)$ and $q-1 = \frac{q-1}{\vp^m-1}{[\vp,m]}(\vp-1)$, we have $2\nmid [\vp,m]/m$, $2\mid j$ and $m$ is odd as $\gcd(m,\ell)=1$. It follows that $\ord_2(q-1)=1$, $p\equiv-1\;(\bmod\;4)$ and $d$ is odd. By Remark \ref{example:P not necessarily a subgroup}, $P$ is a subgroup of $\mathbb{F}_q^\times$ with index $(q-1)/(|K|m/e)=2$.

Conversely, if $p\equiv-1\;(\bmod\;4)$, $d$ is odd and $P$ is the subgroup of $\mathbb{F}_q^\times$ of index $2$, then $|P|$ is odd and $P\neq -P$. If $\widetilde{T}_1 \neq \mathbb{F}_q$, then $|\mathbb{F}_q:\widetilde{T}_1| = q/|\widetilde{T}_1|\geq p$. Hence $2=(q-1)/|P|=(q-1)/|P-1|\geq (q-1)/|\widetilde{T}_1|\geq p(1-1/q)$ as $P-1\subseteq \widetilde{T}_1$, which implies $p=3$, $d=1$ and thus $|P|=1$. Therefore, if $|P|>1$, then $\widetilde{T}_1 = \mathbb{F}_q$ and consequently $\Ga(\DD,\Om,\Xi)$ is connected.
\qed

\subsection{$G_{\bf 0} \unrhd \Sp(n, q)$, $n \geq 4$, $n$ is even, $u =q^n = p^d$}
\label{subsec:G0,Sp(n,q),n>2}

This case contributes to (g) in Table \ref{tab0a}. Denote the underlying symplectic space by $(V, f)$, where $V=\mathbb{F}_q^n$ and $f$ is a symplectic form. Let $\bfe_1:=(1,0,\ldots,0)$ and $H :=\Sp(n, q)\unlhd G_{\bf 0}$. Define $C_\a:=\{{\bf z} \in V \setminus \la{\bfe_1}\ra \;|\;f({\bf z}, {\bfe_1})= \a\}$, $\a \in \mathbb{F}_q$. By Witt's Lemma each $C_\a$ is an $H_{\bfe_1}$-orbit on $V \setminus \la{\bfe_1}\ra$. Moreover, $|C_0|= q^{n-1}-q$ and $|C_\a|= q^{n-1}$ for $\a \in \mathbb{F}_q^\times$.

A typical element of $G\leq \AGammaL(n, q)$ is $\t(A, \bc, \vp):\bx\mapsto\bx^\vp A+\bc$, where $A\in \GL(n, q)$, $\bc\in V$ and $\vp$ is a field automorphism of $\mathbb{F}_q$ acting componentwise on $\bx$. If $\t(A, \bc, \vp)\in G_{\la \bfe_1\ra}$, then $\bc\in{\la \bfe_1\ra}$ and $\row_1(A)\in {\la \bfe_1\ra}$, where $\row_1(A)$ is the first row of $A$. With respect to the bijection $\rho:\la\bfe_1\ra\rightarrow\mathbb{F}_q$, $a\bfe_1\mapsto a$, each $\t(A, \bc, \vp)\in G_{\la \bfe_1\ra}$ induces a permutation $\widetilde{\t}(a, c, \vp):\mathbb{F}_q\rightarrow\mathbb{F}_q$, $z\mapsto az^\vp + c$, where $a\bfe_1=\row_1(A)$ and $c\bfe_1=\bc$. Define $\eta:G_{\la\bfe_1\ra}\rightarrow\AGammaL(1, q)$, $\t(A, \bc, \vp)\mapsto\widetilde{\t}(a, c, \vp)$. Then $\eta$ is a homomorphism.

\begin{lemma}
\label{lem:Sp(n,q)}
Suppose that $P$ is a nontrivial imprimitive block of $G_{\bf 0}$ on $V \setminus \{{\bf 0}\}$ containing $\bfe_1$. Then
\begin{itemize}
\item[\rm (a)] $P\subseteq\la\bfe_1\ra$ and $\rho(P)$ is a subgroup of $\mathbb{F}_q^\times$;
\item[\rm (b)] $\DD:=(V, L^G)$ is a $2$-$(q^n, |P|+1, \l)$ design admitting $G$ as a group of automorphisms, where $L:=P\cup\{\bzero\}$, and $\Om:=(\bzero, L)^G$ is a 1-feasible $G$-orbit on the flag set of $\DD$; moreover, $\l=1$ if and only if $\rho(P)\cup \{0\}$ is a subfield of $\mathbb{F}_q$;
\item[\rm (c)] $\Ga(\DD, \Om, \Xi)$ is disconnected.
\end{itemize}
\end{lemma}

\proof
We first prove that $C_0 \nsubseteq P$. Suppose otherwise. Suppose that $P$ includes $j-1$ orbits of $H_{\bfe_1}$ of length $q^{n-1}$ $(1\leq j <q+1)$ and contains $\ell$ elements in $\la{\bfe_1}\ra$ $(1 \leq \ell <q)$. Then $|P| = jq^{n-1}+\ell-q=\gcd(jq^{n-1}+\ell-q,q^n-1)=\gcd(q^n-1, q^2-\ell q-j)$. If $q^2-\ell q-j \neq 0$, then $jq^{n-1}+\ell-q \leq q^2-\ell q-j$, which is impossible as $n\geq 4$. If $q^2-\ell q-j = 0$, then $j=q$, $\ell=q-1$, and thus $P = V \setminus \{{\bf 0}\}$, violating the condition $({u-1})/{|P|} \geq 2$. Therefore, $C_0 \nsubseteq P$.

Suppose that $P$ includes $j$ orbits of $H_{\bfe_1}$ of length $q^{n-1}$ $(0\leq j <q)$ and contains $\ell$ elements in $\la{\bfe_1}\ra$ $(1 \leq \ell <q)$. Then $|P| = jq^{n-1}+\ell=\gcd(jq^{n-1}+\ell,q^n-1)=\gcd(q^n-1, \ell q+j) \leq \ell q+j$. Since $n\geq 4$, it follows that  $j=0$ and $P\subseteq \la {\bfe_1} \ra$. Since $P$ is an imprimitive block of $G_{\la\bfe_1\ra, \bzero}$ on $\la\bfe_1\ra\setminus\{\bzero\}$, $\rho(P)$ is an imprimitive block of $\eta(G_{\la\bfe_1\ra, \bzero})$ on $\mathbb{F}_q^\times$ containing $1$.

Let $\bfe_1=\ba_1, \bb_1, \ba_2, \bb_2, \ldots, \ba_t, \bb_t$ be a symplectic basis of $(V, f)$, where $(\ba_i, \bb_i)$ is a hyperbolic pair, $i=1,2,\ldots, t$. With respect to this basis, for any $a\in \mathbb{F}_q^\times$, we have $\t(A(a), 0, \id)\in \Sp(n, q)\cap G_{\la\bfe_1\ra, \bzero}$, where $A(a):=\diag(a, 1/a, 1, 1, \ldots, 1, 1)$. Hence $\eta(G_{\la\bfe_1\ra, \bzero})$ contains $\GL(1, q)$ and by the discussion in Section \ref{subsec:AGammaL(1, q)}, $\rho(P)=\{a\in \mathbb{F}_q^\times\;|\;a\bfe_1\in P\}$ is a subgroup of the multiplicative group $\mathbb{F}_q^\times$.

Since $P$ and $\la\bfe_1\ra\setminus\{\bzero\}$ are both imprimitive blocks of $G_\bzero$ on $V\setminus\{\bzero\}$, we have $G_{\bzero, P}\leq G_{\bzero, \la\bfe_1\ra}$. Let $g$ be an element of $G_{\la\bfe_1\ra}$ interchanging $\bzero$ and $\bfe_1$. Then $J:=\la G_{\bzero, P}, g\ra\leq G_{\la\bfe_1\ra}\neq G$, and thus the $G$-flag graph $\Ga(\DD, \Om, \Xi)$ is disconnected by Lemma \ref{lem:connectedness of G-flag graph}.

The $\l$ of $\DD$ is $1$ or $|P|+1$ by Lemma \ref{lem:2 lambda}. Since $\eta(G_{\la\bfe_1\ra})$ is $2$-transitive on $\mathbb{F}_q$, the argument in Section \ref{subsec:AGammaL(1, q)} shows that $\l=1$ if and only if $\rho(P)\cup \{0\}$ is a subfield of $\mathbb{F}_q$.
\qed

\subsection{$G_{\bf 0} \unrhd \SL(2, q)=\Sp(2, q)$, $u =q^2 = p^d$}

This case contributes to the second possibility in (g) in Table \ref{tab0a}. Denote the underlying symplectic space by $(V, f)$, where $V=\mathbb{F}_q^2$ and $f$ is a symplectic form. Let $\bfe_1:=(1,0)$ and $H :=\Sp(2, q)=\SL(2, q)\unlhd G_{\bf 0}$. Define $C_\a$ as in Section \ref{subsec:G0,Sp(n,q),n>2}. Then $C_0 = \emptyset$ and $C_\a=\la{\bfe_1}\ra+{\bf z}_\a$ for $\a \in \mathbb{F}_q^\times$, where $\bz_\a\in C_\a$. Denote all $1$-subspaces of $V$ by $U=\la{\bfe_1}\ra$, $U_1$, $\ldots$, $U_q$.

\begin{lemma}
\label{lem:SL(2,q)=Sp(2,q)}
Suppose that $P$ is a nontrivial imprimitive block of $G_{\bf 0}$ on $V \setminus \{{\bf 0}\}$ containing $\bfe_1$. Then
\begin{itemize}
\item[\rm (a)] $P\subseteq\la\bfe_1\ra$ and $\{a\in \mathbb{F}_q^\times\;|\;a\bfe_1\in P\}$ is a subgroup of $\mathbb{F}_q^\times$;
\item[\rm (b)] $\DD:=(V, L^G)$ is a $2$-$(q^2, |P|+1, \l)$ design admitting $G$ as a group of automorphisms, where $L:=P\cup\{\bzero\}$, and $\Om:=(\bzero, L)^G$ is a 1-feasible $G$-orbit on the flag set of $\DD$; moreover, $\l=1$ if and only if $\{a\in \mathbb{F}_q^\times\;|\;a\bfe_1\in P\}\cup \{0\}$ is a subfield of $\mathbb{F}_q$;
\item[\rm (c)] $\Ga(\DD, \Om, \Xi)$ is disconnected.
\end{itemize}
\end{lemma}

\proof
Suppose that $P\nsubseteq\la {\bfe_1} \ra$. Write $P=(U+\bz_{\a_1})\cup\cdots\cup(U+\bz_{\a_t})\cup E\;(1\leq t < q)$, where $E$ is a subset of $\la {\bfe_1} \ra$ of size $\ell\;(1\leq \ell<q)$, $\a_1, \ldots, \a_t$ are distinct elements of $\mathbb{F}_q^\times$ and $\bz_{\a_j}\in C_{\a_j}$ for $1 \le j \le t$. Since $H$ is transitive on the set of $1$-subspaces of $V$, there exists $\g \in H$ such that $U^\g = U_1$. Hence $P^\g=(U_1+{\bf z}_{\a_1}^\g)\cup\cdots\cup(U_1+{\bf z}_{\a_t}^\g)\cup E^\g$. Since $U$ and $U_1$ are not parallel, $P^\g \cap P \neq \emptyset$ and thus $P=P^\g \supseteq U_1+{\bf z}_{\a_1}^\g$. Since $|(U_1+{\bf z}_{\a_1}^\g)\cap U|=1$ and $|(U_1+{\bf z}_{\a_1}^\g)\cap(U+{\bf z}_{\a_j})|=1$ for $1 \le j \le t$, we have $t+1\geq |U_1+{\bf z}_{\a_1}^\g| =q$ and thus $t=q-1$. Now $|P| = q^2-q+\ell=\gcd(q^2-q+\ell, q^2-1)=\gcd(q^2-1, q-\ell-1)$. Thus $\ell=q-1$ and $P=V \setminus \{{\bf 0}\}$, violating the condition $({u-1})/{|P|} \geq 2$. Therefore, $P\subseteq \la {\bfe_1} \ra$.

Similar to the treatment in Section \ref{subsec:G0,Sp(n,q),n>2}, one can show that the set $\{a\in \mathbb{F}_q^\times\;|\;a\bfe_1\in P\}$ is a subgroup of $\mathbb{F}_q^\times$. Moreover, the $G$-flag graph $\Ga(\DD, \Om, \Xi)$ is disconnected, and $\l=1$ if and only if $\{a\in \mathbb{F}_q^\times\;|\;a\bfe_1\in P\}\cup \{0\}$ is a subfield of $\mathbb{F}_q$.
\qed

\subsection{$G_{\bf 0} \unrhd \SL(n, q)$, $n \geq 3$, $u = q^n = p^d$}
\label{ssec:SLnq}

This case contributes to the second possibility in (g) in Table \ref{tab0a}. Let $P$ be a nontrivial imprimitive block of $G_{\bf 0}$ on $V \setminus \{{\bf 0}\}$ containing $\bfe_1:=(1,0,\ldots,0)$, where $V=\mathbb{F}_q^n$. Since $V \setminus \la {\bfe_1}\ra$ is a $G_{{\bf 0}, \bfe_1}$-orbit of length $q^n -q$, if $P$ includes this orbit, then $P = V \setminus \{{\bf 0}\}$ as $|P|$ is a divisor of $|V \setminus \{{\bf 0}\}|=q^n-1$, but this violates the condition $({u-1})/{|P|} \geq 2$. If $P$ does not include $V \setminus \la {\bfe_1}\ra$, then $P\subseteq \la \bfe_1 \ra$ and similar to the argument in Section \ref{subsec:G0,Sp(n,q),n>2} one can show that $\{a\in \mathbb{F}_q^\times\;|\;a\bfe_1\in P\}$ is a subgroup of the multiplicative group $\mathbb{F}_q^\times$. Set $L:=P\cup\{\bzero\}$, $\DD:=(V, L^G)$ and $\Om:=(\bzero, L)^G$. Then $\DD$ is a $2$-$(q^n, |P|+1, \l)$ design admitting $G$ as a group of automorphisms in its natural action, where $\l = 1$ or $|P|+1$, and $\Om$ is a 1-feasible $G$-orbit on the flag set of $\DD$. Moreover, the $G$-flag graph $\Ga(\DD, \Om, \Xi)$ is disconnected, and $\l=1$ if and only if $\{a\in \mathbb{F}_q^\times\;|\;a\bfe_1\in P\}\cup \{0\}$ is a subfield of $\mathbb{F}_q$.

\subsection{$G_{\bf 0} \unrhd G_2(q)$, $u = q^6 = p^d$, $q >2$, $q$ is even}
\label{subsec:G2(q)}

This case contributes to the third possibility in (g) in Table \ref{tab0a}. Suppose that $P$ is a nontrivial imprimitive block of $G_{\bf 0}$ on $V \setminus \{{\bf 0}\}$ containing $\bfe_1:=(1,0,0,0,0,0)$, where $V=\mathbb{F}_q^6$. Then $P$ is also an imprimitive block of $G_2(q)$ on $V \setminus \{{\bf 0}\}$ and $P$ is the union of some ${G_2(q)}_{\bfe_1}$-orbits on $V \setminus \{{\bf 0}\}$. We are going to determine all possible orbit-lengths of ${G_2(q)}_{\bfe_1}$ on $V \setminus \{{\bf 0}\}$ by using the knowledge of $G_2(q)$ from \cite[p.122, Section 4.3.4]{Wilson}.

Take a basis $\{\bx_1, \bx_2, \ldots, \bx_8\}$ of the octonion algebra $\mathbb{O}$ over $\mathbb{F}_q$ with multiplication given by Table \ref{tab:4 Multiplication table A}, or equivalently by Table \ref{tab:5 Multiplication table B}, where $\bfe:=\bx_4+\bx_5$ is the identity element of $\mathbb{O}$ (since the characteristic is $2$, we omit the signs).

\begin{table}[!ht]
\begin{minipage}{0.5\textwidth}
\centering
\scalebox{0.7}[0.8]{
\begin{tabular}{c|c|c|c|c|c|c|c|c}
${}$ & $\bx_1$ & $\bx_2$ & $\bx_3$ & $\bx_4$ & $\bx_5$ & $\bx_6$ & $\bx_7$ & $\bx_8$\\\hline
$\bx_1$ & ${\bf 0}$ & ${\bf 0}$ & ${\bf 0}$ & ${\bf 0}$ & $\bx_1$ & $\bx_2$ & $\bx_3$ & $\bx_4$\\\hline
$\bx_2$ & ${\bf 0}$ & ${\bf 0}$ & $\bx_1$ & $\bx_2$ & ${\bf 0}$ & ${\bf 0}$ & $\bx_5$ & $\bx_6$\\\hline
$\bx_3$ & ${\bf 0}$ & $\bx_1$ & ${\bf 0}$ & $\bx_3$ & ${\bf 0}$ & $\bx_5$ & ${\bf 0}$ & $\bx_7$\\\hline
$\bx_4$ & $\bx_1$ & ${\bf 0}$ & ${\bf 0}$ & $\bx_4$ & ${\bf 0}$ & $\bx_6$ & $\bx_7$ & ${\bf 0}$\\\hline
$\bx_5$ & ${\bf 0}$ & $\bx_2$ & $\bx_3$ & ${\bf 0}$ & $\bx_5$ & ${\bf 0}$ & ${\bf 0}$ & $\bx_8$\\\hline
$\bx_6$ & $\bx_2$ & ${\bf 0}$ & $\bx_4$ & ${\bf 0}$ & $\bx_6$ & ${\bf 0}$ & $\bx_8$ & ${\bf 0}$\\\hline
$\bx_7$ & $\bx_3$ & $\bx_4$ & ${\bf 0}$ & ${\bf 0}$ & $\bx_7$ & $\bx_8$ & ${\bf 0}$ & ${\bf 0}$\\\hline
$\bx_8$ & $\bx_5$ & $\bx_6$ & $\bx_7$ & $\bx_8$ & ${\bf 0}$ & ${\bf 0}$ & ${\bf 0}$ & ${\bf 0}$\\
\end{tabular}}
\caption{Multiplication table of $\mathbb{O}$
\label{tab:4 Multiplication table A}}
\end{minipage}
\begin{minipage}{0.5\textwidth}
\centering
\scalebox{0.6}[0.8]{
\begin{tabular}{c|c|c|c|c|c|c|c|c}
${}$    & $\bfe$ & $\bx_1$ & $\bx_8$ & $\bx_2$ & $\bx_7$ & $\bx_3$ & $\bx_6$ & $\bx_4$\\\hline
$\bfe$  & $\bfe$ & $\bx_1$ & $\bx_8$ & $\bx_2$ & $\bx_7$ & $\bx_3$ & $\bx_6$ & $\bx_4$\\\hline
$\bx_1$ & $\bx_1$ & ${\bf 0}$ & $\bx_4$ & ${\bf 0}$ & $\bx_3$ & ${\bf 0}$ & $\bx_2$ & ${\bf 0}$\\\hline
$\bx_8$ & $\bx_8$ & $\bfe+\bx_4$ & ${\bf 0}$ & $\bx_6$ & ${\bf 0}$ & $\bx_7$ & ${\bf 0}$ & $\bx_8$\\\hline
$\bx_2$ & $\bx_2$ & ${\bf 0}$ & $\bx_6$ & ${\bf 0}$ & $\bfe+\bx_4$ & $\bx_1$ & ${\bf 0}$ & $\bx_2$\\\hline
$\bx_7$ & $\bx_7$ & $\bx_3$ & ${\bf 0}$ & $\bx_4$ & ${\bf 0}$ & ${\bf 0}$ & $\bx_8$ & ${\bf 0}$\\\hline
$\bx_3$ & $\bx_3$ & ${\bf 0}$ & $\bx_7$ & $\bx_1$ & ${\bf 0}$ & ${\bf 0}$ & $\bfe+\bx_4$ & $\bx_3$\\\hline
$\bx_6$ & $\bx_6$ & $\bx_2$ & ${\bf 0}$ & ${\bf 0}$ & $\bx_8$ & $\bx_4$ & ${\bf 0}$ & ${\bf 0}$\\\hline
$\bx_4$ & $\bx_4$ & $\bx_1$ & ${\bf 0}$ & ${\bf 0}$& $\bx_7$ & ${\bf 0}$ & $\bx_6$ & $\bx_4$\\
\end{tabular}}
\caption{Multiplication table of $\mathbb{O}$
\label{tab:5 Multiplication table B}}
\end{minipage}
\end{table}
There is a quadratic form $N$ and an associated bilinear form $f$ satisfying
\begin{equation*}
\begin{aligned}
N(\bx_i)=0
\end{aligned}
\;\text{and}\;
\begin{gathered}
f(\bx_i,\bx_j)=
\begin{cases}
0, & \text{$i+j\neq9,$}\\
1, & \text{$i+j=9,$}
\end{cases}
\end{gathered}
\;\;\text{ $i, j=$ $1,2,\ldots,8$.}
\end{equation*}
The group $G_2(q)$ is the automorphism group of $\mathbb{O}$. Since $G_2(q)$ preserves the multiplication table of $\mathbb{O}$, one can verify that it preserves $N$ and $f$. Moreover, $G_2(q)$ induces a faithful action on $\bfe^\bot/{\la\bfe\ra}$, where
$\bfe^\bot=\la\bx_1, \bx_8,\bx_2,\bx_7,\bx_3,\bx_6,\bfe\ra$. Hence $G_2(q)$ can be embedded into $\Sp(6, q)$. There is also a symmetric trilinear form $t$ on $\bfe^\bot$, that is, $t(\bz_1,\bz_2,\bz_3)=t(\bz_{1^\s},\bz_{2^\s},\bz_{3^\s})$ for any $\s\in\Sym(\{1,2,3\})$ and $\bz_1,\bz_2,\bz_3\in \bfe^\bot$ (see \cite[p.420]{Aschbacher} for the definition of a symmetric trilinear form), defined by
\begin{equation}
\label{equ:trilinear form}
t(\bfe,\bx_1,\bx_8)=t(\bfe,\bx_2,\bx_7)=t(\bfe,\bx_3,\bx_6)=t(\bx_1,\bx_6,\bx_7)=t(\bx_2,\bx_3,\bx_8)=1
\end{equation}
\noindent and otherwise $t$ is zero on the basis vectors. It is straightforward to verify that
\begin{equation}
\label{equ:trilinear form equal f(xy,z)}
t(\bx,\by,\bz)=f(\bx\by,\bz), \text{\;for \;}\bx,\by,\bz \in \bfe^\bot.
\end{equation}
\noindent Since $G_2(q)$ preserves the multiplication and $f$, it also preserves $t$.

Denote by $\la\bx\ra$ the subspace of $\mathbb{O}$ spanned by $\bx$, and $\la\overline{\bx}\ra$ the subspace of $\bfe^\bot/{\la\bfe\ra}$ spanned by $\overline{\bx}$, where $\overline{\bx}=$ $\bx+\la\bfe\ra$. If $W$ is a subspace of $\mathbb{O}$, then we use $\overline{W}$ to denote $\la W,\bfe\ra/\la\bfe\ra$. The actions of $G_2(q)$ on $\bfe^\bot/{\la\bfe\ra}$ and $V$ are permutation isomorphic.

It is known that ${G_2(q)}_{\la\overline{\bx}_1\ra}$ has four orbits on the set of $1$-subspaces of $\bfe^\bot/{\la\bfe\ra}$ (see \cite[Lemma 3.1]{Cooperstein79} and \cite[p.72]{Kantor85}), which have lengths $1$, $q^5$, $q(q+1)$, $q^3(q+1)$ and are represented by $\la\overline{\bx}_1\ra$, $\la\overline{\bx}_8\ra$, $\la\overline{\bx}_2\ra$, $\la\overline{\bx}_7\ra$, respectively. (Since $\overline{\bx}_8$ is not perpendicular to $\overline{\bx}_1$ while $\overline{\bx}_2$ and $\overline{\bx}_7$ are perpendicular to $\overline{\bx}_1$, $\la\overline{\bx}_8\ra$ and $\la\overline{\bx}_2\ra$ are in distinct ${G_2(q)}_{\la\overline{\bx}_1\ra}$-orbits, and $\la\overline{\bx}_8\ra$ and $\la\overline{\bx}_7\ra$ are in distinct ${G_2(q)}_{\la\overline{\bx}_1\ra}$-orbits. Moreover, $\la\overline{\bx}_2\ra$ and $\la\overline{\bx}_7\ra$ are also in distinct ${G_2(q)}_{\la\overline{\bx}_1\ra}$-orbits. In fact, if this is not the case, say $\vp(\la\overline{\bx}_2\ra)=\la\overline{\bx}_7\ra$ for some $\vp \in {G_2(q)}_{\la\overline{\bx}_1\ra}$, then $\vp(\bx_1)=a\bx_1 + \ell\bfe$, $\vp(\bx_2)=c\bx_7 + s\bfe$ for some $a, c \neq 0$ and so ${\bf 0}=\vp(\bx_1)\vp(\bx_2)=(a\bx_1 + \ell\bfe)(c\bx_7 + s\bfe)=ac\bx_3+ \ell c\bx_7 + as\bx_1 + \ell s\bfe$, which is a contradiction as $ac\neq0$ and $\bx_1$, $\bx_3$, $\bx_7$, $\bfe$ are linearly independent.)

\begin{lemma}
\label{lem:orbit length, G2(q)}
(\cite[Lemma 4.8]{FXZ}) Let $\ba \in V\setminus \{\bzero\}$. Then ${G_2(q)}_\ba$ has $q-1$ orbits of length $1$, $q-1$ orbits of length $q^5$, one orbit of length $q(q^2-1)$ and one orbit of length $q^3(q^2-1)$ on $V\setminus \{\bzero\}$.
\end{lemma}

\begin{lemma}
\label{lem:G2q}
Let $G\leq \AGammaL(6, q)$ be 2-transitive on $V$ such that $G_{\bf 0} \unrhd G_2(q)$. Suppose that $P$ is a nontrivial imprimitive block of $G_{\bf 0}$ on $V \setminus \{{\bf 0}\}$ and let $\bfe_1 \in P$. Then the following hold:
\begin{itemize}
\item[\rm (a)] $P \subseteq \la\bfe_1\ra$ and $\{a\in \mathbb{F}_q^\times\;|\;a\bfe_1\in P\}$ is a subgroup of $\mathbb{F}_q^\times$;
\item[\rm (b)] setting $L:=P\cup\{\bzero\}$, $\DD:=(V, L^G)$ is a $2$-$(q^6, |P|+1, \l)$ design admitting $G$ as a group of automorphisms, where $\l=1$ or $|P|+1$, and $\Om:=(\bzero, L)^G$ is a 1-feasible $G$-orbit on the flag set of $\DD$; moreover, $\l=1$ if and only if $\{a\in \mathbb{F}_q^\times\;|\;a\bfe_1\in P\}\cup \{0\}$ is a subfield of $\mathbb{F}_q$;
\item[\rm (c)] $\Ga(\DD, \Om, \Xi)$ is disconnected.
\end{itemize}
\end{lemma}

\proof
Since $P$ is the union of some ${G_2(q)}_{\bfe_1}$-orbits on $V\setminus \{{\bf 0}\}$, we have a few combinations to consider. We prove that only one possibility can actually occur.

First, if $P$ includes both the orbit of length $q(q^2-1)$ and the orbit of length $q^3(q^2-1)$, then similar to the argument in Section \ref{subsec:G0,Sp(n,q),n>2} one can prove that $P=V \setminus \{{\bf 0}\}$, contradicting the assumption that $P$ is a nontrivial block of $V \setminus \{{\bf 0}\}$.

Next assume that $P$ includes the orbit of length $q^3(q^2-1)$, $i-1$ orbits of length $q^5$ for some $1\leq i <q+1$ and $\ell$ orbits of length $1$ for some $1\leq \ell <q$, but $P$ does not include the orbit of length $q(q^2-1)$. Then $|P| = iq^5-q^3+\ell=\gcd(iq^5-q^3+\ell, q^6-1)=\gcd(iq^5-q^3+\ell +\ell(q^6-1), q^6-1)=\gcd(q^3(\ell q^3+iq^2-1), q^6-1)=\gcd(\ell q^3+iq^2-1, q^6-1)$. Since $0<\ell q^3+iq^2-1$, we have $iq^5-q^3+\ell\leq \ell q^3+iq^2-1\leq q^4-1$, which is impossible.

Now assume that $P$ includes the orbit of length $q(q^2-1)$, $i$ orbits of length $q^5$ for some $0\leq i < q$ and $\ell$ orbits of length $1$ for some $1\leq \ell <q$, but $P$ does not include the orbit of length $q^3(q^2-1)$. Then
\begin{align}
\label{equ:gcd}
&|P| = iq^5+q^3-q+\ell=\gcd(iq^5+q^3-q+\ell, q^6-1)=\gcd(iq^5+q^3-q+\ell +\ell(q^6-1), q^6-1)\notag\\
&=\gcd(\ell q^5+iq^4+q^2-1, q^6-1)=\gcd(\ell q^5+iq^4+q^2-1, q^6-1-(\ell q^5+iq^4+q^2-1))\notag\\
&=\gcd(\ell q^5+iq^4+q^2-1, q^4-\ell q^3-iq^2-1).
\end{align}

\noindent Since $0<q^2-1\leq q^4-\ell q^3-iq^2-1\leq q^4-q^3-1$, we have $iq^5+q^3-q+\ell\leq q^4-\ell q^3-iq^2-1\leq q^4-q^3-1$, which implies $i=0$. Thus \eqref{equ:gcd} gives $|P| = q^3-q+\ell = \gcd(\ell q^5+q^2-1, q^4-\ell q^3-1)=\gcd(q^2(\ell q^3-q^2+\ell q+1), q^4-\ell q^3-1)=\gcd(\ell q^3-q^2+\ell q+1, q^4-\ell q^3-1)=\gcd(\ell q^3-q^2+\ell q+1, q^3-q+\ell) = \gcd(q^2-2\ell q+(\ell^2-1), q^3-q+\ell)$. Since $0\leq q^2-2\ell q+(\ell^2-1)=(\ell-q)^2-1\leq q^2-2q$, if $q^2-2\ell q+(\ell^2-1)\neq0$, then $q^3-q+\ell\leq q^2-2\ell q+(\ell^2-1)\leq q^2-2q$, which is impossible. Hence $q^2-2\ell q+(\ell^2-1)=0$, $\ell=q-1$ and $|P|=q^3-1$. We claim that this $P$ is not an imprimitive block of $G_{\bf 0}$ on $V \setminus \{{\bf 0}\}$. As in \cite[Section 1]{Aschbacher G2}, for $\bx\in \bfe^\bot$, define
\begin{align*}
\bx\mathit{\Delta}:=\{\by\in\bfe^\bot\;|\;t(\bx,\by,\bz)=0 \text{\;\;for any\;}\bz\in \bfe^\bot\}.
\end{align*}
Let $\vp \in {G_2(q)}_{\overline{\bx}_1}$ and suppose $\vp(\bx_1)=\bx_1+c\bfe$. Let $\by=a\bfe+a_1\bx_1+a_8\bx_8+a_2\bx_2+a_7\bx_7+a_3\bx_3+a_6\bx_6$ be a vector in $\bfe^\bot$. If $\by\in (\bx_1+c\bfe)\mathit{\Delta}$, then $t(\bx_1+c\bfe, \by, \bz) = 0$ for any $\bz\in \bfe^\bot$. When $\bz = \bfe,\bx_2,\bx_3$, we get $a_8=0$, $a_7=0$, $a_6=0$, and thus $\by \in W := \la \bfe,\bx_1,\bx_2,\bx_3\ra$. It follows that $(\bx_1+c\bfe)\mathit{\Delta}\subseteq W$. In addition, if $c=0$, then $\bx_1\mathit{\Delta}=\la \bx_1,\bx_2, \bx_3\ra$. Hence $\vp(\bx_2)\in (\bx_1+c\bfe)\mathit{\Delta}\subseteq W$, and $\vp(\overline{\bx}_2)\in \overline{W}$, which implies that the ${G_2(q)}_{\overline{\bx}_1}$-orbit containing $\overline{\bx}_2$ is included in $\overline{W}$. Since by Lemma \ref{lem:orbit length, G2(q)} this orbit has length $q(q^2-1)$, it must be $\overline{W}\setminus \la \overline{\bx}_1\ra$, and $P=\overline{W}\setminus\{\overline{\bzero}\}$, where $\overline{\bzero}$ denotes the zero vector in $\bfe^\bot/\la\bfe\ra$.

Suppose that $\overline{W}\setminus\{\overline{\bzero}\}$ is an imprimitive block of ${G_2(q)}$ on $\overline{\bfe^\bot} \setminus \{\overline{\bzero}\}$. Then there exists $\psi\in G_2(q)$ stabilizing $W$ such that $\psi(\overline{\bx}_1) =\overline{\bx}_2$. Let $\psi(\bx_1)=\bx_2+j\bfe$ for some $j\in \mathbb{F}_q$. Then $0 = t(\bx_1, \bw, \psi^{-1}(\bz))=t(\psi(\bx_1), \psi(\bw), \bz)=t(\bx_2+j\bfe, \psi(\bw), \bz)$ for any $\bw\in \bx_1\mathit{\Delta}$ and $\bz\in \bfe^\bot$. Hence $\psi(\bw)\in (\bx_2+j\bfe)\mathit{\Delta}$ and $\psi(\bx_1\mathit{\Delta})\subseteq (\bx_2+j\bfe)\mathit{\Delta}$. If $j=0$, then $\psi(\la\bx_1,\bx_2,\bx_3\ra)\subseteq W\cap \la\bx_1,\bx_2,\bx_6\ra=\la\bx_1,\bx_2\ra$ as $\bx_2\mathit{\Delta}=\la\bx_1,\bx_2,\bx_6\ra$ and $\psi$ stabilizes $W$, a contradiction. If $j\neq0$, then $(\bx_2+j\bfe)\mathit{\Delta}\subseteq\la\bfe, \bx_1,\bx_2,\bx_6\ra$ and $\psi(\la\bx_1,\bx_2,\bx_3\ra)\subseteq\la\bfe, \bx_1,\bx_2,\bx_6\ra\cap\la\bfe, \bx_1,\bx_2,\bx_3\ra=\la\bfe,\bx_1,\bx_2\ra$. Hence $\psi(\la\bx_1,\bx_2,\bx_3\ra)=\la\bfe,\bx_1,\bx_2\ra$ and $\bfe=\psi^{-1}(\bfe)\in \la\bx_1,\bx_2,\bx_3\ra$, a contradiction. Therefore, $\overline{W}\setminus\{\overline{\bzero}\}$ is not an imprimitive block of ${G_2(q)}$ on $\overline{\bfe^\bot} \setminus \{\overline{\bzero}\}$.

Therefore, the only possibility is that $P$ includes neither the orbit of length $q(q^2-1)$ nor the orbit of length $q^3(q^2-1)$. Similar to the treatment for $C_0 \nsubseteq P$ in Section \ref{subsec:G0,Sp(n,q),n>2}, one can show that $P \subseteq \la\bfe_1\ra$. Let $\bfe_1 = \bc_1, \bc_2, \ldots, \bc_6$ be a basis of $V$, and $\mu:\bfe^\bot/\la\bfe\ra\rightarrow V$ a linear map that maps $\overline{\bx}_1, \overline{\bx}_8, \overline{\bx}_2, \overline{\bx}_7, \overline{\bx}_3, \overline{\bx}_6$ to $\bc_1, \bc_2, \bc_3, \bc_4, \bc_5, \bc_6$, respectively. Up to permutation isomorphism, we only need to consider the action of $G_2(q)$ on $V$ defined by $\bc^g:=\mu((\mu^{-1}(\bc))^g)$ for $\bc\in V$ and $g\in G_2(q)$.

Now $G\leq \AGammaL(6, q)$. Let $\rho$ and $\eta$ be defined as in Section \ref{subsec:G0,Sp(n,q),n>2}. Similar to the proof of Lemma \ref{lem:orbit length, G2(q)}, one can verify that, for any $a\in \mathbb{F}_q^\times$, we have $\t(A(a), 0, \id)\in G_2(q)\cap G_{\la\bfe_1\ra, \bzero}$, where $A(a):=\diag(1/a, a, 1, 1, 1/a, a)$, with respect to the basis $\bc_1, \bc_2, \ldots, \bc_6$. Hence $\eta(G_{\la\bfe_1\ra, \bzero})$ contains $\GL(1, q)$ and the discussion in Section \ref{subsec:AGammaL(1, q)} shows that $\rho(P)=\{a\in \mathbb{F}_q^\times\;|\;a\bfe_1\in P\}$ is a subgroup of $\mathbb{F}_q^\times$. Since $G_{\bzero, P}\leq G_{\bzero, \la\bfe_1\ra}$, $J:=\la G_{\bzero, P}, g\ra\leq G_{\la\bfe_1\ra}\neq G$, where $g$ is an element in $G_{\la\bfe_1\ra}$ interchanging $\bzero$ and $\bfe_1$. Hence,  by Lemma \ref{lem:connectedness of G-flag graph}, the $G$-flag graph $\Ga(\DD, \Om, \Xi)$ is disconnected. The analysis in Section \ref{subsec:AGammaL(1, q)} shows that $\l=1$ if and only if $\{a\in \mathbb{F}_q^\times\;|\;a\bfe_1\in P\}\cup\{0\}$ is a subfield of $\mathbb{F}_q$.
\qed

\subsection{$G_{\bf 0} = G_2(2)' \cong \PSU(3, 3)$, $u = 2^6$}
\label{subsec:PSU(3,3)}

Suppose that $P$ is a nontrivial imprimitive block of $G_{\bf 0}$ on $V \setminus \{{\bf 0}\}$, where $V=\mathbb{F}_2^6$. Let ${\ba} \in P$. Since $1 < |P| < |V|$ and the orbit-lengths of $G_{\bzero, \ba}$ on $V\setminus\la\ba\ra$ are $2(2+1)$, $2^3(2+1)$, $2^4$ and $2^4$ (see \cite[p.72]{Kantor85}), $P$ is the union of $\{\ba\}$ and the $G_{\bzero, \ba}$-orbit of length $2(2+1)$. Similar to the proof of Lemma \ref{lem:G2q}, one can show that this $P$ is not an imprimitive block of ${G_2(2)'}$ on $V\setminus \{\bzero\}$.

\subsection{$G_{\bf 0} \cong A_6$ or $A_7$, $u = 2^4$}

Suppose that $P$ is a nontrivial imprimitive block of $G_{\bf 0}$ on $V \setminus \{{\bf 0}\}$, where $V=\mathbb{F}_2^4$. Let ${\ba} \in P$. Then $G_{\bzero, \ba}$ is transitive on $V\setminus \{\bzero, \ba\}$ when $G_{\bf 0} \cong A_7$, and $G_{\bzero, \ba}$ has orbit-lengths $6$ and $8$ on $V\setminus \{\bzero, \ba\}$ when $G_{\bf 0} \cong A_6$ (see \cite[p.72]{Kantor85}). Hence there is no $2$-design as in Lemma \ref{lem:1 basic} admitting $G$ as a group of automorphisms.

\subsection{The remaining cases}
\label{subsec:the remaining cases}

If $\DD$ is a $1$-design with point set $V$ admitting $G$ as a group of automorphisms and $\Om$ is a 1-feasible $G$-orbit on the set of flags of $\DD$, then the triple $(G,\DD,\Om)$ is said to be {\em proper} on $V$. Denote by $\Pi(V)$ the set of all proper triples on $V$.

The purpose of this subsection is to determine all proper triples in cases (vii)-(ix) in the opening paragraph of Section \ref{sec:affine case} up to flag-isomorphism to be defined as follows, and thus determine all possible flag graphs, contributing to (h), (i) and (j) in Table~\ref{tab0a}.

\begin{definition}
\label{def:G-flag-isomorphic}
{\em Let $(G_i,\DD_i,\Om_i)$ be a proper triple on $V_i$, $i=1,2$. If there exists a bijection $\rho:V_1\rightarrow V_2$ such that the action of $G_1$ on $V_1$ is permutation isomorphic to the action of $G_2$ on $V_2$ with respect to $\rho$, and the action of $G_1$ on $\Om_1$ is permutation isomorphic to the action of $G_2$ on $\Om_2$ with respect to the bijection from $\Om_1$ to $\Om_2$ induced by $\rho$, then $(G_1,\DD_1,\Om_1)$ and $(G_2,\DD_2,\Om_2)$ are said to be \emph{flag-isomorphic} with respect to the \emph{flag-isomorphism} $\rho$.

Flag-isomorphism is an equivalence relation on any subfamily $\Pi$ of $\Pi(V)$. A subset of $\Pi$ that has exactly one proper triple from each equivalence class (of this equivalence relation) is called a {\em representative subset} of $\Pi$.
}
\end{definition}

In what follows we set $V:=\mathbb{F}_p^d$, where $d\geq2$ and $p$ is a prime, and denote
$$
\Pi(d,p):= \left\{(G,\DD,\Om)\in \Pi(\mathbb{F}_p^d)\;|\;G\leq \AGL(d,p)\text{\;is $2$-transitive on\;} \mathbb{F}_p^d\right\}.
$$
Recall that $T$ denotes the group of translations of $V$.
Let $(G,\DD,\Om),(\widetilde{G},\widetilde{\DD},\widetilde{\Om})\in \Pi(d,p)$ be flag-isomorphic with respect to some $\rho\in \Sym(V)$. Then $\rho^{-1}G\rho=\widetilde{G}$, $\rho^{-1}T\rho=T$ and $\rho^{-1}G_{\bzero}\rho=\widetilde{G}_{\bzero}$. By \cite[p.110, Corollary 4.2B]{Dixon-Mortimer}, $\rho\in N_{\Sym(V)}(T)=\AGL(d,p)$. Since $\widetilde{G}_{\bzero}=\rho^{-1}G_{\bzero}\rho=(G^\rho)_{\bzero^\rho}=\widetilde{G}_{\bzero^\rho}$ and $\widetilde{G}_{\bzero}$ is transitive on $V\setminus\{\bzero\}$, we have $\bzero^\rho=\bzero$ and thus $\rho\in \GL(d,p)$. Hence $G_\bzero$ and $\widetilde{G}_\bzero$ are conjugate in $\GL(d,p)$, and up to flag-isomorphism it suffices to consider one representative in the conjugacy class of subgroups of $\GL(d,p)$ containing $G_\bzero$. We may thus assume $G=\widetilde{G}$ in the following. Then $\rho \in N_{\GL(d,p)}(G_\bzero)$. Since $\rho$ is a flag-isomorphism, we have $\Om^\rho=\widetilde{\Om}$, or equivalently $(\bzero,L)^\rho=(\bzero,\widetilde{L})$ for some $(\bzero,{L})\in \Om$ and $(\bzero,\widetilde{L})\in \widetilde{\Om}$ as $\rho$ normalizes $G_\bzero$. Denote $H:=G_{\bzero,L}$ and $\widetilde{H}:={G}_{\bzero,\widetilde{L}}$. Then $H^\rho=\widetilde{H}$ and so $H$ and $\widetilde{H}$ are conjugate in $N_{\GL(d,p)}(G_\bzero)$. Up to flag-isomorphism it suffices to consider one representative in the conjugacy class of subgroups of $N_{\GL(d,p)}(G_\bzero)$ containing $H$. We may thus assume $H=\widetilde{H}$ in the following.

\begin{lemma}
\label{lem:same system of blocks}
Let $I$ be a finite group acting transitively on a set $\Delta$. Suppose that $H \le J \leq I$ such that $J$ is transitive on $\Delta$, $J_x\leq H$ and $J_y\leq H$ for two points $x,y$ in $\Delta$. Set $P:=x^H$ and $Q:=y^H$. Then there exists $\s\in Y:=N_{I}(J)$ such that $(P^J)^\s=Q^J$ if and only if $x$ and $y$ are in the same $N_{Y}(H)$-orbit on $\Delta$. In particular, two systems of blocks $P^J$ and $Q^J$ are identical if and only if $x$ and $y$ are in the same $N_J(H)$-orbit on $\Delta$.
\end{lemma}

\proof
Assume that $(P^J)^\s=Q^J$ for some $\s\in N_I(J)$. Then $y^{Hg}=x^{H\s}$ and $y=x^{\t\xi}$ for some $\t\in H$ and $g\in J$, where $\xi:=\s g^{-1}$. We need to show that $\xi\in N_Y(H)$. Set $z:=x^\t=y^{\xi^{-1}}$. Then $z^H = x^H= y^{Hg\s^{-1}}= z^{\xi H \xi^{-1}}$, and we set $M= z^H$. Since $H_z=(H_x)^\t=(J_x)^\t=J_z=(J_y)^{\xi^{-1}}=(H_y)^{\xi^{-1}}=(\xi H\xi^{-1})_z$, we have $H=J_z H=J_M=J_z H^{\xi^{-1}}=H^{\xi^{-1}}$.

Conversely, if $y=x^\s$ for some $\s\in N_Y(H)$, then $Q=y^H=x^{H\s}=P^\s$ and so $(P^J)^\s=Q^J$.
\qed

\smallskip
By Lemma \ref{lem:same system of blocks}, $L\setminus\{\bzero\}=x^H$ and $\widetilde{L}\setminus\{\bzero\}=y^H$ for some $x,y$ that are in the same $N_{Y}(H)$-orbit on $V$, where $Y=N_{\GL(d,p)}(G_\bzero)$. Therefore, we obtain the following result.

\begin{lemma}
\label{lem:maximal representative subset}
Let $V:=\mathbb{F}_p^d$, where $d\geq2$ and $p$ is a prime. Suppose that $\GL(d,p)\leq \Sym(V\setminus\{\bzero\})$ has exactly $\ell$ conjugacy classes $\CC_1, \ldots, \CC_\ell$ of transitive subgroups. For $1 \le i \le \ell$, let $G_i\in \CC_i$ and $N_i:=N_{\GL(d,p)}(G_i)$, and suppose that $N_i$ has exactly $n(i)$ conjugacy classes $\HH_{i,1}, \ldots, \HH_{i,n(i)}$ of subgroups of $G_i$ containing some point stabilizer of ${G_i}$. For $1\leq i\leq \ell$ and $1 \leq j \leq n(i)$, let $H_{i,j}\in \HH_{i,j}$, $M_{i,j}:=N_{N_i}(H_{i,j})$, $\Delta_{i,j} := \{\bx\in V\setminus\{\bzero\}\;|\;({G_i})_{\bx}\leq H_{i,j}\}$, and suppose that the $M_{i,j}$-orbits on $\Delta_{i,j}$ are $\bx_{i,j,k}^{M_{i,j}}$, $1 \le k \le s(i,j)$, where $\bx_{i,j,k} \in \Delta_{i,j}$ for each $k$. For $1\leq i\leq \ell$, $1 \leq j \leq n(i)$ and $1 \le k \le s(i,j)$, set $\DD_{i,j,k} := (V,L_{i,j,k}^{T \rtimes G_i})$ and $\Om_{i,j,k} := (\bzero, L_{i,j,k})^{T \rtimes G_i}$, where $L_{i,j,k} := P_{i,j,k} \cup\{\bzero\}$ with $P_{i,j,k}:=\bx_{i,j,k}^{H_{i,j}}$. Then
$$
\{(T \rtimes G_i,\DD_{i,j,k},\Om_{i,j,k})\;|\;1\leq i\leq \ell, 1 \leq j \leq n(i), 1 \le k \le s(i,j)\}
$$
is a representative subset of $\Pi(d,p)$.
\end{lemma}

Based on Lemma \ref{lem:maximal representative subset} and by using \magma, we have computed a representative subset $S$ of each of the following sets:
\begin{align*}
\Pi_{6,3}&:=\{(G,\DD,\Om)\in \Pi(6,3)\;|\;G_\bzero\cong \SL(2,13)\},\\
\Pi_{4,3}^1&:=\{(G,\DD,\Om)\in \Pi(4,3)\;|\;G_\bzero\unrhd\SL(2, 5)\},\\
\Pi_{4,3}^2&:=\{(G,\DD,\Om)\in \Pi(4,3)\;|\;G_\bzero\unrhd E\},\\
\Pi_{2,p}&:=\{(G,\DD,\Om)\in \Pi(2,p)\;|\;G_\bzero\unrhd\SL(2, 5)\text{\;or\;}G_{\bf 0}\unrhd\SL(2, 3)\},
\end{align*}
\noindent where $E$ is an extraspecial group of order $32$ and $p=5,7,11,19,23,29,59$. Note that $\DD$ is a $2$-$(u, r+1, \l)$ design in each case, where $u = 3^6$ for $\Pi_{6,3}$, $u = 3^4$ for $\Pi_{4,3}^1$ and $\Pi_{4,3}^2$, and $u = p^2$ for $\Pi_{2,p}$. Our computational results are summarized in Tables \ref{tab1}-\ref{tab:p=59}, where the first row in each table gives the common value of $|G_\bzero|$ for those proper triples in $S$ whose group entries are pairwise conjugate in $\AGL(d,p)$, and the last row gives the number of proper triples in $S$ whose $r$ and $\l$ values are the same and the group entries are pairwise conjugate in $\AGL(d,p)$.

In the case where $G_\bzero \cong \SL(2,13)$ and $u = 3^6$, our computational results are given in Table~\ref{tab1}, and we find that the corresponding $G$-flag graphs are all disconnected. In addition, by \cite[p.73]{Kantor85} we know that $\l=1$ if and only if one of the following occurs: $r=2$ and $\DD=\AG(6,3)$; $r=8$ and $\DD$ is one of the two Hering designs \cite{Hering}; $r=26$ and $\DD$ is the Hering affine plane of order $27$ (see \cite{Hering27} and \cite[p.236]{Dembowski}).

\begin{table}[!ht]
\caption{Representatives of $\Pi_{6,3}$}
\label{tab1}
\centering
\begin{tabular}{|c|c|c|c|c|c|c|c|c|}\hline
$|G_\bzero|$ & \multicolumn{8}{|c|}{$2184$}\\\hline
$r+1$ & $3$ & $5$ & $9$ & $9$ & $14$ & $27$ & $27$ & $53$\\\hline
$\l$       & $1$ & $5$ & $9$ & $1$ & $14$ & $27$ & $1$  & $53$\\\hline
No.        & $1$ & $3$ & $3$ & $2$ & $2$  & $1$  & $1$  & $2$\\
\hline
\end{tabular}
\end{table}

In the case where $G_{\bf 0} \unrhd \SL(2, 5)$, $d=4$ and $p=3$, the results are given in Table~\ref{tab2} and the corresponding $G$-flag graphs are all disconnected. In addition, by \cite{Kantor85}, $\l=1$ if and only if one of the following occurs: $r=8$ and $\DD=\AG(2,9)$; $r=8$ and $\DD$ is the exceptional nearfield plane (see \cite[pp.214, 232, 236]{Dembowski}); $r=2$ and $\DD=\AG(4,3)$.

\begin{table}[!ht]
\caption{Representatives of $\Pi_{4,3}^1$}
\label{tab2}
\centering
\begin{tabular}{*{18}{|c}|}\hline
$|G_\bzero|$ & \multicolumn{5}{|c|}{$240$} & \multicolumn{4}{|c|}{$480$} & \multicolumn{4}{|c|}{$480$} &\multicolumn{4}{|c|}{$960$}\\\hline
$r+1$ & $3$ & $5$ & $9$ &  $17$ & $41$ & $3$ & $5$ & $9$ & $41$ & $3$ & $5$ & $9$ & $41$ & $3$ & $5$ & $9$ & $41$\\\hline
$\l$       & $1$ & $5$ & $1$ &  $17$ & $41$ & $1$ & $5$ & $1$ & $41$ & $1$ & $5$ & $1$ & $41$ & $1$ & $5$ & $1$ & $41$\\\hline
No.        & $1$ & $3$ & $2$ &  $1$  & $1$  & $1$ & $2$ & $1$ & $1$  & $1$ & $1$ & $1$ & $1$  & $1$ & $1$ & $1$ & $1$\\
\hline
\end{tabular}
\end{table}

In the case where $G_\bzero \unrhd E$,  $d=4$ and $p=3$, where $E$ is an extraspecial group of order $32$, the results are give in Table~\ref{tab3} and the corresponding $G$-flag graphs are all disconnected. By \cite{Kantor85}, $\l=1$ if and only if one of the following occurs: $r=8$ and $\DD$ is the exceptional nearfield plane (see \cite[pp.214, 232, 236]{Dembowski}); $r=2$ and $\DD=\AG(4,3)$.

\begin{table}[!ht]
\caption{Representatives of $\Pi_{4,3}^2$}
\label{tab3}
\centering
\begin{tabular}{*{21}{|c}|}\hline
$|G_\bzero|$ & \multicolumn{5}{|c|}{$160$} & \multicolumn{5}{|c|}{$320$} & \multicolumn{4}{|c|}{$640$} & \multicolumn{3}{|c|}{$1920$} & \multicolumn{3}{|c|}{$3840$}\\\hline
$r+1$  & $3$ & $5$ & $9$ & $9$ & $17$ & $3$ & $5$ & $9$ & $9$ & $17$ & $3$ & $5$ & $9$ & $17$ & $3$ & $9$ & $17$ & $3$ & $9$ & $17$\\\hline
$\l$   & $1$ & $5$ & $9$ & $1$ & $17$ & $1$ & $5$ & $9$ & $1$ & $17$ & $1$ & $5$ & $1$ & $17$ & $1$ & $1$ & $17$ & $1$ & $1$ & $17$\\\hline
No.    & $1$ & $3$ & $2$ & $1$ & $1$  & $1$ & $2$ & $1$ & $1$ & $1$  & $1$ & $1$ & $1$ & $1$  & $1$ & $1$ & $1$  & $1$ & $1$ & $1$\\
\hline
\end{tabular}
\end{table}

Finally, in the case where $G_{\bf 0} \unrhd \SL(2, 5)$ or $G_{\bf 0} \unrhd \SL(2, 3)$, $d=2$ and $p = 5, 7, 11, 19, 23, 29$ or $59$, the results are given in Tables \ref{tab:p=5}-\ref{tab:p=59}, respectively, and the corresponding $G$-flag graphs are all disconnected. By \cite{Kantor85}, $\l=1$ if and only if $\DD=\AG(2,p)$.

\begin{table}[!ht]
\caption{Representatives of $\Pi_{2,5}$}
\label{tab:p=5}
\centering
\begin{tabular}{*{20}{|c}|}\hline
$|G_\bzero|$ & \multicolumn{6}{|c|}{$24$} & \multicolumn{4}{|c|}{$48$} & \multicolumn{3}{|c|}{$96$} & \multicolumn{2}{|c|}{$120$} & \multicolumn{2}{|c|}{$240$} & \multicolumn{2}{|c|}{$480$}\\\hline
$r+1$ & $3$ & $4$ & $5$ & $5$ & $7$ & $9$ & $3$ & $5$ & $5$ & $9$ & $3$ & $5$ & $9$ & $3$ & $5$ & $3$ & $5$ & $3$ & $5$ \\\hline
$\l$  & $3$ & $4$ & $5$ & $1$ & $7$ & $9$ & $3$ & $1$ & $5$ & $9$ & $3$ & $1$ & $9$ & $3$ & $1$ & $3$ & $1$ & $3$ & $1$ \\\hline
No.   & $1$ & $1$ & $1$ & $1$ & $1$ & $1$ & $1$ & $1$ & $1$ & $1$ & $1$ & $1$ & $1$ & $1$ & $1$ & $1$ & $1$ & $1$ & $1$ \\\hline
\end{tabular}
\end{table}

\begin{table}[!ht]
\caption{Representatives of $\Pi_{2,7}$}\label{tab:p=7}
\centering
\begin{tabular}{*{16}{|c}|}\hline
$|G_\bzero|$ & \multicolumn{9}{|c|}{$48$} & \multicolumn{6}{|c|}{$144$}\\\hline
$r+1$ & $3$ & $4$ & $5$ & $7$ & $7$ & $9$ & $13$ & $17$ & $25$ & $3$ & $4$ & $7$ & $9$ & $13$ & $25$ \\\hline
$\l$  & $3$ & $4$ & $5$ & $7$ & $1$ & $9$ & $13$ & $17$ & $25$ & $3$ & $4$ & $1$ & $9$ & $13$ & $25$ \\\hline
No.   & $1$ & $2$ & $3$ & $1$ & $1$ & $3$ & $2$  & $1$  & $1$  & $1$ & $1$ & $1$ & $1$ & $1$  & $1$  \\\hline
\end{tabular}
\end{table}

\begin{table}[!ht]
\caption{Representatives of $\Pi_{2,11}$}\label{tab:p=11}
\centering
\scalebox{0.8}[0.9]{
\begin{tabular}{*{24}{|c}|}\hline
$|G_\bzero|$ & \multicolumn{11}{|c|}{$120$} & \multicolumn{12}{|c|}{$120$}\\\hline
$r+1$ & $3$ & $4$ & $5$ & $6$ & $7$ & $9$ & $11$ & $11$ & $13$ & $21$ & $25$ & $3$ & $4$ & $5$ & $6$ & $7$ & $9$ & $11$ & $16$ & $21$ & $25$ & $31$ & $41$\\\hline
$\l$  & $3$ & $4$ & $5$ & $6$ & $7$ & $9$ & $11$ & $1$  & $13$ & $21$ & $25$ & $3$ & $4$ & $5$ & $6$ & $7$ & $9$ & $1$  & $16$ & $21$ & $25$ & $31$ & $41$\\\hline
No.   & $1$ & $2$ & $3$ & $2$ & $2$ & $1$ & $1$  & $1$  & $2$  & $2$  & $1$  & $1$ & $3$ & $2$ & $1$ & $3$ & $1$ & $1$  & $3$  & $2$  & $1$  & $3$ & $1$\\\hline
\hline
$|G_\bzero|$ & \multicolumn{12}{|c|}{$240$} & \multicolumn{4}{|c|}{$600$} & \multirow{4}*{}  & \multirow{4}*{}  & \multirow{4}*{}  & \multirow{4}*{}  & \multirow{4}*{}  & \multirow{4}*{}  & \multirow{4}*{} \\\cline{1-17}
$r+1$ & $3$ & $4$ & $5$ & $6$ & $7$ & $9$ & $11$ & $16$ & $21$ & $25$ & $31$ & $41$ & $3$ & $6$ & $11$ & $21$ & & & & & & & \\\cline{1-17}
$\l$  & $3$ & $4$ & $5$ & $6$ & $7$ & $9$ & $1$  & $16$ & $21$ & $25$ & $31$ & $41$ & $3$ & $6$ & $1$  & $21$ & & & & & & & \\\cline{1-17}
No.   & $1$ & $2$ & $1$ & $1$ & $2$ & $1$ & $1$  & $2$  & $1$  & $1$  & $2$  & $1$  & $1$ & $1$ & $1$  & $1$  & & & & & & & \\\hline
\end{tabular}}
\end{table}

\begin{table}[!ht]
\caption{Representatives of $\Pi_{2,19}$}\label{tab:p=19}
\centering
\begin{tabular}{*{12}{|c}|}\hline
$|G_\bzero|$ & \multicolumn{11}{|c|}{$1080$}\\\hline
$r+1$ & $3$ & $4$ & $7$ & $9$ & $10$ & $13$ & $19$ & $25$ & $37$ & $73$ & $121$ \\\hline
$\l$  & $3$ & $4$ & $7$ & $9$ & $10$ & $13$ & $1$  & $25$ & $37$ & $73$ & $121$ \\\hline
No.   & $1$ & $1$ & $1$ & $2$ & $1$  & $1$  & $1$  & $2$  & $1$  & $2$  & $1$   \\\hline
\end{tabular}
\end{table}

\begin{table}[!ht]
\caption{Representatives of $\Pi_{2,23}$}\label{tab:p=23}
\centering
\begin{tabular}{*{19}{|c}|}\hline
$|G_\bzero|$ & \multicolumn{18}{|c|}{$528$}\\\hline
$r+1$ & $3$ & $4$ & $5$ & $7$ & $9$ & $12$ & $13$ & $17$ & $23$ & $25$ & $34$ & $45$ & $49$ & $67$ & $89$ & $133$ & $177$ & $265$ \\\hline
$\l$  & $3$ & $4$ & $5$ & $7$ & $9$ & $12$ & $13$ & $17$ & $1$  & $25$ & $34$ & $45$ & $49$ & $67$ & $89$ & $133$ & $177$ & $265$ \\\hline
No.   & $1$ & $4$ & $9$ & $4$ & $7$ & $1$  & $4$  & $3$  & $1$  & $1$  & $4$  & $9$  & $1$  & $4$  & $7$  & $4$   & $3$   &   $1$ \\\hline
\end{tabular}
\end{table}

\begin{table}[!ht]
\caption{Representatives of $\Pi_{2,29}$}\label{tab:p=29}
\centering
\scalebox{0.8}[0.9]{
\begin{tabular}{*{24}{|c}|}\hline
$|G_\bzero|$ & \multicolumn{23}{|c|}{$840$}\\\hline
$r+1$ & $3$ & $4$ & $5$ & $6$ & $7$ & $8$ & $9$ & $11$ & $13$ & $15$ & $21$ & $22$ & $25$ & $29$ & $29$ & $36$ & $43$ & $57$ & $71$ & $85$ & $121$ & $141$ & $169$\\\hline
$\l$  & $3$ & $4$ & $5$ & $6$ & $7$ & $8$ & $9$ & $11$ & $13$ & $15$ & $21$ & $22$ & $25$ & $29$ & $1$  & $36$ & $43$ & $57$ & $71$ & $85$ & $121$ & $141$ & $169$ \\\hline
No.   & $1$ & $6$ & $9$ & $4$ & $6$ & $1$ & $3$ & $4$  & $6$  & $1$  & $4$  & $6$  & $3$  & $8$  & $1$  & $4$  & $6$  & $3$  & $4$  & $6$  & $1$   & $4$   & $3$  \\\hline
\hline
$|G_\bzero|$ & \multicolumn{23}{|c|}{$1680$}\\\hline
$r+1$ & $3$ & $4$ & $5$ & $6$ & $7$ & $8$ & $9$ & $11$ & $13$ & $15$ & $21$ & $22$ & $25$ & $29$ & $29$ & $36$ & $43$ & $57$ & $71$ & $85$ & $121$ & $141$ & $169$\\\hline
$\l$  & $3$ & $4$ & $5$ & $6$ & $7$ & $8$ & $9$ & $11$ & $13$ & $15$ & $21$ & $22$ & $25$ & $29$ & $1$  & $36$ & $43$ & $57$ & $71$ & $85$ & $121$ & $141$ & $169$ \\\hline
No.   & $1$ & $2$ & $3$ & $2$ & $2$ & $1$ & $1$ & $2$  & $2$  & $1$  & $2$  & $2$  & $1$  & $2$  & $1$  & $2$  & $2$  & $1$  & $2$  & $2$  & $1$   & $2$   & $1$  \\\hline
\end{tabular}}
\end{table}

\begin{table}[!ht]
\caption{Representatives of $\Pi_{2,59}$}\label{tab:p=59}
\centering
\scalebox{0.75}[0.9]{
\begin{tabular}{*{23}{|c}|}\hline
$|G_\bzero|$ & \multicolumn{22}{|c|}{$3480$}\\\hline
$r+1$ & $3$ & $4$  & $5$  & $6$ & $7$  & $9$ & $11$ & $13$ & $21$ & $25$ & $30$ & $59$ & $88$ & $117$ & $121$ & $146$ & $175$ & $233$ & $291$ & $349$ & $581$ & $697$ \\\hline
$\l$  & $3$ & $4$  & $5$  & $6$ & $7$  & $9$ & $11$ & $13$ & $21$ & $25$ & $30$ & $1$  & $88$ & $117$ & $121$ & $146$ & $175$ & $233$ & $291$ & $349$ & $581$ & $697$\\\hline
No.   & $1$ & $10$ & $15$ & $6$ & $10$ & $5$ & $6$  & $10$ & $6$  & $5$  & $1$  & $1$  & $10$ & $15$  & $1$   & $6$   & $10$  & $5$   & $6$   & $10$  & $6$   & $5$ \\\hline
\end{tabular}}
\end{table}

\medskip
We have completed the proof of Theorem \ref{THM:MAIN RESULT}.

\medskip
\medskip
\noindent\textbf{Acknowledgements}~ T. Fang was supported by the China Scholarship Council, X.~G.~Fang by National Science Foundation of China (NSFC 11231008), and S. Zhou by the Australian Research Council (FT110100629) and the MRGSS of the University of Melbourne. The authors would like to thank the anonymous referees for their comments that lead to improvements of presentation.

\newpage
\appendix

\section*{Appendix: Sample \magma codes}
The following \magma codes were used in the case where $G_{\bf 0} \unrhd \SL(2, 5)$ or $G_{\bf 0} \unrhd \SL(2, 3)$, $d=2$ and $p=5,7,11,19,23,29$ or $59$ in Section \ref{subsec:the remaining cases}. For other cases in Section \ref{subsec:the remaining cases}, the \magma codes are similar, and we just give the {\em filtering conditions} by which we get the groups in Section \ref{subsec:the remaining cases} (see ``$/*\ldots*/$'' parts in the following \magma codes).

\medskip
\noindent$/*$ Searching for a representative subset of $\Pi_{2, p}\;*/$

\begin{verbatim}
d:=2;
for p in [5,7,11,19,23,29,59] do

u:=p^d; A,V:=AGL(d,p); C:=GSet(A); g:=A.1;
T:=NormalSubgroups(A:OrderEqual:=u); T:=T[1]`subgroup;
X:=Stabilizer(A,1);

/* The following three lines of codes are called the filtering
   conditions */
L:=Subgroups(X:OrderMultipleOf:=p^2-1);
L:=[a`subgroup:a in L|#Orbits(a`subgroup) eq 2];
L1:=[a:a in L|#[b:b in NormalSubgroups(a:OrderEqual:=120)|IsIsomorphic
    (b`subgroup,SL(2,5)) eq  true]+#[b:b in NormalSubgroups
    (a:OrderEqual:=24)|IsIsomorphic(b`subgroup,SL(2,3)) eq  true] gt 0];

v:=C!1; i:=1; t:=#L1;

while i lt t+1 do
  "";
  "Transitive subgroup of GL(",d,",",p,") containing a normal subgroup
  isomorphic to SL(2,5) or SL(2,3), No.", i, " of order ", #L1[i];
  "";
  G0:=L1[i];
  Y:=Normalizer(X,G0);
  G:=sub<A|T,G0>; f:=#G;
  b:=Round(#G0/(u-1));
  S:=Subgroups(G0:OrderMultipleOf:=b);
  S:=[a`subgroup:a in S|#(a`subgroup) ne b];
  s:=#S;
  S1:=[S[1]]; s1:=2;
  while s1 lt s+1 do m2:=#S1; a1:=1;
    while a1 lt m2+1 do
       if IsConjugate(Y,S1[a1],S[s1]) eq true then break;
          else a1:=a1+1;
       end if;
    end while;
    if a1 gt m2 then S1:=Append(S1,S[s1]);
    end if;
    s1:=s1+1;
  end while;

time  for y in G do
    if [1,2]^y eq [2,1] then g:=y;
       "Find an element in G interchanging",V[1]," and ",V[2];
       "";
       break y;
    end if;
  end for;
  z:=g;

  j:=1; s:=#S1;
  while j lt s do
    "";
    "subgroup of G0 No.", j;
    "";
    H:=S1[j]; h:=#H; m:=1; e:=Round(h/b+1);
    NH:=Normalizer(Y,H); Ob:=Orbits(NH); n:=#Ob;
time    while m lt n+1 do C:=Ob[m];
      for x in C do
        if #Stabilizer(H,x) eq b then  v:=x;
             "the block of the 2-design has size", e;
                 for y in G0 do
                    if 2^y eq v then g:=y^(-1)*z*y;
                       break x;
                    end if;
                 end for;
           else break x;
        end if;
      end for;

      if v eq 1 then m:=m+1;
        else
           M:=sub<A|H,g>; c:=#M;
           if c eq h*e then "lambda equals 1";
                              "";
              else "lambda equals", e;
                  if c eq f then "G-flag graph is connected";
                  end if;
                   "";
           end if;
           m:=m+1; v:=1;
      end if;
    end while;
    j:=j+1;
  end while;

  i:=i+1;
end while;

end for;
\end{verbatim}

The filtering conditions in searching for a representative subset of $\Pi_{6,3}$ are as follows:
\begin{verbatim}
L:=Subgroups(X:OrderEqual:=2184);
L:=[a`subgroup:a in L|#Orbits(a`subgroup) eq 2];
L1:=[a:a in L|IsIsomorphic(a,SL(2,13)) eq  true];
\end{verbatim}

The filtering conditions in searching for a representative subset of $\Pi_{4,3}^1$ are as follows:
\begin{verbatim}
L:=Subgroups(X:OrderMultipleOf:=80);
L:=[a`subgroup:a in L|#Orbits(a`subgroup) eq 2];
L1:=[a:a in L|#[b:b in NormalSubgroups(a:OrderEqual:=120)|
    IsIsomorphic(b`subgroup,SL(2,5)) eq  true] gt 0];
\end{verbatim}

The filtering conditions in searching for a representative subset of $\Pi_{4,3}^2$ are as follows:
\begin{verbatim}
L:=Subgroups(X:OrderMultipleOf:=80);
L:=[a`subgroup:a in L|#Orbits(a`subgroup) eq 2];
L1:=[a:a in L|#pCore(a,2) eq 32];
\end{verbatim}

\end{document}